\documentclass{llncs}

\usepackage[mathscr]{eucal}      
\usepackage{amssymb,amsfonts,amsmath}             
\usepackage{latexsym}      
\usepackage{amscd}      
\usepackage{xy}      

\usepackage{enumerate}      
    
\input xymatrix      
\input xyarrow

\newdimen\xxxup
\newdimen\yyyup

\def\arrup{\lline\char'066}
\def\arrdn{\lline\char'077}
\def\arrw{\lline\char'033}
\def\arre{\lline\char'055}

\def\SA #1 #2 {\xup=#1\unit \yup=#2\unit
    \killglue\rlap{\kern \xup
    \raise \yup\hbox{\arrdn}}\ignorespaces}

\def\NA #1 #2 {\xup=#1\unit \yup=#2\unit
   \advance \yup by -0.9\unit
      \killglue\rlap{\kern \xup\raise\yup\hbox{\arrup}}\ignorespaces}

\def\WA #1 #2 {\xup=#1\unit \yup=#2\unit
   \advance \yup by -0.5\thickness
   \killglue\rlap{\kern \xup\raise\yup\hbox{\arrw}}\ignorespaces}

\def \EA #1 #2 {\xup=#1\unit \yup=#2\unit
   \advance \yup by -0.5\thickness
   \advance \xup by -0.9\unit
   \killglue\rlap{\hskip \xup\raise\yup\hbox{\arre}}\ignorespaces}

\long\def\NORTH #1 #2 #3 { \yyyup=#3\unit\relax
\xup=#1\unit\yup=#2\unit
\yyup=\yyyup
\advance \yup by -#3\unit
\advance \yyyup by -0.2\unit
\advance \xup by -0.5\thickness
\raise\yup\hbox to 0pt{\hskip \xup \vrule
height \yyyup width \thickness depth 0pt\hss}\ignorespaces}

\long\def\SOUTH #1 #2 #3 {\killglue \yyyup=#3\unit\relax
\xup=#1\unit\yup=#2\unit
\advance \yyyup by -0.2\unit
\advance \yup by 0.2\unit
\raise\yup\hbox to 0pt{\hskip \xup\hskip -0.5\thickness \vrule
height \yyyup width \thickness depth 0pt\hss}\ignorespaces}

\long\def\WEST #1 #2 #3 {\killglue \xxup=#3\unit\relax
\xup=#1\unit\yup=#2\unit
   \advance \yup by -0.5\thickness
\advance \xup by 0.2\unit
\advance \xxup by -0.2\unit
\raise\yup\hbox to 0pt{\hskip \xup \vrule
height \thickness width \xxup depth 0pt\hss}\ignorespaces}

\long\def\EAST #1 #2 #3 {\killglue \xxup=#3\unit\relax
\xup=#1\unit\yup=#2\unit
   \advance \yup by -0.5\thickness
\advance \xup by -\xxup
\advance \xxup by -0.2\unit
\raise\yup\hbox to 0pt{\hskip \xup \vrule
height \thickness width \xxup depth 0pt\hss}\ignorespaces}

\def\sams{\hbox{\lline\char'122}}
\def\ssams{\hbox{\lline\char'152}}
\def\sssams{\hbox{\lline\char'161}}
\def\ssssams{\hbox{\lline\char'172}}
\def\samss{\hbox{\lline\char'125}}
\def\sssamss{\hbox{\lline\char'163}}
\def\samsss{\hbox{\lline\char'116}}
\def\samssss{\hbox{\lline\char'127}}
\def\ssamsss{\hbox{\lline\char'136}}
\def\ssssamsss{\hbox{\lline\char'176}}
\def\sssamssss{\hbox{\lline\char'167}}

\def\spas{\hbox{\lline\char'022}}
\def\sspas{\hbox{\lline\char'052}}
\def\ssspas{\hbox{\lline\char'061}}
\def\sssspas{\hbox{\lline\char'072}}
\def\spass{\hbox{\lline\char'025}}
\def\ssspass{\hbox{\lline\char'063}}
\def\spasss{\hbox{\lline\char'016}}
\def\spassss{\hbox{\lline\char'027}}
\def\sspasss{\hbox{\lline\char'036}}
\def\sssspasss{\hbox{\lline\char'076}}
\def\ssspassss{\hbox{\lline\char'067}}

\def\smas{\hbox{\lline\char'011}}
\def\smass{\hbox{\lline\char'013}}
\def\smasss{\hbox{\lline\char'015}}
\def\smassss{\hbox{\lline\char'017}}
\def\ssmas{\hbox{\lline\char'031}}
\def\ssmasss{\hbox{\lline\char'035}}
\def\sssmas{\hbox{\lline\char'051}}
\def\sssmass{\hbox{\lline\char'053}}
\def\sssmassss{\hbox{\lline\char'057}}
\def\ssssmas{\hbox{\lline\char'071}}
\def\ssssmasss{\hbox{\lline\char'075}}

\def\saps{\hbox{\lline\char'111}}
\def\sapss{\hbox{\lline\char'113}}
\def\sapsss{\hbox{\lline\char'115}}
\def\sapssss{\hbox{\lline\char'117}}
\def\ssaps{\hbox{\lline\char'131}}
\def\ssapsss{\hbox{\lline\char'135}}
\def\sssaps{\hbox{\lline\char'151}}
\def\sssapss{\hbox{\lline\char'153}}
\def\sssapssss{\hbox{\lline\char'157}}
\def\ssssaps{\hbox{\lline\char'171}}
\def\ssssapsss{\hbox{\lline\char'175}}

\def\RNEAR #1 #2 {\xup=#1\unit \yup=#2\unit  
\advance \xup by -0.8\unit \advance \yup by -0.8\unit
\killglue\rlap{\kern \xup\raise\yup\hbox{\spas}}\ignorespaces}

\def\RNEARR #1 #2 {\xup=#1\unit \yup=#2\unit
  \advance \xup by -0.8\wd\csname82\endcsname
\advance \yup by -0.8\ht\csname82\endcsname
      \killglue\rlap{\kern \xup\raise\yup\hbox{\spass}}\ignorespaces}
\def\RNEARRR #1 #2 {\xup=#1\unit \yup=#2\unit
 \advance \xup by -0.8\wd\csname83\endcsname \advance \yup by -0.8\ht\csname83\endcsname
           \killglue\rlap{\kern
\xup\raise\yup\hbox{\spasss}}\ignorespaces}
\def\RRNEARRR #1 #2 {\xup=#1\unit \yup=#2\unit
   \advance \xup by -0.8\wd\csname23\endcsname \advance \yup by -0.8\ht\csname23\endcsname
           \killglue\rlap{\kern \xup\raise\yup\hbox{\sspasss}}\ignorespaces}
\def\RRRNEAR #1 #2 {\xup=#1\unit \yup=#2\unit
    \advance \xup by -0.8\wd\csname31\endcsname \advance \yup by -0.8\ht\csname31\endcsname
           \killglue\rlap{\kern \xup\raise\yup\hbox{\ssspas}}\ignorespaces}
\def\RRRNEARR #1 #2 {\xup=#1\unit \yup=#2\unit
    \advance \xup by -0.8\wd\csname32\endcsname \advance \yup by
-0.8\ht\csname32\endcsname
      \killglue\rlap{\kern \xup\raise\yup\hbox{\ssspass}}\ignorespaces}
\def\RRRNEARRRR #1 #2 {\xup=#1\unit \yup=#2\unit
 \advance \xup by -0.8\wd\csname34\endcsname
 \advance \yup by -0.8\ht\csname34\endcsname
          \killglue\rlap{\kern \xup\raise\yup\hbox{\ssspassss}}\ignorespaces}
\def\RNEARRRR #1 #2 {\xup=#1\unit \yup=#2\unit
 \advance \xup by -0.8\wd\csname84\endcsname
\advance \yup by -0.8\ht\csname84\endcsname
          \killglue\rlap{\kern \xup\raise\yup\hbox{\spassss}}\ignorespaces}
\def\RRRRNEAR #1 #2 {\xup=#1\unit \yup=#2\unit
  \advance \xup  by -0.8\wd\csname41\endcsname
\advance \yup by -0.8\ht\csname41\endcsname           \killglue\rlap{\kern \xup\raise\yup\hbox{\sssspas}}\ignorespaces}
\def\RRNEAR #1 #2 {\xup=#1\unit \yup=#2\unit     \advance \xup
by -0.8\wd\csname21\endcsname
\advance \yup by -0.8\ht\csname21\endcsname
\killglue\rlap{\kern \xup\raise\yup\hbox{\sspas}}\ignorespaces}

\def\RRRRNEARRR #1 #2 {\xup=#1\unit \yup=#2\unit     \advance
\xup by -0.8\wd\csname43\endcsname \advance \yup by -0.8\ht\csname43\endcsname            \killglue\rlap{\kern \xup\raise\yup\hbox{\sssspasss}}\ignorespaces}

\def\RNWAR #1 #2 {\xup=#1\unit \yup=#2\unit  
\advance \yup by -\unit         \killglue\rlap{\kern \xup\raise\yup\hbox{\saps}}\ignorespaces} 
\def\RNWARR #1 #2 {\xup=#1\unit \yup=#2\unit
\advance \yup by -\ht\csname82\endcsname          \killglue\rlap{\kern \xup\raise\yup\hbox{\sapss}}\ignorespaces}
\def\RNWARRR #1 #2 {\xup=#1\unit \yup=#2\unit
\advance \yup by -\ht\csname83\endcsname            \killglue\rlap{\kern \xup\raise\yup\hbox{\sapsss}}\ignorespaces}
\def\RNWARRRR #1 #2 {\xup=#1\unit \yup=#2\unit
\advance \yup by -\ht\csname84\endcsname            \killglue\rlap{\kern \xup\raise\yup\hbox{\sapssss}}\ignorespaces}
\def\RRNWAR #1 #2 {\xup=#1\unit \yup=#2\unit
\advance \yup by -\ht\csname21\endcsname            \killglue\rlap{\kern \xup\raise\yup\hbox{\ssaps}}\ignorespaces}
\def\RRNWARRR #1 #2 {\xup=#1\unit \yup=#2\unit
\advance \yup by -\ht\csname23\endcsname            \killglue\rlap{\kern \xup\raise\yup\hbox{\ssapsss}}\ignorespaces}
\def\RRRNWAR #1 #2 {\xup=#1\unit \yup=#2\unit
\advance \yup by -\ht\csname31\endcsname            \killglue\rlap{\kern \xup\raise\yup\hbox{\sssaps}}\ignorespaces}
\def\RRRNWARR #1 #2 {\xup=#1\unit \yup=#2\unit
\advance \yup by -\ht\csname32\endcsname            \killglue\rlap{\kern \xup\raise\yup\hbox{\sssapss}}\ignorespaces}
\def\RRRRNWAR #1 #2 {\xup=#1\unit \yup=#2\unit
\advance \yup by -\ht\csname41\endcsname            \killglue\rlap{\kern \xup\raise\yup\hbox{\ssssaps}}\ignorespaces}

\def\RRRNWARRRR #1 #2 {\xup=#1\unit \yup=#2\unit
\advance \yup by -\ht\csname45\endcsname
        \killglue\rlap{\kern \xup\raise\yup\hbox{\sssapssss}}\ignorespaces}

\def\RRRRNWARRR #1 #2 {\xup=#1\unit \yup=#2\unit
\advance \yup by -\ht\csname43\endcsname
    \killglue\rlap{\kern \xup\raise\yup\hbox{\ssssapsss}}\ignorespaces}

\def\RSEAR #1 #2 {\xup=#1\unit \yup=#2\unit
   \advance \xup by -\unit
      \killglue\rlap{\kern \xup\raise\yup\hbox{\sams}}\ignorespaces}
\def\RSEARR #1 #2 {\xup=#1\unit \yup=#2\unit     \advance \xup by -\wd\csname82\endcsname
      \killglue\rlap{\kern \xup\raise\yup\hbox{\samss}}\ignorespaces}
\def\RSEARRR #1 #2 {\xup=#1\unit \yup=#2\unit    \advance \xup by -\wd\csname83\endcsname
      \killglue\rlap{\kern \xup\raise\yup\hbox{\samsss}}\ignorespaces}
\def\RSEARRRR #1 #2 {\xup=#1\unit \yup=#2\unit    \advance \xup
by -\wd\csname84\endcsname \killglue\rlap{\kern \xup\raise\yup\hbox{\samssss}}\ignorespaces}
\def\RRSEAR #1 #2 {\xup=#1\unit \yup=#2\unit     \advance \xup by
-\wd\csname21\endcsname
      \killglue\rlap{\kern \xup\raise\yup\hbox{\ssams}}\ignorespaces}
\def\RRSEARRR #1 #2 {\xup=#1\unit \yup=#2\unit     \advance \xup
by -\wd\csname23\endcsname
      \killglue\rlap{\kern \xup\raise\yup\hbox{\ssamsss}}\ignorespaces}
\def\RRRSEAR #1 #2 {\xup=#1\unit \yup=#2\unit     \advance \xup by -\wd\csname31\endcsname
      \killglue\rlap{\kern \xup\raise\yup\hbox{\sssams}}\ignorespaces}
\def\RRRSEARR #1 #2 {\xup=#1\unit \yup=#2\unit     \advance \xup
by -\wd\csname32\endcsname
      \killglue\rlap{\kern \xup\raise\yup\hbox{\sssamss}}\ignorespaces}
\def\RRRSEARRRR #1 #2 {\xup=#1\unit \yup=#2\unit     \advance
\xup by -\wd\csname34\endcsname
      \killglue\rlap{\kern \xup\raise\yup\hbox{\sssamssss}}\ignorespaces}
\def\RRRRSEAR #1 #2 {\xup=#1\unit \yup=#2\unit     \advance \xup
by -\wd\csname41\endcsname
      \killglue\rlap{\kern \xup\raise\yup\hbox{\ssssams}}\ignorespaces}
\def\RRRRSEARRR #1 #2 {\xup=#1\unit \yup=#2\unit     \advance
\xup by -\wd\csname43\endcsname  \killglue\rlap{\kern \xup\raise\yup\hbox{\ssssamsss}}\ignorespaces}

\def\RSWAR #1 #2 {\xup=#1\unit \yup=#2\unit
      \killglue\rlap{\kern \xup\raise\yup\hbox{\smas}}\ignorespaces}
\def\RSWARR #1 #2 {\xup=#1\unit \yup=#2\unit
      \killglue\rlap{\kern \xup\raise\yup\hbox{\smass}}\ignorespaces}
\def\RSWARRR #1 #2 {\xup=#1\unit \yup=#2\unit
      \killglue\rlap{\kern \xup\raise\yup\hbox{\smasss}}\ignorespaces}
\def\RSWARRRR #1 #2 {\xup=#1\unit \yup=#2\unit
      \killglue\rlap{\kern \xup\raise\yup\hbox{\smassss}}\ignorespaces}
\def\RRSWAR #1 #2 {\xup=#1\unit \yup=#2\unit
      \killglue\rlap{\kern \xup\raise\yup\hbox{\ssmas}}\ignorespaces}
\def\RRSWARRR #1 #2 {\xup=#1\unit \yup=#2\unit
      \killglue\rlap{\kern \xup\raise\yup\hbox{\ssmasss}}\ignorespaces}
\def\RRRSWAR #1 #2 {\xup=#1\unit \yup=#2\unit
      \killglue\rlap{\kern \xup\raise\yup\hbox{\sssmas}}\ignorespaces}
\def\RRRSWARR #1 #2 {\xup=#1\unit \yup=#2\unit
      \killglue\rlap{\kern \xup\raise\yup\hbox{\sssmass}}\ignorespaces}
\def\RRRSWARRRR #1 #2 {\xup=#1\unit \yup=#2\unit
      \killglue\rlap{\kern \xup\raise\yup\hbox{\sssmassss}}\ignorespaces}
\def\RRRRSWAR #1 #2 {\xup=#1\unit \yup=#2\unit
      \killglue\rlap{\kern \xup\raise\yup\hbox{\ssssmas}}\ignorespaces}
\def\RRRRSWARRR #1 #2 {\xup=#1\unit \yup=#2\unit
      \killglue\rlap{\kern \xup\raise\yup\hbox{\ssssmasss}}\ignorespaces}

\long\def\RSWR #1 #2 #3 {\killglue  
\multcnt=#3\relax
 \xup=#1\unit \yup=#2\unit  \xxup=\unit \yyup=\unit
\advance \xup by 0.2\unit \advance \yup by 0.2\unit
\whilenum \multcnt > 1 \do%
   {\raise \yup
     \hbox to 0pt{\hskip \xup \rlap{\hbox{\sps}}\hss}%
\ifnum \multcnt = 2%
     \advance\xup by 0.8\xxup      \advance\yup by 0.8\yyup
     \raise\yup
     \hbox to 0pt{\hskip \xup \rlap{\hbox{\sps}}\hss}%
    \fi
\advance\multcnt by -1
\advance\xup by \xxup \advance \yup by \yyup}\ignorespaces}

\long\def\RSWRR #1 #2 #3 {\killglue  
\multcnt=#3\relax
 \xup=#1\unit \yup=#2\unit   \xxup=\wd\csname82\endcsname \yyup=\ht\csname82\endcsname
\advance \xup by 0.2\xxup \advance \yup by 0.2\yyup
\whilenum \multcnt > 1 \do%
   {\raise \yup
     \hbox to 0pt{\hskip \xup \rlap{\hbox{\spss}}\hss}%
\ifnum \multcnt = 2%
     \advance\xup by 0.8\xxup  \advance\yup by 0.8\yyup
     \raise\yup
     \hbox to 0pt{\hskip \xup \rlap{\hbox{\spss}}\hss}%
    \fi
\advance\multcnt by -1  \advance\xup by \xxup
\advance \yup by \yyup}\ignorespaces}

\long\def\RSWRRR #1 #2 #3 {\killglue  
\multcnt=#3 \relax
 \xup=#1\unit \yup=#2\unit  \xxup=\wd\csname83\endcsname \yyup=\ht\csname83\endcsname
\advance \xup by 0.2\xxup
\advance \yup by 0.2\yyup
\whilenum \multcnt > 1 \do%
   {\raise \yup
     \hbox to 0pt{\hskip \xup \rlap{\hbox{\spsss}}\hss}%
\ifnum \multcnt = 2%
     \advance\xup by 0.8\xxup
     \advance\yup by 0.8\yyup
     \raise\yup
     \hbox to 0pt{\hskip \xup \rlap{\hbox{\spsss}}\hss}%
    \fi
\advance\multcnt by -1
\advance\xup by \xxup
\advance \yup by \yyup}\ignorespaces}

\long\def\RSWRRRR #1 #2 #3 {\killglue  
\multcnt=#3 \relax
 \xup=#1\unit \yup=#2\unit  \xxup=\wd\csname84\endcsname \yyup=\ht\csname84\endcsname
\advance \xup by 0.2\xxup
\advance \yup by 0.2\yyup
\whilenum \multcnt > 1 \do%
   {\raise \yup
     \hbox to 0pt{\hskip \xup \rlap{\hbox{\spssss}}\hss}%
\ifnum \multcnt = 2%
     \advance\xup by 0.8\xxup
     \advance\yup by 0.8\yyup
     \raise\yup
     \hbox to 0pt{\hskip \xup \rlap{\hbox{\spssss}}\hss}%
    \fi
\advance \multcnt by -1
\advance \xup by \xxup
\advance \yup by \yyup}\ignorespaces}

\long\def\RRSWR #1 #2 #3 {\killglue
\multcnt=#3\relax
 \xup=#1\unit \yup=#2\unit
 \xxup=\wd\csname21\endcsname \yyup=\ht\csname21\endcsname
\advance \xup by 0.2\xxup \advance \yup by 0.2\yyup
\whilenum \multcnt > 1 \do%
   {\raise \yup
     \hbox to 0pt{\hskip \xup \rlap{\hbox{\ssps}}\hss}%
\ifnum \multcnt = 2%
     \advance\xup by 0.8\xxup \advance\yup by 0.8\yyup
     \raise\yup
     \hbox to 0pt{\hskip \xup \rlap{\hbox{\ssps}}\hss}%
    \fi
\advance\multcnt by -1
\advance \xup by \xxup
\advance \yup by \yyup}\ignorespaces}

\long\def\RRSWRRR #1 #2 #3 {\killglue
\multcnt=#3\relax
 \xup=#1\unit \yup=#2\unit  \xxup=\wd\csname23\endcsname \yyup=\ht\csname23\endcsname
\advance \xup by 0.2\xxup
\advance \yup by 0.2\yyup
\whilenum \multcnt > 1 \do%
   {\raise \yup
     \hbox to 0pt{\hskip \xup \rlap{\hbox{\sspsss}}\hss}%
\ifnum \multcnt = 2%
     \advance\xup by 0.8\xxup
     \advance\yup by 0.8\yyup
     \raise\yup
     \hbox to 0pt{\hskip \xup \rlap{\hbox{\sspsss}}\hss}%
    \fi
\advance\multcnt by -1
\advance\xup by \xxup
\advance \yup by \yyup}\ignorespaces}

\long\def\RRRSWR #1 #2 #3 {\killglue
\multcnt=#3\relax
 \xup=#1\unit \yup=#2\unit  \xxup=\wd\csname31\endcsname \yyup=\ht\csname31\endcsname
\advance \xup by 0.2\xxup
\advance \yup by 0.2\yyup
\whilenum \multcnt > 1 \do%
   {\raise \yup
     \hbox to 0pt{\hskip \xup \rlap{\hbox{\sssps}}\hss}%
\ifnum \multcnt = 2%
     \advance\xup by 0.8\xxup
     \advance\yup by 0.8\yyup
     \raise\yup
     \hbox to 0pt{\hskip \xup \rlap{\hbox{\sssps}}\hss}%
    \fi
\advance\multcnt by -1
\advance\xup by \xxup
\advance \yup by \yyup}\ignorespaces}

\long\def\RRRSWRR #1 #2 #3 {\killglue
\multcnt=#3\relax
 \xup=#1\unit \yup=#2\unit  \xxup=\wd\csname32\endcsname \yyup=\ht\csname32\endcsname
\advance \xup by 0.2\xxup
\advance \yup by 0.2\yyup
\whilenum \multcnt > 1 \do%
   {\raise \yup
     \hbox to 0pt{\hskip \xup \rlap{\hbox{\ssspss}}\hss}%
\ifnum \multcnt = 2%
     \advance\xup by 0.8\xxup
     \advance\yup by 0.8\yyup
     \raise\yup
     \hbox to 0pt{\hskip \xup \rlap{\hbox{\ssspss}}\hss}%
    \fi
\advance\multcnt by -1
\advance\xup by \xxup
\advance \yup by \yyup}\ignorespaces}

\long\def\RRRSWRRRR #1 #2 #3 {\killglue
\multcnt=#3\relax
 \xup=#1\unit \yup=#2\unit  \xxup=\wd\csname34\endcsname \yyup=\ht\csname34\endcsname
\advance \xup by 0.2\xxup
\advance \yup by 0.2\yyup
\whilenum \multcnt > 1 \do%
   {\raise \yup
     \hbox to 0pt{\hskip \xup \rlap{\hbox{\ssspssss}}\hss}%
\ifnum \multcnt = 2%
     \advance\xup by 0.8\xxup
     \advance\yup by 0.8\yyup
     \raise\yup
     \hbox to 0pt{\hskip \xup \rlap{\hbox{\ssspssss}}\hss}%
    \fi
\advance\multcnt by -1
\advance\xup by \xxup
\advance \yup by \yyup}\ignorespaces}

\long\def\RRRRSWR #1 #2 #3 {\killglue
\multcnt=#3\relax
 \xup=#1\unit \yup=#2\unit  \xxup=\wd\csname41\endcsname \yyup=\ht\csname41\endcsname
\advance \xup by 0.2\xxup
\advance \yup by 0.2\yyup
\whilenum \multcnt > 1 \do%
   {\raise \yup
     \hbox to 0pt{\hskip \xup \rlap{\hbox{\ssssps}}\hss}%
\ifnum \multcnt = 2%
     \advance\xup by 0.8\xxup
     \advance\yup by 0.8\yyup
     \raise\yup
     \hbox to 0pt{\hskip \xup \rlap{\hbox{\ssssps}}\hss}%
    \fi
\advance\multcnt by -1
\advance\xup by \xxup
\advance \yup by \yyup}\ignorespaces}

\long\def\RRRRSWRRR #1 #2 #3 {\killglue
\multcnt=#3\relax
 \xup=#1\unit \yup=#2\unit  \xxup=\wd\csname43\endcsname \yyup=\ht\csname43\endcsname
\advance \xup by 0.2\xxup
\advance \yup by 0.2\yyup
\whilenum \multcnt > 1 \do%
   {\raise \yup
     \hbox to 0pt{\hskip \xup \rlap{\hbox{\sssspsss}}\hss}%
\ifnum \multcnt = 2%
     \advance\xup by 0.8\xxup
     \advance\yup by 0.8\yyup
     \raise\yup
     \hbox to 0pt{\hskip \xup \rlap{\hbox{\sssspsss}}\hss}%
    \fi
\advance\multcnt by -1
\advance\xup by \xxup
\advance \yup by \yyup}\ignorespaces}

\long\def\RNER #1 #2 #3 {\killglue  
\multcnt=#3\relax
 \xup=#1\unit \yup=#2\unit
 \xxup=\unit \yyup=\unit
\advance \xup by -#3\unit
\advance \yup by -#3\unit
\whilenum \multcnt > 1 \do%
   {\raise \yup
     \hbox to 0pt{\hskip \xup \rlap{\hbox{\sps}}\hss}%
\ifnum \multcnt = 2%
     \advance\xup by \xxup
     \advance\yup by \yyup
     \raise\yup
     \hbox to 0pt{\hskip \xup \rlap{\hbox{\sps}}\hss}%
 \fi
\advance\multcnt by -1
\advance\xup by \xxup
\advance \yup by \yyup}\ignorespaces}

\long\def\RNERR #1 #2 #3 {\killglue  
\multcnt=#3\relax
 \xup=#1\unit \yup=#2\unit
 \xxup=\wd\csname82\endcsname \yyup=\ht\csname82\endcsname
\advance \xup by -#3\xxup
\advance \yup by -#3\yyup
\whilenum \multcnt > 1 \do%
   {\raise \yup
     \hbox to 0pt{\hskip \xup \rlap{\hbox{\spss}}\hss}%
\ifnum \multcnt = 2%
     \advance\xup by \xxup
     \advance\yup by \yyup
     \raise\yup
     \hbox to 0pt{\hskip \xup \rlap{\hbox{\spss}}\hss}%
 \fi
\advance\multcnt by -1
\advance\xup by \xxup
\advance \yup by \yyup}\ignorespaces}

\long\def\RNERRR #1 #2 #3 {\killglue  
\multcnt=#3\relax
 \xup=#1\unit \yup=#2\unit
 \xxup=\wd\csname83\endcsname \yyup=\ht\csname83\endcsname
\advance \xup by -#3\xxup
\advance \yup by -#3\yyup
\whilenum \multcnt > 1 \do%
   {\raise \yup
     \hbox to 0pt{\hskip \xup \rlap{\hbox{\spsss}}\hss}%
\ifnum \multcnt = 2%
     \advance\xup by \xxup
     \advance\yup by \yyup
     \raise\yup
     \hbox to 0pt{\hskip \xup \rlap{\hbox{\spsss}}\hss}%
 \fi
\advance\multcnt by -1
\advance\xup by \xxup
\advance \yup by \yyup}\ignorespaces}

\long\def\RNERRRR #1 #2 #3 {\killglue  
\multcnt=#3\relax
 \xup=#1\unit \yup=#2\unit
 \xxup=\wd\csname84\endcsname \yyup=\ht\csname84\endcsname
\advance \xup by -#3\xxup
\advance \yup by -#3\yyup
\whilenum \multcnt > 1 \do%
   {\raise \yup
     \hbox to 0pt{\hskip \xup \rlap{\hbox{\spssss}}\hss}%
\ifnum \multcnt = 2%
     \advance\xup by \xxup
     \advance\yup by \yyup
     \raise\yup
     \hbox to 0pt{\hskip \xup \rlap{\hbox{\spssss}}\hss}%
 \fi
\advance\multcnt by -1
\advance\xup by \xxup
\advance \yup by \yyup}\ignorespaces}

\long\def\RRNERRR #1 #2 #3 {\killglue  
\multcnt=#3\relax
 \xup=#1\unit \yup=#2\unit
 \xxup=\wd\csname23\endcsname \yyup=\ht\csname23\endcsname
\advance \xup by -#3\xxup
\advance \yup by -#3\yyup
\whilenum \multcnt > 1 \do%
   {\raise \yup
     \hbox to 0pt{\hskip \xup \rlap{\hbox{\sspsss}}\hss}%
\ifnum \multcnt = 2%
     \advance\xup by \xxup
     \advance\yup by \yyup
     \raise\yup
     \hbox to 0pt{\hskip \xup \rlap{\hbox{\sspsss}}\hss}%
 \fi
\advance\multcnt by -1
\advance\xup by \xxup
\advance \yup by \yyup}\ignorespaces}

\long\def\RRRNER #1 #2 #3 {\killglue  
\multcnt=#3\relax
 \xup=#1\unit \yup=#2\unit
 \xxup=\wd\csname31\endcsname \yyup=\ht\csname31\endcsname
\advance \xup by -#3\xxup
\advance \yup by -#3\yyup
\whilenum \multcnt > 1 \do%
   {\raise \yup
     \hbox to 0pt{\hskip \xup \rlap{\hbox{\sssps}}\hss}%
\ifnum \multcnt = 2%
     \advance\xup by \xxup
     \advance\yup by \yyup
     \raise\yup
     \hbox to 0pt{\hskip \xup \rlap{\hbox{\sssps}}\hss}%
 \fi
\advance\multcnt by -1
\advance\xup by \xxup
\advance \yup by \yyup}\ignorespaces}

\long\def\RRRNERR #1 #2 #3 {\killglue  
\multcnt=#3\relax
 \xup=#1\unit \yup=#2\unit
 \xxup=\wd\csname32\endcsname \yyup=\ht\csname32\endcsname
\advance \xup by -#3\xxup
\advance \yup by -#3\yyup
\whilenum \multcnt > 1 \do%
   {\raise \yup
     \hbox to 0pt{\hskip \xup \rlap{\hbox{\ssspss}}\hss}%
\ifnum \multcnt = 2%
     \advance\xup by \xxup
     \advance\yup by \yyup
     \raise\yup
     \hbox to 0pt{\hskip \xup \rlap{\hbox{\ssspss}}\hss}%
 \fi
\advance\multcnt by -1
\advance\xup by \xxup
\advance \yup by \yyup}\ignorespaces}

\long\def\RRRNERRRR #1 #2 #3 {\killglue  
\multcnt=#3\relax
 \xup=#1\unit \yup=#2\unit
 \xxup=\wd\csname34\endcsname \yyup=\ht\csname34\endcsname
\advance \xup by -#3\xxup
\advance \yup by -#3\yyup
\whilenum \multcnt > 1 \do%
   {\raise \yup
     \hbox to 0pt{\hskip \xup \rlap{\hbox{\ssspssss}}\hss}%
\ifnum \multcnt = 2%
     \advance\xup by \xxup
     \advance\yup by \yyup
     \raise\yup
     \hbox to 0pt{\hskip \xup \rlap{\hbox{\ssspssss}}\hss}%
 \fi
\advance\multcnt by -1
\advance\xup by \xxup
\advance \yup by \yyup}\ignorespaces}

\long\def\RRRRNER #1 #2 #3 {\killglue  
\multcnt=#3\relax
 \xup=#1\unit \yup=#2\unit
 \xxup=\wd\csname41\endcsname \yyup=\ht\csname41\endcsname
\advance \xup by -#3\xxup
\advance \yup by -#3\yyup
\whilenum \multcnt > 1 \do%
   {\raise \yup
     \hbox to 0pt{\hskip \xup \rlap{\hbox{\ssssps}}\hss}%
\ifnum \multcnt = 2%
     \advance\xup by \xxup
     \advance\yup by \yyup
     \raise\yup
     \hbox to 0pt{\hskip \xup \rlap{\hbox{\ssssps}}\hss}%
 \fi
\advance\multcnt by -1
\advance\xup by \xxup
\advance \yup by \yyup}\ignorespaces}

\long\def\RRRRNERRR #1 #2 #3 {\killglue  
\multcnt=#3\relax
 \xup=#1\unit \yup=#2\unit
 \xxup=\wd\csname43\endcsname \yyup=\ht\csname43\endcsname
\advance \xup by -#3\xxup
\advance \yup by -#3\yyup
\whilenum \multcnt > 1 \do%
   {\raise \yup
     \hbox to 0pt{\hskip \xup \rlap{\hbox{\sssspsss}}\hss}%
\ifnum \multcnt = 2%
     \advance\xup by \xxup
     \advance\yup by \yyup
     \raise\yup
     \hbox to 0pt{\hskip \xup \rlap{\hbox{\sssspsss}}\hss}%
 \fi
\advance\multcnt by -1
\advance\xup by \xxup
\advance \yup by \yyup}\ignorespaces}

\long\def\RSER #1 #2 #3 {\killglue  
\multcnt=#3\relax
 \xup=#1\unit \yup=#2\unit  \xxup=\unit \yyup=\unit
\advance \yup by -\unit
\advance \xup by -\multcnt\unit
\advance \yup by \multcnt\unit
\whilenum \multcnt > 1 \do%
{\raise \yup
\hbox to 0pt{\hskip \xup \rlap{\hbox{\sms}}\hss}%
\ifnum \multcnt = 2%
\advance\xup by 0.8\xxup
\advance\yup by -0.8\yyup \raise%
\yup\hbox to 0pt{\hskip \xup \rlap{\hbox{\sms}}\hss}%
\fi \advance\multcnt by -1 \advance\xup by
\xxup \advance \yup by -\yyup}\ignorespaces}

\long\def\RSERR #1 #2 #3 {\killglue
\multcnt=#3\relax
 \xup=#1\unit \yup=#2\unit  \xxup=\wd\csname82\endcsname \yyup=\ht\csname82\endcsname
\advance \yup by -\unit
\advance \xup by -\multcnt\wd\csname82\endcsname
\advance \yup by \multcnt\ht\csname82\endcsname
\whilenum \multcnt > 1 \do%
{\raise \yup
\hbox to 0pt{\hskip \xup \rlap{\hbox{\smss}}\hss}%
\ifnum \multcnt = 2%
\advance\xup by 0.8\xxup
\advance\yup by -0.8\yyup \raise%
\yup\hbox to 0pt{\hskip \xup \rlap{\hbox{\smss}}\hss}%
\fi \advance\multcnt by -1 \advance\xup by
\xxup \advance \yup by -\yyup}\ignorespaces}

\long\def\RSERRR #1 #2 #3 {\killglue
\multcnt=#3\relax
 \xup=#1\unit \yup=#2\unit  \xxup=\wd\csname83\endcsname \yyup=\ht\csname83\endcsname
\advance \yup by -\unit
\advance \xup by -\multcnt\xxup
\advance \yup by \multcnt\unit
\whilenum \multcnt > 1 \do%
{\raise \yup\hbox to 0pt{\hskip \xup \rlap{\hbox{\smsss}}\hss}%
\ifnum \multcnt = 2%
\advance\xup by 0.8\xxup
\advance\yup by -0.8\yyup \raise%
\yup\hbox to 0pt{\hskip \xup \rlap{\hbox{\smsss}}\hss}%
\fi \advance\multcnt by -1 \advance\xup by \xxup \advance \yup by
-\yyup}\ignorespaces}

\long\def\RSERRRR #1 #2 #3 {\killglue
\multcnt=#3\relax
 \xup=#1\unit \yup=#2\unit  \xxup=\wd\csname84\endcsname \yyup=\ht\csname84\endcsname
\advance \yup by -\unit
\advance \xup by -\multcnt\xxup
\advance \yup by \multcnt\unit
\whilenum \multcnt > 1 \do%
{\raise \yup\hbox to 0pt{\hskip \xup \rlap{\hbox{\smssss}}\hss}%
\ifnum \multcnt = 2%
\advance\xup by 0.8\xxup
\advance\yup by -0.8\yyup \raise%
\yup\hbox to 0pt{\hskip \xup \rlap{\hbox{\smssss}}\hss}%
\fi \advance\multcnt by -1 \advance\xup by \xxup \advance \yup by
-\yyup}\ignorespaces}

\long\def\RRSER #1 #2 #3 {\killglue
\multcnt=#3\relax
 \xup=#1\unit \yup=#2\unit\xxup=\wd\csname21\endcsname
\yyup=\ht\csname21\endcsname
\advance \yup by -\yyup
\advance \xup by -\multcnt\xxup
\advance \yup by \multcnt\yyup
\whilenum \multcnt > 1 \do%
{\raise \yup\hbox to 0pt{\hskip \xup \rlap{\hbox{\ssms}}\hss}%
\ifnum \multcnt = 2%
\advance\xup by 0.8\xxup
\advance\yup by -0.8\yyup \raise%
\yup\hbox to 0pt{\hskip \xup \rlap{\hbox{\ssms}}\hss}%
\fi \advance\multcnt by -1 \advance\xup by \xxup \advance \yup by
-\yyup}\ignorespaces}

\long\def\RRSERRR #1 #2 #3 {\killglue
\multcnt=#3\relax
 \xup=#1\unit \yup=#2\unit  \xxup=\wd\csname23\endcsname \yyup=\ht\csname23\endcsname
\advance \yup by -\yyup
\advance \xup by -\multcnt\xxup
\advance \yup by \multcnt\yyup
\whilenum \multcnt > 1 \do%
{\raise \yup\hbox to 0pt{\hskip \xup \rlap{\hbox{\ssmsss}}\hss}%
\ifnum \multcnt = 2%
\advance\xup by 0.8\xxup
\advance\yup by -0.8\yyup \raise%
\yup\hbox to 0pt{\hskip \xup \rlap{\hbox{\ssmsss}}\hss}%
\fi \advance\multcnt by -1 \advance\xup by \xxup \advance \yup by
-\yyup}\ignorespaces}

\long\def\RRRSER #1 #2 #3 {\killglue
\multcnt=#3\relax
 \xup=#1\unit \yup=#2\unit  \xxup=\wd\csname31\endcsname \yyup=\ht\csname31\endcsname
\advance \yup by -\yyup
\advance \xup by -\multcnt\xxup
\advance \yup by \multcnt\yyup
\whilenum \multcnt > 1 \do%
{\raise \yup\hbox to 0pt{\hskip \xup \rlap{\hbox{\sssms}}\hss}%
\ifnum \multcnt = 2%
\advance\xup by 0.8\xxup
\advance\yup by -0.8\yyup \raise%
\yup\hbox to 0pt{\hskip \xup \rlap{\hbox{\sssms}}\hss}%
\fi \advance\multcnt by -1 \advance\xup by \xxup \advance \yup by
-\yyup}\ignorespaces}

\long\def\RRRSERR #1 #2 #3 {\killglue
\multcnt=#3\relax
 \xup=#1\unit \yup=#2\unit  \xxup=\wd\csname32\endcsname \yyup=\ht\csname32\endcsname
\advance \yup by -\yyup
\advance \xup by -\multcnt\xxup
\advance \yup by \multcnt\yyup
\whilenum \multcnt > 1 \do%
{\raise \yup\hbox to 0pt{\hskip \xup \rlap{\hbox{\sssmss}}\hss}%
\ifnum \multcnt = 2%
\advance\xup by 0.8\xxup
\advance\yup by -0.8\yyup \raise%
\yup\hbox to 0pt{\hskip \xup \rlap{\hbox{\sssmss}}\hss}%
\fi \advance\multcnt by -1 \advance\xup by \xxup \advance \yup by
-\yyup}\ignorespaces}

\long\def\RRRSERRRR #1 #2 #3 {\killglue
\multcnt=#3\relax
 \xup=#1\unit \yup=#2\unit  \xxup=\wd\csname34\endcsname \yyup=\ht\csname34\endcsname
\advance \yup by -\yyup
\advance \xup by -\multcnt\xxup
\advance \yup by \multcnt\yyup
\whilenum \multcnt > 1 \do%
{\raise \yup\hbox to 0pt{\hskip \xup \rlap{\hbox{\sssmssss}}\hss}%
\ifnum \multcnt = 2%
\advance\xup by 0.8\xxup
\advance\yup by -0.8\yyup \raise%
\yup\hbox to 0pt{\hskip \xup \rlap{\hbox{\sssmssss}}\hss}%
\fi \advance\multcnt by -1 \advance\xup by \xxup \advance \yup by
-\yyup}\ignorespaces}

\long\def\RRRRSER #1 #2 #3 {\killglue
\multcnt=#3\relax
 \xup=#1\unit \yup=#2\unit  \xxup=\wd\csname41\endcsname \yyup=\ht\csname41\endcsname
\advance \yup by -\yyup
\advance \xup by -\multcnt\xxup
\advance \yup by \multcnt\yyup
\whilenum \multcnt > 1 \do%
{\raise \yup\hbox to 0pt{\hskip \xup \rlap{\hbox{\ssssms}}\hss}%
\ifnum \multcnt = 2%
\advance\xup by 0.8\xxup
\advance\yup by -0.8\yyup \raise%
\yup\hbox to 0pt{\hskip \xup \rlap{\hbox{\ssssms}}\hss}%
\fi \advance\multcnt by -1 \advance\xup by \xxup \advance \yup by
-\yyup}\ignorespaces}

\long\def\RRRRSERRR #1 #2 #3 {\killglue
\multcnt=#3\relax
 \xup=#1\unit \yup=#2\unit  \xxup=\wd\csname43\endcsname \yyup=\ht\csname43\endcsname
\advance \yup by -\yyup
\advance \xup by -\multcnt\xxup
\advance \yup by \multcnt\yyup
\whilenum \multcnt > 1 \do%
{\raise \yup\hbox to 0pt{\hskip \xup \rlap{\hbox{\ssssmsss}}\hss}%
\ifnum \multcnt = 2%
\advance\xup by 0.8\xxup
\advance\yup by -0.8\yyup \raise%
\yup\hbox to 0pt{\hskip \xup \rlap{\hbox{\ssssmsss}}\hss}%
\fi \advance\multcnt by -1 \advance\xup by \xxup \advance \yup by
-\yyup}\ignorespaces}

\long\def\RNWR #1 #2 #3 {\killglue  
\multcnt=#3\relax
 \xup=#1\unit \yup=#2\unit
 \xxup=\unit \yyup=\unit
\advance \xup by 0.2\unit   
\advance \yup by  -0.2\unit    
\advance \yup by -\unit
\whilenum \multcnt > 1 \do{\raise%
\yup\hbox to 0pt{\hskip \xup \rlap{\hbox{\sms}}\hss}\ifnum
\multcnt = 2
\advance\xup by 0.8\unit 
\advance\yup by -0.8\unit 
\raise%
\yup\hbox to 0pt{\hskip \xup \rlap{\hbox{\sms}}\hss}%
\fi \advance\multcnt by -1 \advance\xup by \xxup \advance \yup by -\yyup}\ignorespaces}

\long\def\RNWRR #1 #2 #3 {\killglue
\multcnt=#3\relax
 \xup=#1\unit \yup=#2\unit
 \xxup=\wd\csname82\endcsname \yyup=\ht\csname82\endcsname
\advance \xup by 0.2\xxup   
\advance \yup by  -0.2\yyup    
\advance \yup by -\yyup
\whilenum \multcnt > 1 \do{\raise%
\yup\hbox to 0pt{\hskip \xup \rlap{\hbox{\smss}}\hss}\ifnum
\multcnt = 2
\advance\xup by 0.8\xxup 
\advance\yup by -0.8\yyup 
\raise%
\yup\hbox to 0pt{\hskip \xup \rlap{\hbox{\smss}}\hss}%
\fi \advance\multcnt by -1 \advance\xup by \xxup \advance \yup by -\yyup}\ignorespaces}

\long\def\RNWRRR #1 #2 #3 {\killglue
\multcnt=#3\relax
 \xup=#1\unit \yup=#2\unit
 \xxup=\wd\csname83\endcsname \yyup=\ht\csname83\endcsname
\advance \xup by 0.2\xxup   
\advance \yup by  -0.2\yyup    
\advance \yup by -\yyup
\whilenum \multcnt > 1 \do{\raise%
\yup\hbox to 0pt{\hskip \xup \rlap{\hbox{\smsss}}\hss}\ifnum
\multcnt = 2
\advance\xup by 0.8\xxup 
\advance\yup by -0.8\yyup 
\raise%
\yup\hbox to 0pt{\hskip \xup \rlap{\hbox{\smsss}}\hss}%
\fi \advance\multcnt by -1 \advance\xup by \xxup \advance \yup by -\yyup}\ignorespaces}

\long\def\RNWRRRR #1 #2 #3 {\killglue
\multcnt=#3\relax
 \xup=#1\unit \yup=#2\unit
 \xxup=\wd\csname84\endcsname \yyup=\ht\csname84\endcsname
\advance \xup by 0.2\xxup   
\advance \yup by  -0.2\yyup    
\advance \yup by -\yyup
\whilenum \multcnt > 1 \do{\raise%
\yup\hbox to 0pt{\hskip \xup \rlap{\hbox{\smssss}}\hss}\ifnum
\multcnt = 2
\advance\xup by 0.8\xxup 
\advance\yup by -0.8\yyup 
\raise%
\yup\hbox to 0pt{\hskip \xup \rlap{\hbox{\smssss}}\hss}%
\fi \advance\multcnt by -1 \advance\xup by \xxup \advance \yup by -\yyup}\ignorespaces}

\long\def\RRNER #1 #2 #3 {\killglue  
\multcnt=#3\relax
 \xup=#1\unit \yup=#2\unit
 \xxup=\wd\csname21\endcsname \yyup=\ht\csname21\endcsname
\advance \xup by -#3\xxup
\advance \yup by -#3\yyup
\whilenum \multcnt > 1 \do%
   {\raise \yup
     \hbox to 0pt{\hskip \xup \rlap{\hbox{\ssps}}\hss}%
\ifnum \multcnt = 2%
     \advance\xup by \xxup
     \advance\yup by \yyup
     \raise\yup
     \hbox to 0pt{\hskip \xup \rlap{\hbox{\ssps}}\hss}%
 \fi
\advance\multcnt by -1
\advance\xup by \xxup
\advance \yup by \yyup}\ignorespaces}

\long\def\RRNWR #1 #2 #3 {\killglue
\multcnt=#3\relax
 \xup=#1\unit \yup=#2\unit
 \xxup=\wd\csname21\endcsname \yyup=\ht\csname21\endcsname
\advance \xup by 0.2\xxup   
\advance \yup by  -0.2\yyup    
\advance \yup by -\yyup
\whilenum \multcnt > 1 \do{\raise%
\yup\hbox to 0pt{\hskip \xup \rlap{\hbox{\ssms}}\hss}\ifnum
\multcnt = 2
\advance\xup by 0.8\xxup 
\advance\yup by -0.8\yyup 
\raise%
\yup\hbox to 0pt{\hskip \xup \rlap{\hbox{\ssms}}\hss}%
\fi \advance\multcnt by -1 \advance\xup by \xxup
 \advance \yup by -\yyup}\ignorespaces}

\long\def\RRNWRRR #1 #2 #3 {\killglue
\multcnt=#3\relax
 \xup=#1\unit \yup=#2\unit
 \xxup=\wd\csname23\endcsname \yyup=\ht\csname23\endcsname
\advance \xup by 0.2\xxup   
\advance \yup by  -0.2\yyup    
\advance \yup by -\yyup
\whilenum \multcnt > 1 \do{\raise%
\yup\hbox to 0pt{\hskip \xup \rlap{\hbox{\ssmsss}}\hss}\ifnum
\multcnt = 2
\advance\xup by 0.8\xxup 
\advance\yup by -0.8\yyup 
\raise%
\yup\hbox to 0pt{\hskip \xup \rlap{\hbox{\ssmsss}}\hss}%
\fi \advance\multcnt by -1 \advance\xup by \xxup \advance \yup by -\yyup}\ignorespaces}

\long\def\RRRNWR #1 #2 #3 {\killglue
\multcnt=#3\relax
 \xup=#1\unit \yup=#2\unit
 \xxup=\wd\csname31\endcsname \yyup=\ht\csname31\endcsname
\advance \xup by 0.2\xxup   
\advance \yup by  -0.2\yyup    
\advance \yup by -\yyup
\whilenum \multcnt > 1 \do{\raise%
\yup\hbox to 0pt{\hskip \xup \rlap{\hbox{\sssms}}\hss}\ifnum
\multcnt = 2
\advance\xup by 0.8\xxup 
\advance\yup by -0.8\yyup 
\raise%
\yup\hbox to 0pt{\hskip \xup \rlap{\hbox{\sssms}}\hss}%
\fi \advance\multcnt by -1 \advance\xup by \xxup \advance \yup by -\yyup}\ignorespaces}

\long\def\RRRNWRR #1 #2 #3 {\killglue
\multcnt=#3\relax
 \xup=#1\unit \yup=#2\unit
 \xxup=\wd\csname32\endcsname \yyup=\ht\csname32\endcsname
\advance \xup by 0.2\xxup   
\advance \yup by  -0.2\yyup    
\advance \yup by -\yyup
\whilenum \multcnt > 1 \do{\raise%
\yup\hbox to 0pt{\hskip \xup \rlap{\hbox{\sssmss}}\hss}\ifnum
\multcnt = 2
\advance\xup by 0.8\xxup 
\advance\yup by -0.8\yyup 
\raise%
\yup\hbox to 0pt{\hskip \xup \rlap{\hbox{\sssmss}}\hss}%
\fi \advance\multcnt by -1 \advance\xup by \xxup \advance \yup by -\yyup}\ignorespaces}

\long\def\RRRNWRRRR #1 #2 #3 {\killglue
\multcnt=#3\relax
 \xup=#1\unit \yup=#2\unit
 \xxup=\wd\csname34\endcsname \yyup=\ht\csname34\endcsname
\advance \xup by 0.2\xxup   
\advance \yup by  -0.2\yyup    
\advance \yup by -\yyup
\whilenum \multcnt > 1 \do{\raise%
\yup\hbox to 0pt{\hskip \xup \rlap{\hbox{\sssmssss}}\hss}\ifnum
\multcnt = 2
\advance\xup by 0.8\xxup 
\advance\yup by -0.8\yyup 
\raise%
\yup\hbox to 0pt{\hskip \xup \rlap{\hbox{\sssmssss}}\hss}%
\fi \advance\multcnt by -1 \advance\xup by \xxup \advance \yup by -\yyup}\ignorespaces}

\long\def\RRRRNWR #1 #2 #3 {\killglue
\multcnt=#3\relax
 \xup=#1\unit \yup=#2\unit
 \xxup=\wd\csname41\endcsname \yyup=\ht\csname41\endcsname
\advance \xup by 0.2\xxup   
\advance \yup by  -0.2\yyup    
\advance \yup by -\yyup
\whilenum \multcnt > 1 \do{\raise%
\yup\hbox to 0pt{\hskip \xup \rlap{\hbox{\ssssms}}\hss}\ifnum
\multcnt = 2
\advance\xup by 0.8\xxup 
\advance\yup by -0.8\yyup 
\raise%
\yup\hbox to 0pt{\hskip \xup \rlap{\hbox{\ssssms}}\hss}%
\fi \advance\multcnt by -1 \advance\xup by \xxup \advance \yup by -\yyup}\ignorespaces}

\long\def\RRRRNWRRR #1 #2 #3 {\killglue
\multcnt=#3\relax
 \xup=#1\unit \yup=#2\unit
 \xxup=\wd\csname43\endcsname \yyup=\ht\csname43\endcsname
\advance \xup by 0.2\xxup   
\advance \yup by  -0.2\yyup    
\advance \yup by -\yyup
\whilenum \multcnt > 1 \do{\raise%
\yup\hbox to 0pt{\hskip \xup \rlap{\hbox{\ssssmsss}}\hss}\ifnum
\multcnt = 2
\advance\xup by 0.8\xxup 
\advance\yup by -0.8\yyup 
\raise%
\yup\hbox to 0pt{\hskip \xup \rlap{\hbox{\ssssmsss}}\hss}%
\fi \advance\multcnt by -1 \advance\xup by \xxup \advance \yup by -\yyup}\ignorespaces}

\expandafter\newbox\csname210\endcsname
\expandafter\newbox\csname211\endcsname
\expandafter\newbox\csname212\endcsname
\expandafter\newbox\csname213\endcsname
\expandafter\newbox\csname220\endcsname
\expandafter\newbox\csname221\endcsname
\expandafter\newbox\csname222\endcsname
\expandafter\newbox\csname223\endcsname
\expandafter\newbox\csname244\endcsname
\expandafter\newbox\csname245\endcsname
\expandafter\newbox\csname246\endcsname
\expandafter\newbox\csname247\endcsname

\expandafter\newbox\csname111\endcsname
\expandafter\newbox\csname81\endcsname
\expandafter\newbox\csname77\endcsname
\expandafter\newbox\csname82\endcsname
\expandafter\newbox\csname83\endcsname
\expandafter\newbox\csname84\endcsname
\expandafter\newbox\csname85\endcsname
\expandafter\newbox\csname86\endcsname
\expandafter\newbox\csname21\endcsname
\expandafter\newbox\csname23\endcsname
\expandafter\newbox\csname25\endcsname
\expandafter\newbox\csname31\endcsname
\expandafter\newbox\csname32\endcsname
\expandafter\newbox\csname34\endcsname
\expandafter\newbox\csname35\endcsname
\expandafter\newbox\csname41\endcsname
\expandafter\newbox\csname43\endcsname
\expandafter\newbox\csname45\endcsname
\expandafter\newbox\csname51\endcsname
\expandafter\newbox\csname52\endcsname
\expandafter\newbox\csname53\endcsname
\expandafter\newbox\csname54\endcsname
\expandafter\newbox\csname56\endcsname
\expandafter\newbox\csname61\endcsname
\expandafter\newbox\csname65\endcsname
\expandafter\newbox\csname112\endcsname
\expandafter\newbox\csname113\endcsname
\expandafter\newbox\csname114\endcsname
\expandafter\newbox\csname115\endcsname
\expandafter\newbox\csname116\endcsname
\expandafter\newbox\csname121\endcsname
\expandafter\newbox\csname123\endcsname
\expandafter\newbox\csname125\endcsname
\expandafter\newbox\csname131\endcsname
\expandafter\newbox\csname132\endcsname
\expandafter\newbox\csname134\endcsname
\expandafter\newbox\csname135\endcsname
\expandafter\newbox\csname141\endcsname
\expandafter\newbox\csname143\endcsname
\expandafter\newbox\csname145\endcsname
\expandafter\newbox\csname151\endcsname
\expandafter\newbox\csname152\endcsname
\expandafter\newbox\csname153\endcsname
\expandafter\newbox\csname154\endcsname
\expandafter\newbox\csname156\endcsname
\expandafter\newbox\csname161\endcsname
\expandafter\newbox\csname165\endcsname
\expandafter\newbox\csname244\endcsname
\expandafter\newbox\csname245\endcsname
\expandafter\newbox\csname246\endcsname
\expandafter\newbox\csname247\endcsname

\font\lline=line10
\font\cir=lcircle10 
\def\lcirc{\hbox{\cir\char'143}}
\def\scirc{\hbox{\cir\char'142}}
\def\sblob{\hbox{\cir\char'162}}
\def\lblob{\hbox{\cir\char'163}}

\newdimen\unit  \unit=10pt
\newdimen\scor
\newdimen\lcor
\newdimen\cor
\newdimen\radius
\newdimen\sradius
\newdimen\lradius
\newdimen\yyyup

\def\lcircles{\let\round=\lcirc \let\rround=\lblob \radius=1.99999pt
\cor=1.14142pt  }

\def\scircles{\let\round=\scirc \let\rround=\sblob \radius=1.5pt
\cor=1.06066pt}
\scircles 

\def\point #1 #2
 {\killglue\rlap{\kern#1\unit\raise#2\unit\hbox{\round}}\ignorespaces}
\def\spoint #1 #2
 {\killglue\rlap{\kern#1\unit\raise#2\unit\hbox{\rround}}\ignorespaces}

\def\labelr #1 #2 #3 {\killglue\rlap{\kern#1\unit\raise#2\unit\hbox{\,\,\, $#3
$}}\ignorespaces}
\def\labell #1 #2 #3 {\killglue\rlap{\kern#1\unit\raise#2\unit\hbox{$#3\,\,\
$}}\ignorespaces}

\def\segp{\hbox{\lline\char'000}}
\def\segm{\hbox{\kern-\unit\lline\char'100}}

\setbox\csname111\endcsname=\hbox{\lline\char'100}
\setbox\csname81\endcsname=\hbox{\lline\char'000}
\def\sms{\hbox{\copy\csname111\endcsname}}
\def\sps{\hbox{\copy\csname81\endcsname}}

\setbox\csname244\endcsname=\hbox{\cir\char'044} 
\setbox\csname245\endcsname=\hbox{\cir\char'045}
\setbox\csname246\endcsname=\hbox{\cir\char'046}
\setbox\csname247\endcsname=\hbox{\cir\char'047}
\setbox\csname220\endcsname=\hbox{\cir\char'020} 
\setbox\csname221\endcsname=\hbox{\cir\char'021}
\setbox\csname222\endcsname=\hbox{\cir\char'022}
\setbox\csname223\endcsname=\hbox{\cir\char'023}
\setbox\csname210\endcsname=\hbox{\cir\char'010} 
\setbox\csname211\endcsname=\hbox{\cir\char'011}
\setbox\csname212\endcsname=\hbox{\cir\char'012}
\setbox\csname213\endcsname=\hbox{\cir\char'013}
 \def\trx{\rlap{\kern0.5\wd\csname220\endcsname\raise\ht\csname220\endcsname\hbox{\copy\csname220\endcsname}}}
 \def\brx{\rlap{\kern0.5\wd\csname221\endcsname\lower\ht\csname221\endcsname\hbox{\copy\csname221\endcsname}}}
 \def\blx{\rlap{\kern-0.5\wd\csname222\endcsname\lower\ht\csname222\endcsname\hbox{\copy\csname222\endcsname}}}
 \def\tlx{\rlap{\kern-0.5\wd\csname223\endcsname\raise\ht\csname223\endcsname\hbox{\copy\csname223\endcsname}}}
 \def\trvi{\rlap{\kern0.5\wd\csname210\endcsname\raise\ht\csname210\endcsname\hbox{\copy\csname210\endcsname}}}
 \def\brvi{\rlap{\kern0.5\wd\csname210\endcsname\lower\ht\csname210\endcsname\hbox{\copy\csname211\endcsname}}}
 \def\blvi{\rlap{\kern-0.5\wd\csname210\endcsname\lower\ht\csname210\endcsname\hbox{\copy\csname212\endcsname}}}
 \def\tlvi{\rlap{\kern-0.5\wd\csname210\endcsname\raise\ht\csname210\endcsname\hbox{\copy\csname213\endcsname}}}
 \def\trxx{\rlap{\kern0.5\wd\csname244\endcsname\raise\ht\csname244\endcsname\hbox{\copy\csname244\endcsname}}}
 \def\brxx{\rlap{\kern0.5\wd\csname245\endcsname\lower\ht\csname245\endcsname\hbox{\copy\csname245\endcsname}}}
 \def\blxx{\rlap{\kern-0.5\wd\csname246\endcsname\lower\ht\csname246\endcsname\hbox{\copy\csname246\endcsname}}}
 \def\tlxx{\rlap{\kern-0.5\wd\csname247\endcsname\raise\ht\csname247\endcsname\hbox{\copy\csname247\endcsname}}}

\def\ssps{\hbox{\copy\csname21\endcsname}}
\def\spss{\hbox{\copy\csname82\endcsname}}
\def\spssss{\hbox{\copy\csname84\endcsname}}
\def\spsv{\hbox{\copy\csname85\endcsname}}
\def\spsvi{\hbox{\copy\csname86\endcsname}}
\def\smsv{\hbox{\copy\csname115\endcsname}}
\def\svims{\hbox{\copy\csname161\endcsname}}
\def\svips{\hbox{\copy\csname61\endcsname}}
\def\sssps{\hbox{\copy\csname31\endcsname}}
\def\sssms{\hbox{\copy\csname131\endcsname}}
\def\spsss{\hbox{\copy\csname83\endcsname}}
\def\ssms{\hbox{\copy\csname121\endcsname}}
\def\smss{\hbox{\copy\csname112\endcsname}}
\def\smssss{\hbox{\copy\csname114\endcsname}}
\def\sssmss{\hbox{\copy\csname132\endcsname}}
\def\ssspss{\hbox{\copy\csname32\endcsname}}
\def\sspsss{\hbox{\copy\csname23\endcsname}}
\def\ssmsss{\hbox{\copy\csname123\endcsname}}
\def\sssspsss{\hbox{\copy\csname43\endcsname}}
\def\ssspssss{\hbox{\copy\csname34\endcsname}}
\def\ssssmsss{\hbox{\copy\csname143\endcsname}}
\def\smsss{\hbox{\copy\csname113\endcsname}}
\def\ssssps{\hbox{\copy\csname41\endcsname}}
\def\svips{\hbox{\copy\csname61\endcsname}}
\def\svpss{\hbox{\copy\csname52\endcsname}}
\def\sspsv{\hbox{\copy\csname25\endcsname}}
\def\sssmssss{\hbox{\copy\csname134\endcsname}}
\def\sssmsv{\hbox{\copy\csname135\endcsname}}
\def\svpsss{\hbox{\copy\csname53\endcsname}}
\def\ssmsv{\hbox{\copy\csname125\endcsname}}
\def\smsvi{\hbox{\copy\csname116\endcsname}}
\def\svmss{\hbox{\copy\csname152\endcsname}}
\def\ssssms{\hbox{\copy\csname141\endcsname}}
\def\ssspsv{\hbox{\copy\csname35\endcsname}}
\def\sssmssss{\hbox{\copy\csname134\endcsname}}   
\def\sssspsv{\hbox{\copy\csname45\endcsname}}
\def\ssssmsv{\hbox{\copy\csname145\endcsname}}
\def\svps{\hbox{\copy\csname51\endcsname}}
\def\svms{\hbox{\copy\csname151\endcsname}}
\def\svmsss{\hbox{\copy\csname153\endcsname}}
\def\svpssss{\hbox{\copy\csname54\endcsname}}
\def\svmssss{\hbox{\copy\csname154\endcsname}}
\def\svmsvi{\hbox{\copy\csname156\endcsname}}
\def\svpsvi{\hbox{\copy\csname56\endcsname}}
\def\svims{\hbox{\copy\csname161\endcsname}}
\def\svimsv{\hbox{\copy\csname165\endcsname}}
\def\svipsv{\hbox{\copy\csname65\endcsname}}

\def\thickness{\fontdimen8\lline}

\def\whilenoop#1{}
\def\whilenum#1\do #2{\ifnum #1\relax#2\iwhilenum{#1\relax#2}\fi}
\def\iwhilenum#1{\ifnum #1\let\nextwhile=\iwhilenum
 \else\let\nextwhile=\whilenoop\fi\nextwhile{#1}}

\def\whiledim#1\do #2{\ifdim #1\relax#2\iwhiledim{#1\relax#2}\fi}
\def\iwhiledim#1{\ifdim #1\let\nextwhile=\iwhiledim
\else\let\nextwhile=\whilenoop\fi\nextwhile{#1}}

\def\killglue{\unskip\whiledim \lastskip >0pt\do{\unskip}}

\long\def\putt #1 #2 #3 {\killglue\raise#2\unit\hbox to
0pt{\hskip#1\unit #3\hss}\ignorespaces}  

\long\def\puttc #1 #2 #3 {\killglue\raise#2\unit\hbox to
0pt{\hskip#1\unit \hss #3\hss}\ignorespaces}  

\newcount\multcnt
\newdimen\xup
\newdimen\yup
\newdimen\xxup
\newdimen\yyup
\newdimen\hcor
\newdimen\vcor

\def\nep{\rlap{\kern\hcor\raise\vcor\hbox{\segp}}}
\def\nwp{\rlap{\kern-\hcor\raise\vcor\hbox{\segm}}}
\def\nwm{\rlap{\kern\hcor\lower\vcor\hbox{\segm}}}
\def\nem{\rlap{\kern-\hcor\lower\vcor\hbox{\segp}}}
\def\nwwp{\rlap{\kern-\hcor\raise\vcor\hbox{\sms}}}
\def\nwwm{\rlap{\kern\hcor\lower\vcor\hbox{\sms}}}

\long\def\linep #1 #2 #3 {\killglue
\multcnt=#3\relax \hcor=\cor \vcor=\cor \xup=#1\unit \yup=#2\unit \whilenum
 \multcnt > 1 \do{\raise%
\yup\hbox to 0pt{\hskip \xup \nep\hss}\ifnum \multcnt = 2\advance\xup by \unit
\advance\yup by \unit \raise%
\yup\hbox to 0pt{\hskip \xup \nem\hss}\advance\xup by -\unit
\advance\yup by -\unit \fi \advance\multcnt by -1 \advance\xup by
\unit \advance \yup by \unit}\ignorespaces}

\long\def\slps #1 #2 #3 {\killglue
\multcnt=#3\relax \hcor=0.7071\radius
                  \vcor=0.7071\radius
 \xup=#1\unit \yup=#2\unit
 \xxup=\unit \yyup=\unit
\whilenum \multcnt > 1 \do{\raise%
\yup\hbox to 0pt{\hskip \xup \rlap{\kern\hcor\raise\vcor\hbox{\sps}}\hss}\ifnum
 \multcnt = 2\advance\xup by \xxup
\advance\yup by \yyup \raise%
\yup\hbox to 0pt{\hskip \xup \rlap{\kern-\hcor\lower\vcor\hbox{\sps}}\hss}%
\advance\xup by -\xxup
\advance\yup by -\yyup \fi \advance\multcnt by -1 \advance\xup by
\xxup \advance \yup by \yyup}\ignorespaces}

\long\def\svilpsv #1 #2 #3 {\killglue
\multcnt=#3\relax \hcor=0.78113\radius
                  \vcor=0.65094\radius
 \xup=#1\unit \yup=#2\unit
 \xxup=\wd\csname65\endcsname \yyup=\ht\csname65\endcsname
\whilenum \multcnt > 1 \do{\raise%
\yup\hbox to 0pt{\hskip \xup
 \rlap{\kern\hcor\raise\vcor\hbox{\svipsv}}\hss}\ifnum \multcnt = 2\advance\xup
 by \xxup
\advance\yup by \yyup \raise%
\yup\hbox to 0pt{\hskip \xup \rlap{\kern-\hcor\lower\vcor\hbox{\svipsv}}\hss}%
\advance\xup by -\xxup
\advance\yup by -\yyup \fi \advance\multcnt by -1 \advance\xup by
\xxup \advance \yup by \yyup}\ignorespaces}

\long\def\svlpsvi #1 #2 #3 {\killglue
\multcnt=#3\relax \hcor=0.65094\radius
                  \vcor=0.78113\radius
 \xup=#1\unit \yup=#2\unit
 \xxup=\wd\csname56\endcsname \yyup=\ht\csname56\endcsname
\whilenum \multcnt > 1 \do{\raise%
\yup\hbox to 0pt{\hskip \xup
 \rlap{\kern\hcor\raise\vcor\hbox{\svpsvi}}\hss}\ifnum \multcnt = 2\advance\xup
 by \xxup
\advance\yup by \yyup \raise%
\yup\hbox to 0pt{\hskip \xup \rlap{\kern-\hcor\lower\vcor\hbox{\svpsvi}}\hss}%
\advance\xup by -\xxup
\advance\yup by -\yyup \fi \advance\multcnt by -1 \advance\xup by
\xxup \advance \yup by \yyup}\ignorespaces}

\long\def\svlpssss #1 #2 #3 {\killglue
\multcnt=#3\relax \hcor=0.8006\radius
                  \vcor=0.6405\radius
 \xup=#1\unit \yup=#2\unit
 \xxup=\wd\csname54\endcsname \yyup=\ht\csname54\endcsname
\whilenum \multcnt > 1 \do{\raise%
\yup\hbox to 0pt{\hskip \xup
 \rlap{\kern\hcor\raise\vcor\hbox{\svpssss}}\hss}\ifnum \multcnt = 2\advance\xup
 by \xxup
\advance\yup by \yyup \raise%
\yup\hbox to 0pt{\hskip \xup \rlap{\kern-\hcor\lower\vcor\hbox{\svpssss}}\hss}%
\advance\xup by -\xxup
\advance\yup by -\yyup \fi \advance\multcnt by -1 \advance\xup by
\xxup \advance \yup by \yyup}\ignorespaces}

\long\def\svlps #1 #2 #3 {\killglue
\multcnt=#3\relax \hcor=0.98058\radius
                  \vcor=0.19611\radius
 \xup=#1\unit \yup=#2\unit
 \xxup=\wd\csname51\endcsname \yyup=\ht\csname51\endcsname
\whilenum \multcnt > 1 \do{\raise%
\yup\hbox to 0pt{\hskip \xup \rlap{\kern\hcor\raise\vcor\hbox{\svps}}\hss}\ifnum
 \multcnt = 2\advance\xup by \xxup
\advance\yup by \yyup \raise%
\yup\hbox to 0pt{\hskip \xup \rlap{\kern-\hcor\lower\vcor\hbox{\svps}}\hss}%
\advance\xup by -\xxup
\advance\yup by -\yyup \fi \advance\multcnt by -1 \advance\xup by
\xxup \advance \yup by \yyup}\ignorespaces}

\long\def\sssslpsv #1 #2 #3 {\killglue
\multcnt=#3\relax \hcor=0.6405\radius
                  \vcor=0.80006\radius
 \xup=#1\unit \yup=#2\unit
 \xxup=\wd\csname45\endcsname \yyup=\ht\csname45\endcsname
\whilenum \multcnt > 1 \do{\raise%
\yup\hbox to 0pt{\hskip \xup
 \rlap{\kern\hcor\raise\vcor\hbox{\sssspsv}}\hss}\ifnum \multcnt = 2\advance\xup
 by \xxup
\advance\yup by \yyup \raise%
\yup\hbox to 0pt{\hskip \xup \rlap{\kern-\hcor\lower\vcor\hbox{\sssspsv}}\hss}%
\advance\xup by -\xxup
\advance\yup by -\yyup \fi \advance\multcnt by -1 \advance\xup by
\xxup \advance \yup by \yyup}\ignorespaces}

\long\def\ssslpssss #1 #2 #3 {\killglue
\multcnt=#3\relax \hcor=0.6\radius   
                  \vcor=0.8\radius   
 \xup=#1\unit \yup=#2\unit
 \xxup=\wd\csname34\endcsname \yyup=\ht\csname34\endcsname
\whilenum \multcnt > 1 \do{\raise%
\yup\hbox to 0pt{\hskip \xup
 \rlap{\kern\hcor\raise\vcor\hbox{\ssspssss}}\hss}\ifnum \multcnt =
 2\advance\xup by \xxup
\advance\yup by \yyup \raise%
\yup\hbox to 0pt{\hskip \xup \rlap{\kern-\hcor\lower\vcor\hbox{\ssspssss}}\hss}%
\advance\xup by -\xxup
\advance\yup by -\yyup \fi \advance\multcnt by -1 \advance\xup by
\xxup \advance \yup by \yyup}\ignorespaces}

\long\def\sslpsss #1 #2 #3 {\killglue
\multcnt=#3\relax \hcor=0.5547\radius
                  \vcor=0.8332\radius
 \xup=#1\unit \yup=#2\unit
 \xxup=\wd\csname23\endcsname \yyup=\ht\csname23\endcsname
\whilenum \multcnt > 1 \do{\raise%
\yup\hbox to 0pt{\hskip \xup
 \rlap{\kern\hcor\raise\vcor\hbox{\sspsss}}\hss}\ifnum \multcnt = 2\advance\xup
 by \xxup
\advance\yup by \yyup \raise%
\yup\hbox to 0pt{\hskip \xup \rlap{\kern-\hcor\lower\vcor\hbox{\sspsss}}\hss}%
\advance\xup by -\xxup
\advance\yup by -\yyup \fi \advance\multcnt by -1 \advance\xup by
\xxup \advance \yup by \yyup}\ignorespaces}

\long\def\slpsvi #1 #2 #3 {\killglue
\multcnt=#3\relax \hcor=0.16439\radius
                  \vcor=0.98639\radius
 \xup=#1\unit \yup=#2\unit
 \xxup=\wd\csname86\endcsname \yyup=\ht\csname86\endcsname
\whilenum \multcnt > 1 \do{\raise%
\yup\hbox to 0pt{\hskip \xup
 \rlap{\kern\hcor\raise\vcor\hbox{\spsvi}}\hss}\ifnum \multcnt = 2\advance\xup
 by \xxup
\advance\yup by \yyup \raise%
\yup\hbox to 0pt{\hskip \xup \rlap{\kern-\hcor\lower\vcor\hbox{\spsvi}}\hss}%
\advance\xup by -\xxup
\advance\yup by -\yyup \fi \advance\multcnt by -1 \advance\xup by
\xxup \advance \yup by \yyup}\ignorespaces}

\long\def\slpsv #1 #2 #3 {\killglue
\multcnt=#3\relax \hcor=0.19611\radius
                  \vcor=0.98058\radius
 \xup=#1\unit \yup=#2\unit
 \xxup=\wd\csname85\endcsname \yyup=\ht\csname85\endcsname
\whilenum \multcnt > 1 \do{\raise%
\yup\hbox to 0pt{\hskip \xup \rlap{\kern\hcor\raise\vcor\hbox{\spsv}}\hss}\ifnum
 \multcnt = 2\advance\xup by \xxup
\advance\yup by \yyup \raise%
\yup\hbox to 0pt{\hskip \xup \rlap{\kern-\hcor\lower\vcor\hbox{\spsv}}\hss}%
\advance\xup by -\xxup
\advance\yup by -\yyup \fi \advance\multcnt by -1 \advance\xup by
\xxup \advance \yup by \yyup}\ignorespaces}

\long\def\slpssss #1 #2 #3 {\killglue
\multcnt=#3\relax \hcor=0.24253\radius
                  \vcor=0.97014\radius
 \xup=#1\unit \yup=#2\unit
 \xxup=\wd\csname84\endcsname \yyup=\ht\csname84\endcsname
\whilenum \multcnt > 1 \do{\raise%
\yup\hbox to 0pt{\hskip \xup
 \rlap{\kern\hcor\raise\vcor\hbox{\spssss}}\hss}\ifnum \multcnt = 2\advance\xup
 by \xxup
\advance\yup by \yyup \raise%
\yup\hbox to 0pt{\hskip \xup \rlap{\kern-\hcor\lower\vcor\hbox{\spssss}}\hss}%
\advance\xup by -\xxup
\advance\yup by -\yyup \fi \advance\multcnt by -1 \advance\xup by
\xxup \advance \yup by \yyup}\ignorespaces}

\long\def\sslps #1 #2 #3 {\killglue
\multcnt=#3\relax \hcor=0.89442\radius
                  \vcor=0.4472\radius
 \xup=#1\unit \yup=#2\unit
 \xxup=\wd\csname21\endcsname \yyup=\ht\csname21\endcsname
\whilenum \multcnt > 1 \do{\raise%
\yup\hbox to 0pt{\hskip \xup \rlap{\kern\hcor\raise\vcor\hbox{\ssps}}\hss}\ifnum
 \multcnt = 2\advance\xup by \xxup
\advance\yup by \yyup \raise%
\yup\hbox to 0pt{\hskip \xup \rlap{\kern-\hcor\lower\vcor\hbox{\ssps}}\hss}%
\advance\xup by -\xxup
\advance\yup by -\yyup \fi \advance\multcnt by -1 \advance\xup by
\xxup \advance \yup by \yyup}\ignorespaces}

\long\def\sslpsv #1 #2 #3 {\killglue
\multcnt=#3\relax \hcor=0.37137\radius
                  \vcor=0.92847\radius
 \xup=#1\unit \yup=#2\unit
 \xxup=\wd\csname25\endcsname \yyup=\ht\csname25\endcsname
\whilenum \multcnt > 1 \do{\raise%
\yup\hbox to 0pt{\hskip \xup
 \rlap{\kern\hcor\raise\vcor\hbox{\sspsv}}\hss}\ifnum \multcnt = 2\advance\xup
 by \xxup
\advance\yup by \yyup \raise%
\yup\hbox to 0pt{\hskip \xup \rlap{\kern-\hcor\lower\vcor\hbox{\sspsv}}\hss}%
\advance\xup by -\xxup
\advance\yup by -\yyup \fi \advance\multcnt by -1 \advance\xup by
\xxup \advance \yup by \yyup}\ignorespaces}

\long\def\slms #1 #2 #3 {\killglue
\multcnt=#3\relax \hcor=0.7071\radius
                  \vcor=0.7071\radius
 \xup=#1\unit \yup=#2\unit
 \xxup=\unit \yyup=\unit   \advance\yup by -\yyup\whilenum \multcnt > 1
 \do{\raise%
\yup\hbox to 0pt{\hskip \xup \rlap{\kern\hcor\lower\vcor\hbox{\sms}}\hss}\ifnum
 \multcnt = 2\advance\xup by \xxup
\advance\yup by -\yyup \raise%
\yup\hbox to 0pt{\hskip \xup \rlap{\kern-\hcor\raise\vcor\hbox{\sms}}\hss}%
\advance\xup by -\xxup%
\advance\yup by -\yyup \fi \advance\multcnt by -1 \advance\xup by
\xxup \advance \yup by -\yyup}\ignorespaces}

\long\def\slmssss #1 #2 #3 {\killglue
\multcnt=#3\relax \hcor=0.24253\radius
                  \vcor=0.97014\radius
 \xup=#1\unit \yup=#2\unit
 \xxup=\wd\csname84\endcsname \yyup=\ht\csname84\endcsname   \advance\yup by -\yyup\whilenum \multcnt > 1 \do{\raise%
\yup\hbox to 0pt{\hskip \xup \rlap{\kern\hcor\lower\vcor\hbox{\smssss}}\hss}%
\ifnum \multcnt = 2\advance\xup by \xxup
\advance\yup by -\yyup \raise%
\yup\hbox to 0pt{\hskip \xup \rlap{\kern-\hcor\raise\vcor\hbox{\smssss}}\hss}%
\advance\xup by -\xxup%
\advance\yup by -\yyup \fi \advance\multcnt by -1 \advance\xup by
\xxup \advance \yup by -\yyup}\ignorespaces}

\long\def\slmsv #1 #2 #3 {\killglue
\multcnt=#3\relax \hcor=0.19611\radius
                  \vcor=0.98058\radius
 \xup=#1\unit \yup=#2\unit
 \xxup=\wd\csname85\endcsname \yyup=\ht\csname85\endcsname
 \advance\yup by -\yyup
\whilenum \multcnt > 1 \do{\raise%
\yup\hbox to 0pt{\hskip \xup
\rlap{\kern\hcor\lower\vcor\hbox{\smsv}}\hss}%
\ifnum \multcnt = 2\advance\xup by \xxup
\advance\yup by -\yyup \raise%
\yup\hbox to 0pt{\hskip \xup \rlap{\kern-\hcor\raise\vcor\hbox{\smsv}}\hss}%
\advance\xup by -\xxup%
\advance\yup by -\yyup \fi \advance\multcnt by -1 \advance\xup by
\xxup \advance \yup by -\yyup}\ignorespaces}

\long\def\slmsvi #1 #2 #3 {\killglue
\multcnt=#3\relax \hcor=0.16439\radius
                  \vcor=0.98639\radius
 \xup=#1\unit \yup=#2\unit
 \xxup=\wd\csname86\endcsname \yyup=\ht\csname86\endcsname   \advance\yup by -\yyup
\whilenum \multcnt > 1  \do{\raise%
\yup\hbox to 0pt{\hskip \xup
 \rlap{\kern\hcor\lower\vcor\hbox{\smsvi}}\hss}\ifnum \multcnt = 2\advance\xup
 by \xxup
\advance\yup by -\yyup \raise%
\yup\hbox to 0pt{\hskip \xup \rlap{\kern-\hcor\raise\vcor\hbox{\smsvi}}\hss}%
\advance\xup by -\xxup%
\advance\yup by -\yyup \fi \advance\multcnt by -1 \advance\xup by
\xxup \advance \yup by -\yyup}\ignorespaces}

\long\def\sssslmsv #1 #2 #3 {\killglue
\multcnt=#3\relax \hcor=0.6405\radius
                  \vcor=0.8006\radius
 \xup=#1\unit \yup=#2\unit
 \xxup=\wd\csname45\endcsname \yyup=\ht\csname45\endcsname   \advance\yup by -\yyup
\whilenum \multcnt > 1  \do{\raise%
\yup\hbox to 0pt{\hskip \xup
 \rlap{\kern\hcor\lower\vcor\hbox{\ssssmsv}}\hss}\ifnum \multcnt = 2\advance\xup
 by \xxup
\advance\yup by -\yyup \raise%
\yup\hbox to 0pt{\hskip \xup \rlap{\kern-\hcor\raise\vcor\hbox{\ssssmsv}}\hss}%
\advance\xup by -\xxup%
\advance\yup by -\yyup \fi \advance\multcnt by -1 \advance\xup by
\xxup \advance \yup by -\yyup}\ignorespaces}

\long\def\svlms #1 #2 #3 {\killglue
\multcnt=#3\relax \hcor=0.98058\radius
                  \vcor=0.19611\radius
 \xup=#1\unit \yup=#2\unit
 \xxup=\wd\csname51\endcsname \yyup=\ht\csname51\endcsname\advance\yup by -\yyup
\whilenum \multcnt > 1  \do{\raise%
\yup\hbox to 0pt{\hskip \xup \rlap{\kern\hcor\lower\vcor\hbox{\svms}}\hss}\ifnum
 \multcnt = 2\advance\xup by \xxup
\advance\yup by -\yyup \raise%
\yup\hbox to 0pt{\hskip \xup \rlap{\kern-\hcor\raise\vcor\hbox{\svms}}\hss}%
\advance\xup by -\xxup%
\advance\yup by -\yyup \fi \advance\multcnt by -1 \advance\xup by
\xxup \advance \yup by -\yyup}\ignorespaces}

\long\def\svlmsss #1 #2 #3 {\killglue
\multcnt=#3\relax \hcor=0.8574\radius
                  \vcor=0.51449\radius
 \xup=#1\unit \yup=#2\unit
 \xxup=\wd\csname53\endcsname \yyup=\ht\csname53\endcsname   \advance\yup by -\yyup
\whilenum \multcnt > 1  \do{\raise%
\yup\hbox to 0pt{\hskip \xup%
\rlap{\kern\hcor\lower\vcor\hbox{\svmsss}}\hss}\ifnum \multcnt = 2\advance\xup
 by \xxup%
\advance\yup by -\yyup \raise%
\yup\hbox to 0pt{\hskip \xup \rlap{\kern-\hcor\raise\vcor\hbox{\svmsss}}\hss}%
\advance\xup by -\xxup%
\advance\yup by -\yyup \fi \advance\multcnt by -1 \advance\xup by
\xxup \advance \yup by -\yyup}\ignorespaces}

\long\def\svlmssss #1 #2 #3 {\killglue
\multcnt=#3\relax \hcor=0.8006\radius
                  \vcor=0.6405\radius
 \xup=#1\unit \yup=#2\unit
 \xxup=\wd\csname54\endcsname \yyup=\ht\csname54\endcsname   \advance\yup by -\yyup
\whilenum \multcnt > 1  \do{\raise%
\yup\hbox to 0pt{\hskip \xup
 \rlap{\kern\hcor\lower\vcor\hbox{\svmssss}}\hss}\ifnum \multcnt = 2\advance\xup
 by \xxup
\advance\yup by -\yyup \raise%
\yup\hbox to 0pt{\hskip \xup \rlap{\kern-\hcor\raise\vcor\hbox{\svmssss}}\hss}%
\advance\xup by -\xxup%
\advance\yup by -\yyup \fi \advance\multcnt by -1 \advance\xup by
\xxup \advance \yup by -\yyup}\ignorespaces}

\long\def\svlmsvi #1 #2 #3 {\killglue
\multcnt=#3\relax \hcor=0.65094\radius
                  \vcor=0.78113\radius
 \xup=#1\unit \yup=#2\unit
 \xxup=\wd\csname56\endcsname \yyup=\ht\csname56\endcsname   \advance\yup by -\yyup
\whilenum \multcnt > 1  \do{\raise%
\yup\hbox to 0pt{\hskip \xup
 \rlap{\kern\hcor\lower\vcor\hbox{\svmsvi}}\hss}\ifnum \multcnt = 2\advance\xup
 by \xxup
\advance\yup by -\yyup \raise%
\yup\hbox to 0pt{\hskip \xup \rlap{\kern-\hcor\raise\vcor\hbox{\svmsvi}}\hss}%
\advance\xup by -\xxup%
\advance\yup by -\yyup \fi \advance\multcnt by -1 \advance\xup by
\xxup \advance \yup by -\yyup}\ignorespaces}

\long\def\svilms #1 #2 #3 {\killglue
\multcnt=#3\relax \hcor=0.9863\radius
                  \vcor=0.16439\radius
 \xup=#1\unit \yup=#2\unit
 \xxup=\wd\csname61\endcsname \yyup=\ht\csname61\endcsname
  \advance\yup by -\yyup
\whilenum \multcnt > 1  \do{\raise%
\yup\hbox to 0pt{\hskip \xup
 \rlap{\kern\hcor\lower\vcor\hbox{\svims}}\hss}\ifnum \multcnt = 2\advance\xup
 by \xxup
\advance\yup by -\yyup \raise%
\yup\hbox to 0pt{\hskip \xup \rlap{\kern-\hcor\raise\vcor\hbox{\svims}}\hss}%
\advance\xup by -\xxup%
\advance\yup by -\yyup \fi \advance\multcnt by -1 \advance\xup by
\xxup \advance \yup by -\yyup}\ignorespaces}

\long\def\svilmsv #1 #2 #3 {\killglue
\multcnt=#3\relax \hcor=0.78113\radius
                  \vcor=0.53094\radius
 \xup=#1\unit \yup=#2\unit
 \xxup=\wd\csname65\endcsname \yyup=\ht\csname65\endcsname
  \advance\yup by -\yyup \whilenum \multcnt > 1  \do{\raise%
\yup\hbox to 0pt{\hskip \xup
 \rlap{\kern\hcor\lower\vcor\hbox{\svimsv}}\hss}\ifnum \multcnt = 2\advance\xup
 by \xxup
\advance\yup by -\yyup \raise%
\yup\hbox to 0pt{\hskip \xup \rlap{\kern-\hcor\raise\vcor\hbox{\svimsv}}\hss}%
\advance\xup by -\xxup%
\advance\yup by -\yyup \fi \advance\multcnt by -1 \advance\xup by
\xxup \advance \yup by -\yyup}\ignorespaces}

\long\def\sslms #1 #2 #3 {\killglue
\multcnt=#3\relax \hcor=0.89442\radius
                  \vcor=0.4472\radius
 \xup=#1\unit \yup=#2\unit
 \xxup=\wd\csname21\endcsname \yyup=\ht\csname21\endcsname
  \advance\yup by -\yyup \whilenum \multcnt > 1 \do{\raise%
\yup\hbox to 0pt{\hskip \xup \rlap{\kern\hcor\lower\vcor\hbox{\ssms}}\hss}\ifnum
 \multcnt = 2\advance\xup by \xxup
\advance\yup by -\yyup \raise%
\yup\hbox to 0pt{\hskip \xup \rlap{\kern-\hcor\raise\vcor\hbox{\ssms}}\hss}%
\advance\xup by -\xxup
\advance\yup by -\yyup \fi \advance\multcnt by -1 \advance\xup by
\xxup \advance \yup by -\yyup}\ignorespaces}

\long\def\sssslms #1 #2 #3 {\killglue
\multcnt=#3\relax \hcor=0.97014\radius
                  \vcor=0.24253\radius
 \xup=#1\unit \yup=#2\unit
 \xxup=\wd\csname41\endcsname \yyup=\ht\csname41\endcsname   \advance\yup by -\yyup \whilenum \multcnt > 1 \do{\raise%
\yup\hbox to 0pt{\hskip \xup
 \rlap{\kern\hcor\lower\vcor\hbox{\ssssms}}\hss}\ifnum \multcnt = 2\advance\xup
 by \xxup
\advance\yup by -\yyup \raise%
\yup\hbox to 0pt{\hskip \xup \rlap{\kern-\hcor\raise\vcor\hbox{\ssssms}}\hss}%
\advance\xup by -\xxup
\advance\yup by -\yyup \fi \advance\multcnt by -1 \advance\xup by
\xxup \advance \yup by -\yyup}\ignorespaces}

\long\def\svlmss #1 #2 #3 {\killglue
\multcnt=#3\relax \hcor=0.92847\radius
                  \vcor=0.37139\radius
 \xup=#1\unit \yup=#2\unit
 \xxup=\wd\csname152\endcsname \yyup=\ht\csname152\endcsname   \advance\yup by -\yyup
\whilenum \multcnt > 1 \do{\raise%
\yup\hbox to 0pt{\hskip \xup
 \rlap{\kern\hcor\lower\vcor\hbox{\svmss}}\hss}\ifnum \multcnt = 2\advance\xup
 by \xxup
\advance\yup by -\yyup \raise%
\yup\hbox to 0pt{\hskip \xup \rlap{\kern-\hcor\raise\vcor\hbox{\svmss}}\hss}%
\advance\xup by -\xxup
\advance\yup by -\yyup \fi \advance\multcnt by -1 \advance\xup by
\xxup \advance \yup by -\yyup}\ignorespaces}

\long\def\sslmsv #1 #2 #3 {\killglue
\multcnt=#3\relax \hcor=0.37139\radius
                  \vcor=0.92847\radius
 \xup=#1\unit \yup=#2\unit
 \xxup=\wd\csname125\endcsname \yyup=\ht\csname125\endcsname   \advance\yup by -\yyup
\whilenum \multcnt > 1 \do{\raise%
\yup\hbox to 0pt{\hskip \xup
 \rlap{\kern\hcor\lower\vcor\hbox{\ssmsv}}\hss}\ifnum \multcnt = 2\advance\xup
 by \xxup
\advance\yup by -\yyup \raise%
\yup\hbox to 0pt{\hskip \xup \rlap{\kern-\hcor\raise\vcor\hbox{\ssmsv}}\hss}%
\advance\xup by -\xxup
\advance\yup by -\yyup \fi \advance\multcnt by -1 \advance\xup by
\xxup \advance \yup by -\yyup}\ignorespaces}

\long\def\ssslmssss #1 #2 #3 {\killglue
\multcnt=#3\relax \hcor=0.6\radius   
                  \vcor=0.8\radius   
 \xup=#1\unit \yup=#2\unit
 \xxup=\wd\csname134\endcsname \yyup=\ht\csname134\endcsname
\advance\yup by -\yyup \whilenum \multcnt > 1 \do{\raise%
\yup\hbox to 0pt{\hskip \xup \rlap{\kern\hcor\lower\vcor\hbox{\sssmssss}}\hss}%
\ifnum \multcnt = 2\advance\xup by \xxup
\advance\yup by -\yyup \raise%
\yup\hbox to 0pt{\hskip \xup \rlap{\kern-\hcor\raise\vcor\hbox{\sssmssss}}\hss}%
\advance\xup by -\xxup
\advance\yup by -\yyup \fi \advance\multcnt by -1 \advance\xup by
\xxup \advance \yup by -\yyup}\ignorespaces}

\long\def\slmss #1 #2 #3 {\killglue
\multcnt=#3\relax \vcor=0.89442\radius
                  \hcor=0.4472\radius
 \xup=#1\unit \yup=#2\unit
 \xxup=\wd\csname82\endcsname \yyup=\ht\csname82\endcsname   \advance\yup by -\yyup \whilenum \multcnt > 1 \do{\raise%
\yup\hbox to 0pt{\hskip \xup
\rlap{\kern\hcor\lower\vcor\hbox{\smss}}\hss}\ifnum \multcnt = 2\advance\xup by \xxup
\advance\yup by -\yyup \raise%
\yup\hbox to 0pt{\hskip \xup \rlap{\kern-\hcor\raise\vcor\hbox{\smss}}\hss}%
\advance\xup by -\xxup
\advance\yup by -\yyup \fi \advance\multcnt by -1 \advance\xup by
\xxup \advance \yup by -\yyup}\ignorespaces}

\long\def\ssslpss #1 #2 #3 {\killglue
\multcnt=#3\relax \hcor=0.8332\radius
                  \vcor=0.5547\radius
 \xup=#1\unit \yup=#2\unit
 \xxup=\wd\csname32\endcsname \yyup=\ht\csname32\endcsname
\whilenum \multcnt > 1 \do{\raise%
\yup\hbox to 0pt{\hskip \xup
 \rlap{\kern\hcor\raise\vcor\hbox{\ssspss}}\hss}\ifnum \multcnt = 2\advance\xup
 by \xxup
\advance\yup by \yyup \raise%
\yup\hbox to 0pt{\hskip \xup \rlap{\kern-\hcor\lower\vcor\hbox{\ssspss}}\hss}%
\advance\xup by -\xxup
\advance\yup by -\yyup \fi \advance\multcnt by -1 \advance\xup by
\xxup \advance \yup by \yyup}\ignorespaces}

\long\def\sssslps #1 #2 #3 {\killglue
\multcnt=#3\relax \hcor=0.97014\radius
                  \vcor=0.24253\radius
 \xup=#1\unit \yup=#2\unit
 \xxup=\wd\csname41\endcsname \yyup=\ht\csname41\endcsname
\whilenum \multcnt > 1 \do{\raise%
\yup\hbox to 0pt{\hskip \xup
 \rlap{\kern\hcor\raise\vcor\hbox{\ssssps}}\hss}\ifnum \multcnt = 2\advance\xup
 by \xxup
\advance\yup by \yyup \raise%
\yup\hbox to 0pt{\hskip \xup \rlap{\kern-\hcor\lower\vcor\hbox{\ssssps}}\hss}%
\advance\xup by -\xxup
\advance\yup by -\yyup \fi \advance\multcnt by -1 \advance\xup by
\xxup \advance \yup by \yyup}\ignorespaces}

\long\def\svilps #1 #2 #3 {\killglue
\multcnt=#3\relax \hcor=0.9863\radius
                  \vcor=0.16439\radius
 \xup=#1\unit \yup=#2\unit
 \xxup=\wd\csname61\endcsname \yyup=\ht\csname61\endcsname
\whilenum \multcnt > 1 \do{\raise%
\yup\hbox to 0pt{\hskip \xup
 \rlap{\kern\hcor\raise\vcor\hbox{\svips}}\hss}\ifnum \multcnt = 2\advance\xup
 by \xxup
\advance\yup by \yyup \raise%
\yup\hbox to 0pt{\hskip \xup \rlap{\kern-\hcor\lower\vcor\hbox{\svips}}\hss}%
\advance\xup by -\xxup
\advance\yup by -\yyup \fi \advance\multcnt by -1 \advance\xup by
\xxup \advance \yup by \yyup}\ignorespaces}

\long\def\svlpss #1 #2 #3 {\killglue
\multcnt=#3\relax \hcor=0.92847\radius
                  \vcor=0.37139\radius
 \xup=#1\unit \yup=#2\unit
 \xxup=\wd\csname52\endcsname \yyup=\ht\csname52\endcsname
\whilenum \multcnt > 1 \do{\raise%
\yup\hbox to 0pt{\hskip \xup
 \rlap{\kern\hcor\raise\vcor\hbox{\svpss}}\hss}\ifnum \multcnt = 2\advance\xup
 by \xxup
\advance\yup by \yyup \raise%
\yup\hbox to 0pt{\hskip \xup \rlap{\kern-\hcor\lower\vcor\hbox{\svpss}}\hss}%
\advance\xup by -\xxup
\advance\yup by -\yyup \fi \advance\multcnt by -1 \advance\xup by
\xxup \advance \yup by \yyup}\ignorespaces}

\long\def\ssslpsv #1 #2 #3 {\killglue
\multcnt=#3\relax \hcor=0.55148\radius
                  \vcor=0.8574\radius
 \xup=#1\unit \yup=#2\unit
 \xxup=\wd\csname35\endcsname \yyup=\ht\csname35\endcsname
\whilenum \multcnt > 1 \do{\raise%
\yup\hbox to 0pt{\hskip \xup
 \rlap{\kern\hcor\raise\vcor\hbox{\ssspsv}}\hss}\ifnum \multcnt = 2\advance\xup
 by \xxup
\advance\yup by \yyup \raise%
\yup\hbox to 0pt{\hskip \xup \rlap{\kern-\hcor\lower\vcor\hbox{\ssspsv}}\hss}%
\advance\xup by -\xxup
\advance\yup by -\yyup \fi \advance\multcnt by -1 \advance\xup by
\xxup \advance \yup by \yyup}\ignorespaces}

\long\def\svlpsss #1 #2 #3 {\killglue
\multcnt=#3\relax \vcor=0.55148\radius
                  \hcor=0.8574\radius
 \xup=#1\unit \yup=#2\unit
 \xxup=\wd\csname53\endcsname \yyup=\ht\csname53\endcsname
\whilenum \multcnt > 1 \do{\raise%
\yup\hbox to 0pt{\hskip \xup
 \rlap{\kern\hcor\raise\vcor\hbox{\svpsss}}\hss}\ifnum \multcnt = 2\advance\xup
 by \xxup
\advance\yup by \yyup \raise%
\yup\hbox to 0pt{\hskip \xup \rlap{\kern-\hcor\lower\vcor\hbox{\svpsss}}\hss}%
\advance\xup by -\xxup
\advance\yup by -\yyup \fi \advance\multcnt by -1 \advance\xup by
\xxup \advance \yup by \yyup}\ignorespaces}

\long\def\ssslps #1 #2 #3 {\killglue
\multcnt=#3\relax \hcor=0.94868\radius
                  \vcor=0.31622\radius
 \xup=#1\unit \yup=#2\unit
 \xxup=\wd\csname31\endcsname \yyup=\ht\csname31\endcsname
\whilenum \multcnt > 1 \do{\raise%
\yup\hbox to 0pt{\hskip \xup
 \rlap{\kern\hcor\raise\vcor\hbox{\sssps}}\hss}\ifnum \multcnt = 2\advance\xup
 by \xxup
\advance\yup by \yyup \raise%
\yup\hbox to 0pt{\hskip \xup \rlap{\kern-\hcor\lower\vcor\hbox{\sssps}}\hss}%
\advance\xup by -\xxup
\advance\yup by -\yyup \fi \advance\multcnt by -1 \advance\xup by
\xxup \advance \yup by \yyup}\ignorespaces}

\long\def\sssslpsss #1 #2 #3 {\killglue
\multcnt=#3\relax \hcor=0.8\radius
                  \vcor=0.6\radius
 \xup=#1\unit \yup=#2\unit
 \xxup=\wd\csname43\endcsname \yyup=\ht\csname43\endcsname
\whilenum \multcnt > 1 \do{\raise%
\yup\hbox to 0pt{\hskip \xup
 \rlap{\kern\hcor\raise\vcor\hbox{\sssspsss}}\hss}\ifnum \multcnt =
 2\advance\xup by \xxup
\advance\yup by \yyup \raise%
\yup\hbox to 0pt{\hskip \xup \rlap{\kern-\hcor\lower\vcor\hbox{\sssspsss}}\hss}%
\advance\xup by -\xxup
\advance\yup by -\yyup \fi \advance\multcnt by -1 \advance\xup by
\xxup \advance \yup by \yyup}\ignorespaces}

\long\def\sslpsss #1 #2 #3 {\killglue
\multcnt=#3\relax \vcor=0.8332\radius
                  \hcor=0.5547\radius
 \xup=#1\unit \yup=#2\unit
 \xxup=\wd\csname23\endcsname \yyup=\ht\csname23\endcsname
\whilenum \multcnt > 1 \do{\raise%
\yup\hbox to 0pt{\hskip \xup
 \rlap{\kern\hcor\raise\vcor\hbox{\sspsss}}\hss}\ifnum \multcnt = 2\advance\xup
 by \xxup
\advance\yup by \yyup \raise%
\yup\hbox to 0pt{\hskip \xup \rlap{\kern-\hcor\lower\vcor\hbox{\sspsss}}\hss}%
\advance\xup by -\xxup
\advance\yup by -\yyup \fi \advance\multcnt by -1 \advance\xup by
\xxup \advance \yup by \yyup}\ignorespaces}

\long\def\slpss #1 #2 #3 {\killglue
\multcnt=#3\relax \vcor=0.89442\radius
                 \hcor=0.4472\radius
 \xup=#1\unit \yup=#2\unit
 \xxup=\wd\csname82\endcsname \yyup=\ht\csname82\endcsname
\whilenum \multcnt > 1 \do{\raise%
\yup\hbox to 0pt{\hskip \xup \rlap{\kern\hcor\raise\vcor\hbox{\spss}}\hss}\ifnum
 \multcnt = 2\advance\xup by \xxup
\advance\yup by \yyup \raise%
\yup\hbox to 0pt{\hskip \xup \rlap{\kern-\hcor\lower\vcor\hbox{\spss}}\hss}%
\advance\xup by -\xxup
\advance\yup by -\yyup \fi \advance\multcnt by -1 \advance\xup by
\xxup \advance \yup by \yyup}\ignorespaces}

\long\def\slpsss #1 #2 #3 {\killglue
\multcnt=#3\relax \vcor=0.94868\radius
                  \hcor=0.31622\radius
 \xup=#1\unit \yup=#2\unit
 \xxup=\wd\csname83\endcsname \yyup=\ht\csname83\endcsname
\whilenum \multcnt > 1 \do{\raise%
\yup\hbox to 0pt{\hskip \xup
 \rlap{\kern\hcor\raise\vcor\hbox{\spsss}}\hss}\ifnum \multcnt = 2\advance\xup
 by \xxup
\advance\yup by \yyup \raise%
\yup\hbox to 0pt{\hskip \xup \rlap{\kern-\hcor\lower\vcor\hbox{\spsss}}\hss}%
\advance\xup by -\xxup
\advance\yup by -\yyup \fi \advance\multcnt by -1 \advance\xup by
\xxup \advance \yup by \yyup}\ignorespaces}

\long\def\svlmss #1 #2 #3 {\killglue
\multcnt=#3\relax \hcor=0.92847\radius
                  \vcor=0.37139\radius
 \xup=#1\unit \yup=#2\unit
 \xxup=\wd\csname52\endcsname \yyup=\ht\csname52\endcsname \advance\yup by -\yyup
\whilenum \multcnt > 1 \do{\raise%
\yup\hbox to 0pt{\hskip \xup
 \rlap{\kern\hcor\lower\vcor\hbox{\svmss}}\hss}\ifnum \multcnt = 2\advance\xup
 by \xxup
\advance\yup by -\yyup \raise%
\yup\hbox to 0pt{\hskip \xup \rlap{\kern-\hcor\raise\vcor\hbox{\svmss}}\hss}%
\advance\xup by -\xxup
\advance\yup by -\yyup \fi \advance\multcnt by -1 \advance\xup by
\xxup \advance \yup by -\yyup}\ignorespaces}

\long\def\ssslmsv #1 #2 #3 {\killglue
\multcnt=#3\relax \hcor=0.51449\radius
                  \vcor=0.8574\radius
 \xup=#1\unit \yup=#2\unit
 \xxup=\wd\csname35\endcsname \yyup=\ht\csname35\endcsname \advance\yup by -\yyup \whilenum \multcnt > 1 \do{\raise%
\yup\hbox to 0pt{\hskip \xup
 \rlap{\kern\hcor\lower\vcor\hbox{\sssmsv}}\hss}\ifnum \multcnt = 2\advance\xup
 by \xxup
\advance\yup by -\yyup \raise%
\yup\hbox to 0pt{\hskip \xup \rlap{\kern-\hcor\raise\vcor\hbox{\sssmsv}}\hss}%
\advance\xup by -\xxup
\advance\yup by -\yyup \fi \advance\multcnt by -1 \advance\xup by
\xxup \advance \yup by -\yyup}\ignorespaces}

\long\def\ssslmss #1 #2 #3 {\killglue
\multcnt=#3\relax \hcor=0.8332\radius
                  \vcor=0.5547\radius
 \xup=#1\unit \yup=#2\unit
 \xxup=\wd\csname32\endcsname \yyup=\ht\csname32\endcsname \advance\yup by -\yyup \whilenum \multcnt > 1 \do{\raise%
\yup\hbox to 0pt{\hskip \xup
 \rlap{\kern\hcor\lower\vcor\hbox{\sssmss}}\hss}\ifnum \multcnt = 2\advance\xup
 by \xxup
\advance\yup by -\yyup \raise%
\yup\hbox to 0pt{\hskip \xup \rlap{\kern-\hcor\raise\vcor\hbox{\sssmss}}\hss}%
\advance\xup by -\xxup
\advance\yup by -\yyup \fi \advance\multcnt by -1 \advance\xup by
\xxup \advance \yup by -\yyup}\ignorespaces}

\long\def\slmsss #1 #2 #3 {\killglue
\multcnt=#3\relax \vcor=0.94868\radius
                  \hcor=0.31622\radius
 \xup=#1\unit \yup=#2\unit
 \xxup=\wd\csname83\endcsname \yyup=\ht\csname83\endcsname \advance\yup by -\yyup \whilenum \multcnt > 1 \do{\raise%
\yup\hbox to 0pt{\hskip \xup
 \rlap{\kern\hcor\lower\vcor\hbox{\smsss}}\hss}\ifnum \multcnt = 2\advance\xup
 by \xxup
\advance\yup by -\yyup \raise%
\yup\hbox to 0pt{\hskip \xup \rlap{\kern-\hcor\raise\vcor\hbox{\smsss}}\hss}%
\advance\xup by -\xxup
\advance\yup by -\yyup \fi \advance\multcnt by -1 \advance\xup by
\xxup \advance \yup by -\yyup}\ignorespaces}

\long\def\ssslms #1 #2 #3 {\killglue
\multcnt=#3\relax \hcor=0.94868\radius
                  \vcor=0.31622\radius
 \xup=#1\unit \yup=#2\unit
 \xxup=\wd\csname31\endcsname \yyup=\ht\csname31\endcsname \advance\yup by -\yyup \whilenum \multcnt > 1 \do{\raise%
\yup\hbox to 0pt{\hskip \xup
 \rlap{\kern\hcor\lower\vcor\hbox{\sssms}}\hss}\ifnum \multcnt = 2\advance\xup
 by \xxup
\advance\yup by -\yyup \raise%
\yup\hbox to 0pt{\hskip \xup \rlap{\kern-\hcor\raise\vcor\hbox{\sssms}}\hss}%
\advance\xup by -\xxup
\advance\yup by -\yyup \fi \advance\multcnt by -1 \advance\xup by
\xxup \advance \yup by -\yyup}\ignorespaces}

\long\def\sssslmsss #1 #2 #3 {\killglue
\multcnt=#3\relax \hcor=0.8\radius
                  \vcor=0.6\radius
 \xup=#1\unit \yup=#2\unit
 \xxup=\wd\csname43\endcsname \yyup=\ht\csname43\endcsname \advance\yup by -\yyup \whilenum \multcnt > 1 \do{\raise%
\yup\hbox to 0pt{\hskip \xup
 \rlap{\kern\hcor\lower\vcor\hbox{\ssssmsss}}\hss}\ifnum \multcnt =
 2\advance\xup by \xxup
\advance\yup by -\yyup \raise%
\yup\hbox to 0pt{\hskip \xup \rlap{\kern-\hcor\raise\vcor\hbox{\ssssmsss}}\hss}%
\advance\xup by -\xxup
\advance\yup by -\yyup \fi \advance\multcnt by -1 \advance\xup by
\xxup \advance \yup by -\yyup}\ignorespaces}

\long\def\sslmsss #1 #2 #3 {\killglue
\multcnt=#3\relax \vcor=0.8332\radius
                  \hcor=0.5547\radius
 \xup=#1\unit \yup=#2\unit
 \xxup=\wd\csname23\endcsname \yyup=\ht\csname23\endcsname \advance\yup by -\yyup \whilenum \multcnt > 1 \do{\raise%
\yup\hbox to 0pt{\hskip \xup
 \rlap{\kern\hcor\lower\vcor\hbox{\ssmsss}}\hss}\ifnum \multcnt = 2\advance\xup
 by \xxup
\advance\yup by -\yyup \raise%
\yup\hbox to 0pt{\hskip \xup \rlap{\kern-\hcor\raise\vcor\hbox{\ssmsss}}\hss}%
\advance\xup by -\xxup
\advance\yup by -\yyup \fi \advance\multcnt by -1 \advance\xup by
\xxup \advance \yup by -\yyup}\ignorespaces}

\def\vm{\rlap{\lower\radius\hbox{\hskip -\thickness \vrule    
height\unit width\thickness depth0pt}}}

\setbox\csname77\endcsname=\hbox{\lline\char'077}
\def\vp{\rlap{\raise\radius\hbox{\hskip -0.5\thickness \vrule
height\unit width\thickness depth0pt}}}

\long\def\vline #1 #2 #3 {\killglue
\yup=#2\unit\xup=#1\unit\yyyup=#3\unit
\advance \yup by \radius
\advance \yyyup by -2\radius
\raise\yup\hbox to 0pt{\hskip \xup\hskip -0.5\thickness \vrule
height \yyyup width \thickness depth 0pt\hss}\ignorespaces}

\long\def\tvline #1 #2  {\killglue\xup=#1\unit \yup=#2\unit \raise\yup\hbox
to 0pt{\hskip #1\unit \vp\hss}\advance\yup by 5pt
\raise\yup\hbox to 0pt{\hskip \xup \vm\hss}\ignorespaces}

\long\def\tlinep #1 #2 {\hcor=0.7071\radius
                         \vcor=0.7071\radius
                      \putt #1 #2 {\hbox{\putt 0 0 {\nep}
                                         \putt 0.5 0.5 {\nem}
                                     }} }
\long\def\tlinem #1 #2 {\hcor=0.7071\radius
                         \vcor=0.7071\radius
                           \putt #1 #2 {\hbox{\putt -1 0 {\nwwp}
                                         \putt -1.5 0.5 {\nwwm}
                                     }} }
\long\def\linepi #1 #2 {\hcor=0.7071\radius
                         \vcor=0.7071\radius
                      \putt #1 #2 {\hbox{\putt 0 0 {\nep}
                                         \putt 0 0 {\nem}
                                     }} }
\long\def\linemi #1 #2 {\hcor=0.7071\radius
                         \vcor=0.7071\radius
                           \putt #1 #2 {\hbox{\putt -1 0 {\nwwp}
                                         \putt -1 0 {\nwwm}
                                     }} }

\long\def\linem #1 #2 #3 {\killglue\multcnt=#3\relax \xup=#1\unit \yup=#2\unit
 \whilenum \multcnt > 1 \do{\raise%
\yup\hbox to 0pt{\hskip \xup \nwp\hss}\ifnum \multcnt = 2%
\advance\xup by -\unit\advance\yup by \unit\raise%
\yup\hbox to 0pt{\hskip \xup \nwm\hss}\advance\xup by -\unit%
\advance\yup by -\unit\fi\advance\multcnt by -1 \advance\xup by%
-\unit \advance \yup by \unit}\ignorespaces}

\long\def\multiput #1 #2 #3 #4 #5 #6 {\killglue\multcnt=#5\relax
\xup=#1\unit \yup=#2\unit
\whilenum \multcnt > 0 \do
{\raise\yup\hbox to 0pt{\hskip
\xup #6\hss}\advance\multcnt by -1\advance\xup  by #3\unit\advance\yup by
 #4\unit}\ignorespaces}

\setbox\csname82\endcsname=\hbox{\lline\char'001}
\setbox\csname83\endcsname=\hbox{\lline\char'002}
\setbox\csname84\endcsname=\hbox{\lline\char'003}
\setbox\csname85\endcsname=\hbox{\lline\char'004}
\setbox\csname86\endcsname=\hbox{\lline\char'005}
\setbox\csname21\endcsname=\hbox{\lline\char'010}
\setbox\csname23\endcsname=\hbox{\lline\char'012}
\setbox\csname25\endcsname=\hbox{\lline\char'014}
\setbox\csname31\endcsname=\hbox{\lline\char'020}
\setbox\csname32\endcsname=\hbox{\lline\char'021}
\setbox\csname34\endcsname=\hbox{\lline\char'023}
\setbox\csname35\endcsname=\hbox{\lline\char'024}
\setbox\csname41\endcsname=\hbox{\lline\char'030}
\setbox\csname43\endcsname=\hbox{\lline\char'032}
\setbox\csname45\endcsname=\hbox{\lline\char'034}
\setbox\csname51\endcsname=\hbox{\lline\char'040}
\setbox\csname52\endcsname=\hbox{\lline\char'041}
\setbox\csname53\endcsname=\hbox{\lline\char'042}
\setbox\csname54\endcsname=\hbox{\lline\char'043}
\setbox\csname56\endcsname=\hbox{\lline\char'045}
\setbox\csname61\endcsname=\hbox{\lline\char'050}
\setbox\csname65\endcsname=\hbox{\lline\char'054}

\setbox\csname112\endcsname=\hbox{\lline\char'101}
\setbox\csname113\endcsname=\hbox{\lline\char'102}
\setbox\csname114\endcsname=\hbox{\lline\char'103}
\setbox\csname115\endcsname=\hbox{\lline\char'104}
\setbox\csname116\endcsname=\hbox{\lline\char'105}
\setbox\csname121\endcsname=\hbox{\lline\char'110}
\setbox\csname123\endcsname=\hbox{\lline\char'112}
\setbox\csname125\endcsname=\hbox{\lline\char'114}
\setbox\csname131\endcsname=\hbox{\lline\char'120}
\setbox\csname132\endcsname=\hbox{\lline\char'121}
\setbox\csname134\endcsname=\hbox{\lline\char'123}
\setbox\csname135\endcsname=\hbox{\lline\char'124}
\setbox\csname141\endcsname=\hbox{\lline\char'130}
\setbox\csname143\endcsname=\hbox{\lline\char'132}
\setbox\csname145\endcsname=\hbox{\lline\char'134}
\setbox\csname151\endcsname=\hbox{\lline\char'140}
\setbox\csname152\endcsname=\hbox{\lline\char'141}
\setbox\csname153\endcsname=\hbox{\lline\char'142}
\setbox\csname154\endcsname=\hbox{\lline\char'143}
\setbox\csname156\endcsname=\hbox{\lline\char'145}
\setbox\csname161\endcsname=\hbox{\lline\char'150}
\setbox\csname165\endcsname=\hbox{\lline\char'154}

\newcommand{\cat}[1]{\boldsymbol{\mathscr{#1}}}

\newcommand{\CPb}{\cat P_{\kern -2pt\scriptscriptstyle 01}}

\renewcommand{\leq}{\leqslant}         
\renewcommand{\geq}{\geqslant}

\newcommand{\twiddle}[1]{{\smash{\underset{\raise.375ex\hbox{$\smash\sim$}}     
       {#1}}\vphantom{\underline{#1}}}} %

\renewcommand{\nleq}{\nleqslant}   
\renewcommand{\ngeq}{\ngeqslant}

\spnewtheorem*{acknw}{Acknowledgments}{\bfseries}{\rmfamily}   
   
\begin{document}

\title{A view of canonical extension}

\author{Mai Gehrke\inst{1} and Jacob Vosmaer\inst{2}}   
\institute%
{IMAPP\\
Radboud Universiteit Nijmegen\\
The Netherlands\\
\email{%
M.Gehrke@math.ru.nl}     
\and   
ILLC\\
Universiteit van Amsterdam\\
The Netherlands\\%
\email{%
J.Vosmaer@uva.nl}  
  }

\maketitle   

\begin{abstract} 
This is a short survey illustrating some of the essential aspects of the
theory of canonical extensions.  In addition some topological
results about canonical extensions of lattices with additional operations in 
finitely generated varieties are given. In particular, they are doubly algebraic 
lattices and their interval topologies agree with their double Scott topologies 
and make them Priestley topological algebras.
\end{abstract}     
     

\begin{acknw}
  Both authors would like to acknowledge the influence of discussions
  and work with H.A. Priestley on the content of this paper. In
  particular, the fact that something like Theorem~\ref{thm:fgv} might
  hold was first discussed with H. A. Priestley in the process of the
  first author's work on a book-in-preparation on Lattices in
  Logic. The expository material in this paper is also based to a
  large extent on work on the book.  Furthermore, the first author
  acknowledges support from and would like to thank the Programme {\it
    Research in Paris}. The work of the second author was made
  possible by VICI grant
   639.073.501 of the Netherlands Organization for Scientific Research
   (NWO).
\end{acknw}

\section{Introduction}\label{sec:intro}

Associating algebraic models to propositional logics is often achieved by an 
easy transcription of the syntactic
 specifications of such logics. 
This may be through 
the associated Linden\-baum\---Tarski algebras or through a transcription of a 
Gentzen-style calculus. As a consequence, semantic modelling by such algebras 
is often not far removed from the syntactic treatment of the logics. Relational 
semantics on the other hand, when available, are likely to give a significantly different 
and much more powerful tool. This phenomenon is akin to that whereby algebraists 
have exploited topological dualities to great advantage. One twist in the logic setting 
is that the topology doesn't have a natural place in the logic landscape, thus prompting 
logicians simply to discard it. As a consequence we obtain, from an algebra of formulas, 
a topo-relational space and then, by forgetting the topology, simply a relational structure.
The complex algebra of this structure is then an algebra based on a powerset in which 
the original formula algebra embeds. This is the so-called canonical extension. It turns 
out that it is abstractly characterised by three simple properties of the way it extends the 
original algebra and that it is in fact a very natural completion of the algebra. As such it 
provides a tool for studying Stone duality and its generalisations in an algebraic setting 
which is particularly well-suited for the treatment of additional algebraic structure on the 
underlying lattices or Boolean algebras.

The study of canonical  extensions originated in the famous paper of B.~J\'{o}nsson and 
A.~Tarski  \cite{JT51} on Boolean algebras with operators (BAOs).  Amongst BAOs are 
the modal algebras which supply semantic models for modal logics. The theory has since 
been generalised and simplified and now includes the algebraic counterparts of positive 
versions of modal logic, as well as intuitionistically based logics and substructural logics. 
Canonicity, that is, the property of being closed under canonical extension, for a class of 
algebraic models associated with a logic, yields complete relational semantics for the logic 
and even in the absence of canonicity canonical extensions, just like topological duality, 
provide a powerful tool for studying a logic.

This short survey, which corresponds to three tutorial lectures by the
first author in Bakuriani in Fall 2009, 
is based on materials Hilary Priestley and the first author are preparing for our book in the Oxford 
University Press Logic Guides series on Lattices in Logic: duality, correspondence, and 
canonicity. The three lectures focused in turn on: the relationship of canonical extension 
to topological duality and to questions of relational semantics of logics; the flavour and form 
of the basic theory of canonical extensions in their own right; topological methods in the 
theory of canonical extensions. This survey follows the same pattern with the addition of
a section on finitely generated varieties of lattices with additional operations to illustrate
the theory.

Accordingly, in Section~\ref{sec:LatLog} we identify the connection between questions 
about relational semantics in logic, topological duality, and canonical extension. In 
particular, we show that topological duality gives rise to a completion satisfying the 
properties which are the defining properties of canonical extension. In 
Section~\ref{sec:canext} we give the abstract definition of canonical extensions of 
arbitrary lattices. We give a few examples and outline how the abstract properties of 
canonical extensions uniquely determine them thereby actually deriving an alternate 
way of building canonical extensions which does not depend on the axiom of choice.
In Section \ref{sec:addop} we consider additional operations on lattices introducing the topological
approach. We give a few new results on the interaction of the lifting of maps to canonical 
extensions and topological properties of the maps. In the final section we study finitely 
generated lattice varieties. We show that canonical extensions of lattices lying in finitely 
generated lattice varieties are doubly algebraic lattices that are Stone spaces in their 
Scott and dual Scott topologies. We also show that canonical extension is functorial 
on all finitely generated varieties of lattice-based algebras and that the canonical extensions 
are Stone topological algebras in their double Scott topologies.

\section{Canonical extension, duality, and relational semantics}\label{sec:LatLog}

A propositional logic is typically specified by a consequence relation on the formulas 
or compound propositions of the logic. That is, the connectives and their arities are
specified and a set of primitive propositional variables is chosen. The formulas are
then defined inductively by proper application of the connectives. This already is closely
related to algebra as the formulas form the absolutely free algebra in the type of the
connectives over the set of variables. In the syntactic specification of a logic, a calculus
is then given for generating the consequence relation. In many cases this calculus 
corresponds to quotienting the free algebra by an equational theory and thus results
in a free algebra of a variety. For example, classical propositional logic corresponds 
to the variety of Boolean algebras, intuitionistic propositional logic to Heyting algebras, 
modal logic to modal algebras, and the Lambek calculus to ordered residuated
monoids.

In contrast, semantic conceptions of logic are based on some notion of
models and interpretations of the formulas in these. Thus models of classical logic are valuations
specifying the truth of the primitive propositions, and models of modal logics are 
evaluations on Kripke structures. These are objects of a different nature than formulas 
and their quotient algebras. This fundamental difference of sorts becomes very clear 
when considering the meaning of syntax and semantics in applications: in computer 
science applications, formulas and their logical calculi model specification of programs 
whereas their semantics model state-based transition systems. Lines of code and states 
of a machine are objects of completely different physical sorts. A fundamental question 
then is how we can identify the corresponding sort when we are given only one or the other. 
That is, given a syntactic specification, what is the corresponding semantics and vice versa? 
Going from semantics to syntax may be seen as a significant goal of coalgebraic logic. In the 
other direction, mathematics provides a useful tool in the form of topological duality theory. 

Topological duality theory is a fundamental tool of mathematics that allows one to connect 
theories or completely different sorts, namely algebra and topology. Topological duality, 
pioneered by Stone, is central in functional analysis, in algebraic geometry, and also in logic 
where it is at the base of most completeness results (in predicate logic as well as in 
propositional logic). It allows one to build a dual space from a lattice.  In logic
applications, the lattice operations are typically present as they model (some aspect of)
conjunction and disjunction but there are usually other connectives as well. Extended 
Stone or Priestley duality is tailored to this setting. For example, the dual space of a 
Boolean algebra is a Boolean space, that is, a compact $0$-dimensional space, while
the dual of a modal algebra is a Boolean space equipped with a binary relation whose point 
images are closed and such that the inverse image of each clopen is clopen (known as
descriptive general frames). In general in extended duality, distributive lattices with additional 
operations having some property of preserving or reversing joins or meets correspond dually 
to topo-relational spaces where the additional relations are of arity one higher than the arity of 
the corresponding operations and have some basic topological properties. 

While this correspondence provided by extended duality is pertinent, one fundamental 
difficulty in logic and computer science applications is how to understand and deal with the 
topology. There are essentially two solutions to this problem:
\begin{itemize}
\item Simply discard and forget the topology; this is, for example, the approach in modal logic.
\item Restrict to a setting where the topology is fully determined by a first order structure;
this is the case in domain theory where dual spaces carry the Scott topology which is fully 
determined by an order.
\end{itemize}
The second setting recognises topology as having meaning, namely in the form of observability, 
but both raise questions about duality: how it relates to the discrete duality and when a poset 
is spectral in its Scott topology. We will touch on both of these in this article. 

Canonical extension is most obviously related to the first approach of forgetting the 
topology but it is in fact a way, in general, of studying duality in an algebraic setting. 
This is useful not only for forgetting the topology but also for studying additional 
algebraic structure, that is, extended duality and for identifying algebraic settings where 
the topology is order determined.
 
As mentioned above, at its origin, canonical extension is an algebraic way of talking 
about extended Stone duality for Boolean algebras with operators. We illustrate this
with the case of modal algebras \cite{BDV}. The pertinent square is the following.
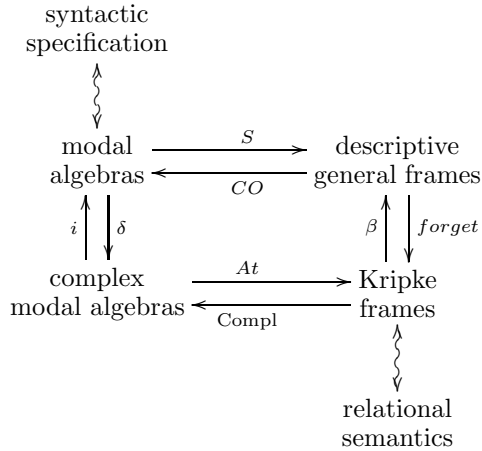
\[ \xymatrix@C=10ex{ 
\txt{syntactic\\specification} \ar@{<~>}[d] &\\
\txt{modal\\algebras} \ar@<1ex>[d]^{\delta}
  \ar@<1ex>[r]^{S}&
  \txt{descriptive\\general frames} \ar@<1ex>[l]^{CO} \ar@<1ex>[d]^{forget}\\
  \txt{complex\\modal algebras} \ar@<1ex>[u]^{i}\ar@<1ex>[r]^{At}&
  \txt{Kripke\\frames} \ar@<1ex>[u]^{\beta}\ar@<1ex>[l]^{\rm Compl} \\
&\txt{relational\\semantics} \ar@{<~>}[u]
 }\] 
Here the upper pair of functors gives the extended Stone duality for modal algebras. 
The dual of a modal algebra is a descriptive general frame $(X,\tau,R)$ and 
forgetting the topology yields a Kripke frame $(X,R)$. Kripke frames also lie in the 
scope of a duality namely the `discrete' duality with complete and atomic Boolean 
algebras with a completely join preserving diamond. The canonical extension is 
obtained concretely by walking around the square from upper left to lower left corner.
That is, given a modal algebra, $A$, we take its dual general descriptive frame,
$(X,\tau, R)$, forget the topology to get the Kripke frame $(X,R)$, and then we 
form the complex algebra, ${\rm Compl}(X,R)=({\mathcal P}(X), \Diamond_R)$ where 
$R$ and $\Diamond_R$ are related by 
\begin{equation}\label{R-def}
\forall x,y\in X \qquad R(x,y)\ \iff\ x\leq\Diamond_R(y).
\end{equation}
Here we identify atoms of $\mathcal{P}(X)$ with elements of $X$. 
The fact that extended topological duality is a duality includes the fact that the original 
modal algebra is isomorphic to the modal algebra of clopen subsets of $(X,\tau,R)$
with the restriction of the operation $\Diamond_R$. Thus we have, for each modal algebra, 
an embedding $A\hookrightarrow {\rm Compl}(X,R)=A^\delta$; this embedding is a concrete incarnation of
what is known as the \emph{canonical extension}. It is clear that to study what happens when 
we `forget the topology', the canonical extension is central. However, what makes the 
canonical extension of general interest are the following two facts:
\begin{itemize}
\item The canonical extension may be abstractly characterised as a certain completion of 
$A$ in a purely complete lattice theoretic setting;
\item We can construct the dual space of $A$ from the canonical extension 
$A\hookrightarrow A^\delta$.
\end{itemize}
This is why we can claim that the theory of canonical extensions may be seen as an algebraic 
formulation of Stone/Priestley duality. 

The second of the two above facts is easy to see: Suppose we have somehow been supplied with 
$\widehat{\ }: A\hookrightarrow A^\delta$, how can we 
reconstruct $X$, $R$, and $\tau$ from this algebraic information? First we apply discrete duality to 
$A^\delta$. That is, we recover $X$ as the atoms of $A^\delta$ and we recover $R$
by using (\ref{R-def}). The topology is generated by the `shadows' of the elements of $A$ on the set $X$, 
that is, by the sets ${\downarrow}\hat{a}\cap At(A^\delta)=\{x\in At(A^\delta)\mid x\leq \widehat{a}\}$ 
where $a$ ranges over $A$.

The abstract characterisation of the embedding $\widehat{\ }: A\hookrightarrow {\rm Compl}(X,R)$ is obtained in 
two tempi. First for the underlying lattice and then for the additional operations. We will return to the additional
operations in Section~\ref{sec:addop} where we see they are natural upper- or lower-semicontinuous envelopes. 
We conclude this section by proving the three properties of 
$\widehat{\ }: A\hookrightarrow {\mathcal P}(X)$ which are used in the abstract definition of canonical extension.
To this end, let $A$ be a Boolean algebra. The Stone space of $A$ is given by 
\begin{align*}
X&=\{\mu\subseteq A\mid \mu\mbox{ is an ultrafilter of }A\} \mbox{ is the set underlying the space,}\\
{\mathcal B}\ &=\{\widehat{a}\mid a\in A\} \mbox{ is a basis for the topology where 
                                                    $\widehat{a}=\{\mu\mid a\in\mu\}$ for $a\in A$.}
\end{align*}
The fundamental result needed to derive properties of dual spaces is Stone's Prime Filter Theorem:
If a filter $F$ and an ideal $I$ of a Boolean algebra $A$ are disjoint then there exists an ultrafilter
$\mu$ of $A$ containing $F$ and disjoint from $I$. 
 Here we use the
fact that since $A$ is a Boolean algebra, $F\subseteq A$ is an
ultrafilter iff it is a prime filter. 
We prove the following three propositions.

\begin{proposition}\label{canextcompl}
Let $A$ be a Boolean algebra and $X$ the dual space of $A$. Then the map
\begin{align*}
\widehat{\ }\ :A&\to{\mathcal P}(X)\\
                   a&\mapsto \widehat{a}=\{\mu\mid a\in\mu\}
\end{align*}
is a lattice completion of $A$.
\end{proposition}

\proof
It is clear that ${\mathcal P}(X)$ is a complete lattice. We have to show that the map $\widehat{\ }$
is a lattice embedding. Since ultrafilters are upsets, it is clear that $\widehat{\ }$ is order preserving.
Thus $\widehat{a\wedge b}\subseteq\widehat{a}\cap\widehat{b}$ and $\widehat{a}\cup\widehat{b}
\subseteq\widehat{a\vee b}$. Also, if $\mu\in \widehat{a}\cap\widehat{b}$ then $a\in\mu$ and 
$b\in\mu$ and thus $a\wedge b\in\mu$ since filters are closed under meet. For the join preservation
note that $\mu\in \widehat{a\vee b}$ implies that $a\vee b\in\mu$ and since ultrafilters are prime 
filters, it follows that $a\in\mu$ or $b\in\mu$ and thus, in either case, $\mu\in\widehat{a}\cup\widehat{b}$.
Finally, if $a,b\in A$ with $a\neq b$ then either $a\nleq b$ or
$b\nleq a$. The former implies that 
 the filter
$F={\uparrow}a$ and the ideal $I={\downarrow}b$ are disjoint. 
Thus there is a $\mu\in X$ with $F\subseteq\mu$ and $I$ disjoint from $\mu$. That is, $\mu\in\widehat{a}$  
but $\mu\not\in\widehat{b}$. By symmetry the same thing happens if  $b\nleq a$.
\qed\smallskip

\begin{proposition}\label{canextdense}
Let $A$ be a Boolean algebra and $X$ the dual space of $A$. Then the image of the map \ 
$\widehat{\ }:A\to{\mathcal P}(X)$ given by $a\mapsto \widehat{a}=\{\mu\mid a\in\mu\}$ is $\bigvee\bigwedge$- 
and $\bigwedge\bigvee$-dense in ${\mathcal P}(X)$. That is, every element of ${\mathcal P}(X)$ is both an 
intersection of unions and a union of intersections of elements in the image of $\ \widehat{\ }$.
\end{proposition}

\proof
This is easily seen by noting that each subset of ${\mathcal P}(X)$ is a union of singletons and for each 
singleton $\{\mu\}$ we have $\{\mu\}=\bigcap\{\widehat{a}\mid
a\in\mu\}$. The rest follows by order duality, 
 using De Morgan's laws.
\qed\smallskip\smallskip

\begin{proposition}\label{canextcompact}
Let $A$ be a Boolean algebra and $X$ the dual space of $A$. The map \ $\widehat{\ }:A\to{\mathcal P}(X)$
given by $a\mapsto \widehat{a}=\{\mu\mid a\in\mu\}$ is such that for any subsets $S$ and $T$ of $A$ with 
$\bigcap\{\widehat{s}\mid s\in S\} \subseteq \bigcup\{\widehat{t}\mid t\in T\}$, there exist finite sets 
$S' \subseteq S$ and $T' \subseteq T$ such that $\bigwedge S' \leq \bigvee T'$ in $A$.\\
\end{proposition}

\proof
This is a straight forward consequence of Stone's Prime Filter Theorem: If the conclusion is false, then
the filter generated by $S$ is disjoint from the ideal generated by $T$ and it follows that there is a prime
filter $\mu\in X$ containing the filter and disjoint from the ideal. It follows that 
$\mu\in\widehat{s}$ for each $s\in S$ but $\mu\not\in\widehat{t}$ for any $t\in T$ thus violating the 
antecedent of the statement of the proposition.
\qed\smallskip

\section{Working with canonical extensions} 
\label{sec:canext}

\noindent 
The philosophy of the canonical extension approach, since its first introduction by J{\'o}nsson 
and Tarski, and its real power, come from the fact that one can work with it abstractly 
{\it without referring to the particular way the canonical extension has been built}, using only a 
few very simple properties, namely what we will call completeness,
compactness, and density. We work in the setting of arbitrary bounded lattices.

\begin{definition}  {\rm(canonical extension)} 
Let $L$ be a lattice. A {\it  canonical extension of}~$L$ is a lattice completion  
$L\hookrightarrow L^\delta $ of $L$ with the following two properties:\\

\noindent{\bf Density:}\  The image of $L$ is $\bigvee\bigwedge$- and $\bigwedge\bigvee$-dense in 
$L^\delta $, that is, every element of $L^\delta $ is both a join of meets and a 
meet of joins of elements from $L$;\\

\noindent{\bf Compactness:} \  Given  any subsets $S$ and $T$ of $L$ with 
$\bigwedge S \leq \bigvee T$ in $L^\delta $, there exist finite sets $S' \subseteq S$ 
and $T' \subseteq T$ such that $\bigwedge S' \leq \bigvee T'$.\\
\end{definition}

The following equivalent formulations of compactness are often useful and are not hard to prove.

\begin{proposition}  {\rm(variants of compactness)}
\label{prop:compactness}
Let $L$ be a lattice and $L'$ a complete lattice. Each of the following conditions on
an embedding $L\hookrightarrow L'$ is equivalent to the compactness
property:\\
\begin{tabular}{cp{11cm}}
(C')& Given  any down-directed subset $S$ of $L$ and any up-directed subset $T$ of $L$ 
with $\bigwedge S \leq \bigvee T$ in $L'$, there exist $s\in S$ and $t\in T$ such that 
$s \leq t$.\\
(C") &Given  any filter $F$ of $L$ and any ideal $I$ of $L$ with $\bigwedge F \leq \bigvee I$ 
in $L'$, we have $F\cap I\ne\emptyset$.
\end{tabular}
\end{proposition}

First we consider a few examples.

\begin{example}{\rm (lattices that are their own canonical extension)}  \label{canext-ex1}  
Let $L$ be a  finite lattice, or more generally a bounded lattice with no 
infinite chains. We claim that the identity $L\hookrightarrow L$ is a 
canonical extension of $L$. This is a completion of $L$ because a 
bounded lattice with no infinite chains is automatically complete; see, 
for example, \cite{ILO2}, Theorem~2.41.
We remark that  a poset has no infinite chains if and only if it satisfies
both (ACC) and (DCC) (sufficiency requires the axiom of choice) and that 
the reason that this forces  completeness of a bounded lattice
is because, in the presence of (ACC), arbitrary non-empty joins 
reduce to finite joins, and dually; more details can be found in \cite{ILO2}; see
Lemma~2.39 and Theorem~2.40.  
It is of course clear that the identity is a dense embedding, and compactness 
follows because every join and meet reduces to finite ones in a lattice with
(ACC) and (DCC) as remarked above. 
We note that the converse is also true. Suppose $L \hookrightarrow L$ is a canonical 
extension and $C\subseteq L$ is a chain in $L$. Then $a=\bigvee C\in L$ must 
exist (since $L$ must be complete), and by compactness, there must be 
$c_1,\ldots, c_n\in L$ with $a\leq c_1\vee\ldots\vee c_n$. Since $C$ is a chain,
this implies there is an $i\in\{1,\ldots,n\}$ with $c_1\vee\ldots\vee c_n=c_i$ and 
thus $a=c_i$ and $L$ satisfies (ACC). 
If the identity on $L$ is a canonical extension
then the same is true for the dual lattice. Thus, by order duality, $L$ also satisfies (DCC)
and thus $L$ has no infinite chains.\end{example}

\begin{example}{\rm (canonical extensions of chains).}  \label{canext-ex2}  
As our next  example we consider the infinite chain $L= \omega \oplus \omega ^\partial$, where 
$P^\partial$ denotes the order dual of a poset $P$. This lattice $L$, which is shown in Fig.~\ref{chain-ex}, 
is the reduct of the MV-chain known as the Chang algebra. 
\begin{figure}
$$
{\hbox{
\putt -10 0 {\hbox{
  \putt -2 0 {\hbox{
  \point 0 0
  \point 0 2
  \point 0 4
  \point 0 14
  \point 0 16
  \point 0 18
  \labelr 0 0 {c_0}
  \labelr 0 2 {c_1}
  \labelr 0 4 {c_2}
  \labelr 0 14 {b_2}
  \labelr 0 16 {b_1}
  \labelr 0 18 {b_0} 
  \vline 0 0 2
  \vline 0 2 2
  \putt 0 5.2 {\hbox{
    \putt 0 0 {$\cdot$}
    \putt 0 .5 {$\cdot$}
    \putt 0 1 {$\cdot$}
     }}
\putt 0 11.2 {\hbox{ 
     \putt 0 0 {$\cdot$}
     \putt 0 .5 {$\cdot$}
     \putt 0 1 {$\cdot$}
     }}
  \vline 0 14 2
  \vline 0 16 2
}}
\putt -3 -1.5  {$L$} 
}}
\putt 0 0 {\hbox{
  \putt -2 0 {\hbox{
  \point 0 0
  \point 0 2
  \point 0 4
  \spoint 0 9
  \point 0 14
  \point 0 16
  \point 0 18
  \labelr 0 0 {c_0}
  \labelr 0 2 {c_1}
  \labelr 0 4 {c_2}
  \labelr 0 9 {z}
  \labelr 0 14 {b_2}
  \labelr 0 16 {b_1}
  \labelr 0 18 {b_0}
  \vline 0 0 2
  \vline 0 2 2
  \putt 0 5.2 {\hbox{
    \putt 0 0 {$\cdot$}
    \putt 0 .5 {$\cdot$}
    \putt 0 1 {$\cdot$}
     }}
\putt 0 11.2 {\hbox{ 
     \putt 0 0 {$\cdot$}
     \putt 0 .5 {$\cdot$}
     \putt 0 1 {$\cdot$}
     }}
  \vline 0 14 2
  \vline 0 16 2
}}
\putt -3 -1.5  {$\overline{L} $}
}}      
\putt 10 0 {\hbox{
  \putt -2 0 {\hbox{
  \point 0 0
  \point 0 2
  \point 0 4
  \spoint 0 8
  \spoint 0 10
  \point 0 14
  \point 0 16
  \point 0 18
  \labelr 0 0 {c_0}
  \labelr 0 2 {c_1}
  \labelr 0 4 {c_2}
  \labelr 0 8 {y=c_\infty}
  \labelr 0 10 {x = b_\infty}
  \labelr 0 14 {b_2}
  \labelr 0 16 {b_1}
  \labelr 0 18 {b_0}
  \vline 0 0 2
  \vline 0 2 2
  \putt 0 5.2 {\hbox{
    \putt 0 0 {$\cdot$}
    \putt 0 .5 {$\cdot$}
    \putt 0 1 {$\cdot$}
     }}
\putt 0 11.2 {\hbox{ 
     \putt 0 0 {$\cdot$}
     \putt 0 .5 {$\cdot$}
     \putt 0 1 {$\cdot$}
     }}
  \vline 0 8 2
  \vline 0 14 2
  \vline 0 16 2
}}
\putt -3 -1.5  {$L^\delta $}
}}                      
}}
$$ 
\caption{}\label{chain-ex}
\end{figure}
We claim that the embedding of $L$ as a subposet of 
the lattice $L^\delta$ as depicted in the figure is a canonical extension
of $L$ but that the embedding of $L$ as a subposet of $\overline{L}$ 
is not. It is clear that both $\overline{L}$ and $L^\delta$ are complete 
(while $L$ is not). Thus the inclusions  $L\hookrightarrow\overline{L}$ 
and $L\hookrightarrow L^\delta$ are both completions of $L$. Further 
it is easy to see that both satisfy the density condition. However, 
$L\hookrightarrow\overline{L}$ is not compact since 
$$\bigwedge_{i=1}^\infty b_i=z\leq z=\bigvee_{i=1}^\infty c_i$$
but no finite meet of $b_i$s gets below a finite join of $c_j$s. It is easy to 
convince oneself that the embedding $L\hookrightarrow L^\delta$ is 
compact. We note that $L\hookrightarrow\overline{L}$ is the MacNeille
completion, i.e., the unique completion of $L$ with the stronger
density property that every element of the completion is obtained 
both as a join of elements from $L$ and as a meet of elements 
from $L$.   
\end{example}

\begin{example}{ \rm (Classical propositional logic example)}  \label{canext-cpl} 
Let $L$ denote the Lindenbaum-Tarski algebra of classical propositional logic, or 
equivalently the free Boolean algebra, on the countable set of variables
$X=\{x_1,x_2,\ldots\}$. Also, let $L_n$ be the classical propositional 
logic on the set $X_n=\{x_1,x_2,\ldots, x_n\}$. It is well known that for 
each $n$ we have $L_n\cong2^{2^{X_n}}$. For infinitely many variables
this is not so, however, we will see that the canonical extension of 
$L$ is the algebra $2^{2^X}$. More precisely, we show that the Boolean 
homomorphism uniquely determined by the freeness of $L$ over $X$ 
and the assignment
\begin{align*}
 L & \overset e \hookrightarrow 2^{2^X}\\
        x_i  & \mapsto \{\,\alpha\in2^X\mid x_i\in\alpha\,\}
\end{align*}
is a canonical extension of $L$. By the very definition of $e$ it is a Boolean
homomorphism. Note that in the finite case
\begin{align*}
e_n:L_n&\to 2^{2^{X_n}}\\
x_i&\mapsto~\{\,\alpha\in2^{X_n}\mid x_i\in\alpha\,\}
\end{align*}
is the standard isomorphism showing that $L_n\cong~2^{2^{X_n}}$. 
For each 
two formulas $\phi$ and $\psi$ there is an $n$ so that $\phi,\psi\in L_n$ 
and for $\phi\in L_n$ we have $e(\phi)\cap2^{X_n}=e_n(\phi)$ and thus $e$ is 
an injection since the $e_n$ are. Thus $e$ is an embedding.
\newline
Next we show that $e$ satisfies the density condition. Since we are dealing 
with Boolean algebras and the embedding preserves negation, it is enough 
to show that every element of $2^{2^X}$ may be obtained as a join of meets of 
elements in the image of $e$. Thinking of $2^{2^X}$ as ${\mathcal P}({\mathcal P}(X))$,   
it suffices to show that $\{\alpha\}$ may be obtained as an intersection of sets in the 
image of $e$ for each $\alpha\in {\mathcal P}(X)$. For $\alpha\in {\mathcal P}(X)$ let 
$$
\phi_n=(\bigwedge [X_n\cap\alpha])\wedge(\bigwedge\{\,\neg x\mid x\in X_n\setminus\alpha\,\})
$$
where `$\setminus$'  denotes the difference of sets, it is then easy to see that
$$
\bigcap_{n=1}^\infty e(\phi_n)=\{\alpha\}.
$$
\newline
Finally we show that $e$ is a compact embedding. Let $S,T\subseteq L$ with 
$\bigcap e(S)\subseteq\bigcup e(T)$. Since we are in a power 
set and $e$ preserves complements, we can rewrite this as 
${\mathcal P}(X)=\bigcup e(\neg S\cup T)$ where $\neg S=\{\,\neg \phi\mid \phi\in S\,\}$. 
Thus we just need to verify the usual topological compactness property. To this end let 
$T$ be any subset of $L$ with ${\mathcal P}(X)=\bigcup e(T)$.
but assume that no finite subcover of ${\mathcal C}=e(T)$ covers ${\mathcal P}(X)$. 
Since each $\phi$ in $T$ may be written as a disjunction of conjunctions of 
literals, %
we may assume without 
loss of generality that each $\phi\in T$ is a conjunction of literals. We define 
a sequence of literals inductively. Let $l_1=x_1$ if $e(x_1)$ cannot be covered 
by a finite subcover of ${\mathcal C}$, otherwise let $l_1=\neg x_1$. 
Note that if both $e(x_1)$ and $e(\neg x_1)$ can be covered by finite subcovers 
of ${\mathcal C}$ then so can ${\mathcal P}(X)$. Thus $l_1$ is not covered by a finite 
subcover of ${\mathcal C}$. For each $n\geq 1$, if $l_1,\ldots, l_n$ have been defined, 
we define $l_{n+1}=l_1\wedge\ldots l_n\wedge x_{n+1}$ if 
$e(l_1\wedge\ldots l_n\wedge x_{n+1})$ cannot be covered by a finite subcover 
of ${\mathcal C}$ and we let $l_{n+1}=l_1\wedge\ldots l_n\wedge\neg x_{n+1}$ 
otherwise. From our assumption, it is not hard to prove by induction on $n$ that 
$e(\bigwedge_{i=1}^n l_i)$ cannot be covered by a finite subcover of ${\mathcal C}$. 
Now let $\alpha=\{\,x_i\mid l_i=x_i\,\}$. Since ${\mathcal C}$ covers ${\mathcal P}(X)$, there 
is some $\phi\in T$ with $\alpha\in e(\phi)$
and thus
$\phi=\bigwedge_{i\in I}l_i$ where $I$ is a finite set of natural numbers. If $I=\emptyset$,
then $\phi=1$ and $e(\phi)={\mathcal P}(X)$ is a singleton subcover of $\mathcal C$. Since we are
assuming no such cover exists, $I\ne\emptyset$. Now let $n=\max(I)$, then 
$\bigwedge_{i=1}^n l_i\leq\phi$ and thus $e(\bigwedge_{i=1}^n l_i)$
is covered by $e(\phi)$ which is a contradiction. 
We conclude that ${\mathcal C}$ must contain a finite subcover of ${\mathcal P}(X)$
thus proving compactness.
\newline
We note that this illustrative example is just a special case of the fact that the canonical extension 
of any Boolean algebra is given by the Stone embedding into the power set of its dual space.
\end{example}

Next we outline the development leading to the uniqueness and existence in general 
of canonical extensions of lattices. The density condition that is part of the abstract 
definition of canonical extension makes it clear that the meet and the join closure of 
$L$ in $L^\delta$ play a central role.

\begin{definition} {\rm(filter and ideal elements)}\ 
Let $L$ be a lattice,  and $L^\delta$ a canonical extension of $L$.
Define
\begin{align*}
      F(L^\delta )&:= \{ \, x \in L^\delta \mid x \text{ is a meet of elements from } L\,\}, \\ \\
      I(L^\delta ) &:= \{ \, y \in L^\delta \mid y \text{ is a join of elements from } L\,\}.   
\end{align*}
We refer to the elements of $F(L^\delta )$ as {\it filter elements} and to the elements
of $I(L^\delta )$ as {\it ideal elements}.
\end{definition}

The rationale for naming these elements filter and ideal elements, respectively, is 
made clear by the following proposition. 

\begin{proposition}
Let $L$ be a lattice,  and $L^\delta$ a canonical extension of $L$. Then 
the poset $F(L^\delta )$ of filter elements of $L^\delta$ is reverse order
isomorphic to the poset ${\text{\bf Filt}}(L)$  of lattice filters of $L$ and the poset 
$I(L^\delta )$ of ideal elements of $L^\delta$ is order isomorphic to the poset 
${\text{\bf Idl}}(L)$  of lattice ideals of $L$.
\end{proposition}

\proof 
We show the claim for the filters. The isomorphism is given by 
$F(L^\delta )\to {\text{\bf Filt}}(L),\  x\mapsto {\uparrow} x \cap L$ and 
${\text{\bf Filt}}(L)\to F(L^\delta ),\ F\mapsto \bigwedge F$. 
It is clear that each $x\in F(L^\delta )$ satisfies $x=\bigwedge ({\uparrow} x\cap L)$. 
To show that $F={\uparrow}(\bigwedge F) \cap L$ compactness is used.
\qed\smallskip

Note that it is a consequence of compactness that the elements of a canonical extension that 
are both filter and ideal elements are exactly the elements of the original lattice. We call these
elements lattice elements.

\begin{proposition}
Let $L$ be a lattice,  and $L^\delta$ a canonical extension of $L$. Then the order on 
the subposet $F(L^\delta )\cup I(L^\delta )$ of $L^\delta$ is uniquely determined by $L$.
\end{proposition}
This follows as we can prove, using density and compactness, that the order is given by
\begin{enumerate}[(i)]
\item $x\leq x'$ if and only if $F_{x'}\subseteq F_x$;
\item $x\leq y$ if and only if $F_x\cap I_y\ne\emptyset$;
\item $y\leq x$ if and only if $a\in I_y$, $b\in F_x$ implies $a\leq b$;
\item $y\leq y'$ if and only if $I_y\subseteq I_{y'}$.
\end{enumerate}
\noindent Here $x,x'$ stand for elements in $F(L^\delta )$; $F_x,F_{x'}$ for the 
corresponding filters and $y,y'$ stand for elements in $I(L^\delta )$; $I_y,I_{y'}$ 
for the corresponding ideals.

Now the uniqueness of the canonical extension follows modulo the well-known abstract 
characterisation of MacNeille completion.

\begin{theorem} {\rm(uniqueness of canonical extensions)}
Let $L$ be a lattice. Then the canonical extension of $L$ is unique up to 
an isomorphism fixing $L$.
\end{theorem}

\proof
It is clear from the above proposition that for any canonical extension $L\hookrightarrow L'$
of $L$, the poset $Int(L')=F(L')\cup I(L')$ is uniquely determined. The MacNeille completion 
of a poset is the unique completion in which the original poset is both join-dense and 
meet-dense. The density condition for canonical extensions tells us that $Int(L')$ is both 
join-dense (because of the filter elements) and meet-dense (because of the ideal elements) 
in $L'$ and thus $L'$ is uniquely determined as the MacNeille completion of $Int(L')$. 
\qed\smallskip

Note that this uniqueness proof also provides a key to existence: one can build the canonical 
extension of any lattice by taking the MacNeille completion of the amalgam of the ideal lattice 
and the order dual of the filter lattice of $L$ according to the four conditions given above. This 
construction has the virtue of not using the axiom of choice. However, by uniqueness, it will 
produce the embedding of $L$ into the dense completion defined by its dual space whenever
the latter exists. 

Remarkably, even in the non-distributive case, the canonical extension of a lattice satisfies a 
restricted complete distributivity condition. We do not give the straight forward proof which may
be found in \cite{GeHa01}.

\begin{proposition}  {\rm(restricted distributivity for canonical extensions)}
Let $L$ be a bounded lattice and $\mathcal Y$ a family of down-directed subsets of $L$, 
viewed as a family of subsets of the canonical extension $L^\delta$ of $L$. Then $\mathcal Y$ 
satisfies the complete $\bigvee\!\!\bigwedge$-distributive law. Dually,  if $\mathcal Y$ is a family 
of up-directed subsets of $L$  then $\mathcal Y$ satisfies the $\bigwedge\!\!\bigvee$-distributive 
law relative to $L^\delta$.
Here $\mathcal Y$ is said to satisfy the complete $\bigvee\!\!\bigwedge$-distributive law provided
$$\bigvee\{\,\bigwedge Y\mid Y\in \mathcal Y\,\}=\bigwedge\{\,\bigvee Z\mid Z\in \mathcal Y^\sharp\,\}$$
where $\mathcal Y^\sharp=\{\,Z\subseteq L\mid Y\cap Z\neq\emptyset\text{ for all }Y\in\mathcal Y\,\}$
and the $\bigwedge\!\!\bigvee$-distributive 
law is defined order dually.
\end{proposition}

From this one can show that the canonical extension of a distributive lattice is distributive and,
using the axiom of choice, that it is completely distributive. Using the axiom of choice one can 
also show that the canonical extension of any lattice is join generated by the set $J^\infty(L^\delta)$
of completely join irreducible elements of the canonical extension. In the distributive setting, these 
of course correspond to the prime filters of the original lattice and we get that $L^\delta$ is 
isomorphic to the upset lattice of  $J^\infty(L^\delta)$. By symmetry, the order dual statements hold
about the collection of completely meet irreducible elements of $L^\delta$, $M^\infty(L^\delta)$.

Given that canonical extensions satisfy some directed infinite distributivity conditions, it is 
natural to wonder whether they must always be continuous lattices. For distributive lattices
this is true but it is not the case in general. We give an example here of a canonical extension 
that is not meet-continuous and thus, as it is a complete lattice, not a continuous lattice, see 
\cite[Proposition I-1.8, p.56]{comp2}.

\begin{example}  
 (A canonical extension that is not continuous)  Let 
$$
L=\{0,1\}\cup\{a_{ij}\mid i,j\in \Bbb N\}
$$
where $0$ is the bottom, $1$ is the top, and 
$$
a_{ij}\geq a_{kl}\quad\iff\quad( i+j\leq k+l \text{ and } i\geq k).
$$
This lattice, see Figure~\ref{fig:non-cotacnext}, is non-distributive as, e.g., 
$1,a_{20},a_{11},a_{00},a_{02}$ form a copy 
of the lattice $N_5$. $L$ satisfies ACC and thus the intermediate structure is isomorphic 
to the filter completion of $L$ which is obtained by adding filter elements $x_i, i\in \Bbb N$ 
with $x_i\leq a_{ij}$ for all $i$ and $j$ (and then $x_i\leq x_k$ whenever $i\leq k$). The 
resulting structure is complete and is thus the canonical extension of $L$. 
\begin{figure}
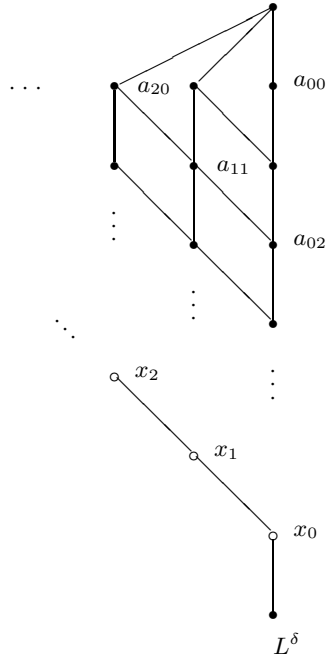

$$
{\hbox{
				\putt 4 0 {\hbox{
  			\spoint 0 23
 			 \spoint 0 20
 			 \labelr 0  20 {a_{00}}
  			\spoint 0 17
 			 \spoint 0 14
   			\labelr 0 14 {a_{02}}
  			\spoint 0 11
  			\vline 0 11 12
   						\putt -.15 8 {\hbox{
    							\putt 0 0 {$\cdot$}
    							\putt 0 .5 {$\cdot$}
   							 \putt 0 1 {$\cdot$}
    					 	}}
 			\point 0 3
  			\labelr 0 3 {x_0}
				\spoint 0 0
				\vline 0 0 3
 \RSER .1 3 3
 \RSER .1 17  3
  \RSER .1 14  3
   \RSER .1 11  3
    \RSWR -3 20  3
  			\spoint -3 20
 			 \spoint -3 17
  			 \labelr -2.9 16.8 {a_{11}}
				  \spoint -3 14
				  \vline -3 14 6
  						 \putt -3.15 11 {\hbox{
   						 \putt 0 0 {$\cdot$}
    						\putt 0 .5 {$\cdot$}
   						 \putt 0 1 {$\cdot$}
     						}}
 			\point -3 6
  			\labelr -3 6 {x_1}
\RRSWR -6 20 6
\RSER -2.9 17 3
\RSER -2.9 14  3
\RSER -2.9 6 3
 			 \spoint -6 20
   			\labelr -5.9 19.8 {a_{20}}
  			\spoint -6 17
  			\vline -6 17 3
  						 \putt -6.15 14 {\hbox{
    						\putt 0 0 {$\cdot$}
    						\putt 0 .5 {$\cdot$}
   						 \putt 0 1 {$\cdot$}
  						  }}
					 \point -6 9
 			 \labelr -6 9 {x_2}
 						\putt -7.5 10.5 {\hbox{
 											 \putt -.2  -.2 {\hbox{       
   											 \putt  0 0 {$\cdot$}  
    											\putt -.3 .3  {$\cdot$}
    											 \putt -.6 .6 {$\cdot$}
												                           }}
												\putt -1.5 9.2 {\hbox{
    								\putt 0 0 {$\cdot$}
    								\putt -.5 0 {$\cdot$}
   								 \putt -1 0 {$\cdot$}
   								 }}
									}}
				\putt 0 -1.5  {$L^\delta $}
				}}   
}}                   
$$ 
\caption{Non-continuous canonical extension}\label{fig:non-cotacnext}
\end{figure}
To see that $L^\delta$ is not meet-continuous note that 
$a_{00}\wedge(\bigvee_{i=0}^\infty x_i)=a_{00}\wedge 1=a_{00}$ while 
$\bigvee_{i=0}^\infty (a_{00}\wedge x_i)=\bigvee_{i=0}^\infty x_0=x_0$.
\end{example}

\section{Morphisms, maps, and additional operations}\label{sec:addop}

In domain theory maps are extended using directed join density. In canonical extensions the
original lattice may be neither meet nor join dense but two layers of joins and meets are needed.
However, by introducing a topology we can translate this to a topological setting in which the 
original lattice is topologically dense in the canonical extension.

\begin{definition}
Let $L$ be a lattice. We denote by $\delta $, $\delta ^{\uparrow }$ and $\delta ^{\downarrow }$ 
the topologies on $L^{\delta}$ having as bases, respectively, the sets of the forms 
${\uparrow}x\cap{\downarrow}y=[x,y]$, ${\uparrow}x=[x,1]$ and ${\downarrow}y=[0,y]$, with 
$x\in F(L^{\delta })$ and $y\in I(L^{\delta})$.
\end{definition}

We will denote the interval topology on any poset by $\iota$ and its one-sided parts, the 
upper topology and the lower topology, by $\iota ^{\uparrow }$ and $\iota ^{\downarrow }$, 
respectively.  Further, we denote the Scott topology by $\sigma^{\uparrow }$, the dual Scott
topology by $\sigma^{\downarrow}$, and the double Scott topology by $\sigma$. We have 
the following basic facts about the topology $\delta$.

\begin{theorem}
\label{deltatop}
Let $L$ be lattice. The $\delta$-topologies are refinements of the $\sigma$-topologies and thus
also of the $\iota$-topologies and the space $(L^{\delta },\delta)$ is Hausdorff. The set $L$ is 
dense in $(L^{\delta },\delta)$ and the elements of $L$ are exactly the isolated points of the space. 
\end{theorem}

\proof
Since the filter elements of a canonical extension join-generate it, by directed joins, it is clear that 
$\sigma^{\uparrow }\subseteq\delta^{\uparrow }$ and by order duality 
$\sigma^{\downarrow }\subseteq\delta^{\downarrow}$  and thus also $\sigma\subseteq\delta$. To 
see that $\delta$ is Hausdorff, let $u,v\in L$ with $u\nleq v$, then there is $x\in F(L^{\delta })$ with 
$x\leq u$ but $x\nleq v$. Now since $x\nleq v$ there is $y\in I(L^{\delta})$ with $v\leq y$ but 
$x\nleq y$. That is, ${\uparrow}x$ and ${\downarrow}y$ are disjoint $\delta$ open sets separating 
$u$ and $v$. The set $A$ is dense in $L^{\delta }$ since each non-empty basic intervals $[x,y]$ 
contains a lattice element by compactness. Finally, for $a\in L$, the interval $[a,a]=\{a\}$ is open, 
and $a$ is therefore isolated. On the other hand, since $L$ is dense in $(L^{\delta },\delta )$, it 
follows that if $\{u\}$ is open then $u\in L$.
\qed\smallskip

Further basic facts about this topology are that it is stable under order duality and that it commutes with
Cartesian product (i.e. is productive). We note also that if $L$ is distributive, then $L^\delta$ is a Priestley 
space in its interval 
topology which is also equal to the double Scott topology and is generated by the complementary pairs
${\uparrow}p,{\downarrow}\kappa(p)$, with $p\in J^\infty(L^\delta)$ and $\kappa(p)\in M^\infty(L^\delta)$ 
given by $p\leq u$ iff $u\nleq\kappa(p)$ for $u\in L^\delta$. In fact, the topology generated
by upsets of elements of $J^\infty(L^\delta)$ and downsets of elements of $M^\infty(L^\delta)$ also plays 
an important role in the theory of canonical extensions in general \cite{GeHaVe}. 

In defining and investigating extensions of maps $f:K\rightarrow L$ between lattices to maps between their 
canonical extensions, we make use of the various topologies on $K^{\delta }$ and $L^{\delta }$. Since
several topologies have been defined on each set, it is often necessary to specify which ones are under 
consideration. In general, if $\tau $ and $\mu $ are topologies on the sets $X$ and $Y$, and if the map 
$f:X\rightarrow Y\,$ is continuous relative to $\tau $ on $X$ and $\mu $ on $Y$, then we write that $f$ is 
$(\tau ,\mu)$-continuous.

\begin{definition}
Let $L$ be a lattice and $C$ a complete lattice. For any map $f:L\rightarrow C$, 
and for all $u\in L^\delta$, we define 
\begin{eqnarray*}
f^\sigma(u)&=&\underline{\lim}_{\delta}f(u)=\bigvee\{\bigwedge f(U\cap L)\mid u\in U\in \delta\}\\
&=&\bigvee\{\bigwedge f([x,y]\cap L)\mid F(L^\delta)\ni x\leq u\leq y\in I(L^\delta)\}, \\[1ex]
f^\pi(u)&=&\overline{\lim }_{\delta}f(u) =\bigwedge \{\bigvee f(U\cap L)\mid u\in U\in \delta \}\\
&=&\bigwedge\{\bigvee f([x,y]\cap L)\mid F(L^\delta)\ni x\leq u\leq y\in I(L^\delta)\}.
\end{eqnarray*}
In particular, for maps $f:L\to M$ between lattices, we define $f^\sigma$ and $f^\pi$ by considering 
the composition of $f$ with the embedding of $M$ in $M^\delta$.
\end{definition}

Note that, as each point of $L$ is isolated in the $\delta$-topology
it follows that both 
of the functions defined above
are extensions of $f$, that is, agree with $f$ on $L$. Also, as the $\delta$ topology commutes with 
products, the lifting of operations is just a special case of lifting of maps.
This is of course the well-known upper and lower envelope constructions from topology and, under 
some restrictions, they are, respectively, the largest $(\delta,\iota ^{\uparrow })$-continuous function 
that is below $f$ on $L$ and  the least $(\delta,\iota ^{\downarrow })$-continuous function that is above 
$f$ on $L$. A careful analysis of when this works is in the second
author's Ph.D. thesis \cite{VoPhD}. Here we record the following facts.

\begin{proposition}\label{prop:univenv}
Let $f:L\to M$ be a map between lattices. Then $f^\sigma:L^\delta\to M^\delta$ is 
$(\delta,\sigma^{\uparrow})$-continuous and thus also $(\delta,\iota ^{\uparrow })$-continuous.
Furthermore,
\begin{enumerate}
\item If $f$ is order preserving or reversing, then $f^\sigma$ is the largest 
$(\delta,\iota^{\uparrow })$-continuous function that is below $f$ on $L$;
\item If $\sigma^{\uparrow}$ has a basis of principal up-sets, i.e.~if
 $M^\delta$ is algebraic, then $f^\sigma$ is the largest 
$(\delta,\sigma^{\uparrow })$-continuous function that is below $f$ on $L$.
\end{enumerate} 
Dual statements hold about $f^\pi$.
\end{proposition}

When the envelopes are the largest $(\delta,\iota ^{\uparrow })$-continuous functions above, 
respectively smallest $(\delta,\iota ^{\downarrow })$-continuous functions below, the original
function we will say that the envelopes of $f$ are universal.
This is the case, by (1), for operations that are monotone (that is, order preserving or
reversing in each coordinate). We shall see, in the next section, that on canonical extensions 
of lattices lying in finitely generated varieties, the Scott topology is equal to the upper topology 
and has a basis of principal up-sets so that the envelopes are universal for any mapping between 
lattices lying in finitely generated lattice varieties.

We give a few examples of extensions of maps.

\begin{example}\label{modalex}  {\rm (of the $\sigma$- and $\pi$-extensions of a
modal operator)}
The following is a notorious example from modal logic. It illustrates that modal axioms may 
fail to be preserved by canonical extension. Let $B$ be the Boolean algebra of finite and co-finite 
subsets of $\Bbb N$ and consider the relation $>$ on $\Bbb N$. The Boolean algebra $B$
is closed under the operation 
$\Diamond(S)=\{n\mid\exists m\ (n>m\text{ and }m\in S\}$ 
since it gives ${\uparrow}(\min(S)+1)$ for any non-empty set $S$. It is straight forward to check
that the modal algebra $(B,\Diamond)$ satisfies the G\"odel-L\"ob axiom:
$$
\Diamond(\neg\Diamond a\wedge a)\geq \Diamond a.
$$
It is clearly true for $\emptyset$ since $\Diamond\emptyset=\emptyset$. For any non-empty set
$S$, we have $\Diamond(S)={\uparrow}(\min(S)+1)$ and thus the complement contains $\min(S)$
and we get $\Diamond(\neg\Diamond(S)\wedge S)=\Diamond (S)$. The canonical extension of $B$ is 
easily seen to be the powerset of 
$\Bbb N_\infty~=~\Bbb N\cup~\{\infty\}$ with the embedding of $B$ into $\mathcal P(\Bbb N_\infty)$ 
which sends each finite subset of $\Bbb N$ to itself and each co-finite subset of $\Bbb N$ to its 
union with $\{\infty\}$. Thus the singleton $\{\infty\}$ is the filter element which is the meet of all 
the co-finite elements of $B$. We have
$$
\Diamond^\sigma(\{\infty\})=\bigcap\{\Diamond(S)\mid S\text{ is co-finite}\}.
$$
Since $\Diamond(S)$ for a co-finite set can be ${\uparrow} n$ for any $n\in\Bbb N$ it follows that
$\Diamond^\sigma(\{\infty\})=\{\infty\}$ and thus 
$\Diamond(\neg\Diamond\{\infty\}\wedge\{\infty\})=\emptyset\ngeq \Diamond\{\infty\}$.
\end{example}         

A map $f$ between lattices is called {\it smooth} provided its $\sigma$- and $\pi$-extensions are equal. 
In this case we denote the extension by $f^\delta$ to stress its order-symmetry. 

\begin{example}\label{non-smoothex}   {\rm (of a non-smooth operation)}         
Let $X$ be an infinite set and let $B$ be the Boolean algebra of all subsets of $X$ which are either
finite or co-finite. Consider the map $f:B^2\to B$ defined by $f(b_1,b_2)=0_B=\emptyset$ if $b_1$
and $b_2$ are disjoint and $f(b_1,b_2)=1_B=X$ otherwise. As in the above example, the canonical 
extension of $B$ is the powerset of $X_\infty~=~X\cup~\{\infty\}$ where $\infty\not\in X$ with the 
embedding of $B$ into $\mathcal P(X_\infty)$ which sends each finite subset of $X$ to itself and each 
co-finite subset of $X$ to its union with $\{\infty\}$.

Let $u\in B^\delta=\mathcal P(X_\infty)$ be a subset of $X$ that is neither finite nor co-finite. We claim that 
$f^\sigma(u,\neg u)=0$ whereas $f^\pi(u,\neg u)=1$. 
$$
f^\sigma(u,\neg u)=\bigvee\{\bigwedge f([\,\overline{s},\overline{t}\ ]\cap B^2)\mid 
F((B^2)^\delta)\ni\overline{s}\leq(u,\neg u)\leq\overline{t}\in I((B^2)^\delta)\}
$$
Note that canonical extension commutes with product so that $(B^2)^\delta=(B^\delta)^2$,
$F((B^2)^\delta)=(F(B^\delta))^2$, and $I((B^2)^\delta)=(I(B^\delta))^2$. Now pick 
$\overline{s}=(s_1,s_2)\in (F(B^\delta))^2$ and $\overline{t}=(t_1,t_2)\in (I(B^\delta))^2$ with 
$\overline{s}\leq(u,\neg u)\leq\overline{t}$. It is not hard to verify that 
$s\in B^\delta=\mathcal P(X_\infty)$ is a filter element if and only if it is finite or contains $\infty$. By choice of 
$u$ we have $\infty\notin u$ and thus $\infty\notin s_1\leq u$ and $s_1$ must be a finite subset of $X$.
That is, $s_1\in [s_1,t_1]\cap B$ is a finite subset of $u$. Now $s_2\leq\neg u\leq\neg s_1\in B$ and it 
follows by compactness that there is $b_2\in B$ with $s_2\leq b_2\leq\neg s_1\wedge t_2\leq t_2$. Since 
$s_1$ and $\neg s_1$ are disjoint, so are $s_1$ and $b_2$ and we have $f(s_1,b_2)=0$. Also 
$(s_1,b_2)\in [\,\overline{s},\overline{t}\,]\cap B^2$ so $\bigwedge f(\,[\overline{s},\overline{t}\,]\cap B^2)=0$ and 
$f^\sigma(u,\neg u)=0$  as claimed.

Now consider
$$
f^\pi(u,\neg u)=\bigwedge\{\bigvee f([\,\overline{s},\overline{t}\,]\cap B^2)\mid 
F((B^2)^\delta)\ni\overline{s}\leq(u,\neg u)\leq\overline{t}\in I((B^2)^\delta)\}
$$
and pick $\overline{s}=(s_1,s_2)\in (F(B^\delta))^2$ and $\overline{t}=(t_1,t_2)\in (I(B^\delta))^2$ 
with $\overline{s}\leq(u,\neg u)\leq\overline{t}$. We have $\neg u\leq t_2\in I(B^\delta)$. Now, an
element  $t\in B^\delta=\mathcal P(X_\infty)$ is an ideal element if and 
only if $t$ is co-finite or doesn't contain $\infty$. By choice of $u$ we have $\infty\in \neg u$ so that 
$\infty\in t_2$ and thus $t_2$ must be co-finite. It follows that $t_2\in[s_2,t_2]\cap B$. Since $u$ is 
not finite, $u\wedge t_2\neq 0$. Let $b\in B$ be any finite non-empty subset of $u\wedge t_2$. 
Then $b\leq u\leq t_1$ and by an argument similar to the one above, we obtain a $b_1\in B$ with 
$s_1\leq s_1\vee b\leq b_1\leq t_1$. Now $(b_1,t_2)\in[\,\overline{s},\overline{t}\,]\cap B^2$ and  
$0\neq b\leq b_1\wedge t_2$ so that $f(b_1,t_2)=1$. It follows that $f^\pi(u,\neg u)=1$.
\end{example}

The fact 
that 
 the universal properties of the upper and lower extensions of a map are asymmetric
with respect to the topology used on the domain and codomain has as consequence that, in 
total generality, extensions do not commute with composition \cite[Ex.2.34]{GeJo04} so that 
canonical extension isn't functorial when considering arbitrary set maps between lattices. The 
paper \cite{GeJo04} analysed the situation in detail and in \cite{GeHa01} some of the results 
were generalised to the lattice setting. A simple general fact encompassing most applications 
in logic is: canonical extension is functorial for homomophisms of algebras that are lattices 
with additional basic operations each of which is order-preserving or -reversing in each of its 
coordinates (such algebras are called monotone lattice expansions). 

Preservation of identities when moving to the canonical extension is also closely tied to 
compositionality of the extension of maps and, as explained in detail in \cite{GeJo04}, 
compositionality results can in many cases be inferred by an analysis of the topological properties 
of the extensions of maps with particular properties. Examples are given in the following theorem. 

\begin{theorem}\label{thm:addop}
Let $K,L,M,N$ be lattices, $h:K\to L$ a lattice homomorphism, and $f:M\to N$ a map with 
universal envelopes.
Then the following hold:
\begin{enumerate}
\item If $f$ has a $(\delta,\iota)$-continuous extension, $\tilde{f}:M^\delta\to N^\delta$, then $f$ is 
smooth and $f^\delta=\tilde{f}$. 
\item $h$ is smooth and $h^\delta:K^\delta\to L^\delta$ is a complete homomorphism and is both
$(\delta,\delta)$- and $(\iota,\iota)$-continuous;
\item If $N=K$ then $(hf)^\sigma=h^\sigma f^\sigma$;
\item If $L=M$ and $h$ is surjective then $(fh)^\sigma=f^\sigma h^\sigma$
\item If $M=M_1\times\ldots\times M_n$ and $f$ preserves joins in each coordinate
(i.e., $f$ is an operator) and $M$ is distributive, then $f^\sigma$ is 
$(\iota^{\uparrow},\iota^{\uparrow})$-continuous.
\end{enumerate}
\end{theorem}

\proof
The facts (1),(3) and (4) are proved for distributive lattices in \cite[Cor.2.17]{GeJo04}, \cite[Lem.3.3]{GeJo04}, 
and \cite[Lem.3.6]{GeJo04}, respectively, and an inspection of the proofs readily shows that they are still valid 
in the lattice setting.

The fact that lattice homomorphisms are smooth and lift to complete lattice homomorphisms is 
proved in \cite{GeHa01}. The fact that $h^\delta$ is $(\delta,\delta)$-continuous is proved for 
distributive lattices in \cite[Thm.2.24(iii),(iv)]{GeJo04} and an inspection of the proof readily shows 
that it is true in the lattice setting as well. The $(\iota,\iota)$-continuity is another matter (see (4) below). 
Let $v\in L^\delta$. For each $u\in K^\delta$, we have
\[
h^\delta(u)\leq v\ \iff\ u\leq (h^\delta)^\sharp(v)
\]
where $(h^\delta)^\sharp$ is the upper adjoint of $h^\delta$. Thus the same holds for the negation of 
these inequalities,
i.e. $(h^\delta)^{-1}(({\downarrow}v)^c)=({\downarrow}(h^\delta)^\sharp(v))^c$, 
where $(\ )^c$ denotes complement, 
and thus 
$h^\delta$ is $(\iota^{\uparrow},\iota^{\uparrow})$-continuous. By symmetry $h^\delta$ is $(\iota^{\downarrow},\iota^{\downarrow})$-continuous.

The proof of (4), which is the cornerstone of the paper \cite{GeJo94}, relies on the fact that 
$\iota^{\uparrow}$ is generated by $\{{\uparrow}p\mid p\in J^\infty(L^\delta)\}$ in the distributive
setting, see e.g. \cite[Lem.4.2]{GeJo94}. 
\qed\smallskip

We now illustrate the use of these tools by proving the following propositions. Note that it is not 
specified in the following propositions whether we are using the $\sigma$- or the $\pi$-extension 
in taking the canonical extensions of the additional operations. The point is that the results hold 
in either case. 

\begin{proposition}\label{prop1}
Let $(A,f)$ and $(B,g)$ be lattices with additional $n$-ary operation with universal envelopes, 
and let $h:(A,f)\rightarrow (B,g)$ be a homomorphism. If $g$ is smooth then $h$ 
lifts to a homomorphism between the canonical extensions. 
\end{proposition}

\proof
Since $h:(A,f)\to(B,g)$ is a homomorphism, we have $hf=gh^{[n]}$ and thus 
$(hf)^\sigma=(gh^{[n]})^\sigma$. Now $(hf)^\sigma=h^\sigma (f)^\sigma$ by 
Theorem~\ref{thm:addop}(3). Note that $g^\delta(h^{[n]})^\delta$ is  
$(\delta,\iota)$-continuous since $(h^{[n]})^\delta$ is $(\delta,\delta)$-continuous
by Theorem~\ref{thm:addop}(2) and $g^\delta$ is $(\delta,\iota)$-continuous by 
hypothesis. Also, $g^\delta(h^{[n]})^\delta$ is an extension of $gh^{[n]}$ so by 
Theorem~\ref{thm:addop}(1), we have $(gh^{[n]})^\sigma=(gh^{[n]})^\pi=g^\delta(h^{[n]})^\delta$. 
That is, $h^\delta (f)^\sigma=g^\delta(h^{[n]})^\delta$ and the homomorphism lifts.
\qed\smallskip

\begin{lemma}\label{lem1}
Let $A$ and $B$ be lattices and  $h:A\twoheadrightarrow B$ a surjective homomorphism. 
Then $h^\delta:A^\delta\twoheadrightarrow B^\delta$ is a $(\delta,\delta)$-open mapping.
\end{lemma}

\proof
Note first that surjective morphisms lift to surjective morphisms \cite{GeHa01}.
If $x$ and $y$ are filter and ideal elements in $A^\delta$, respectively, then clearly $h^\delta(x)$
and $h^\delta(y)$ are filter and ideal elements in $B^\delta$ since $h^\delta$ preserves arbitrary 
meets and joins. Also, using the fact that $h^\delta$ is surjective, it is straight forward to check 
that $h^\delta([x,y])=[h^\delta(x),h^\delta(y)]$ (for this note that if 
$h^\delta(x)\leq h^\delta(u)\leq h^\delta(y)$ then 
$h^\delta(x)\leq h^\delta((u\vee x)\wedge y)\leq h^\delta(y)$ and $x\leq(u\vee x)\wedge y\leq y$). 
Now the result follows as forward image always preserves union.
\qed\smallskip

\begin{proposition}\label{prop2}
Let $(A,f)$ and $(B,g)$ be lattices with additional $n$-ary operation with universal 
envelopes, and let 
$h:(A,f)\twoheadrightarrow (B,g)$ a surjective homomorphism. If $f$ is smooth then 
so is $g$. If the extension of $f$ is $(\iota,\iota)$-continuous and $h^\delta$ sends
$\iota$-open $h^\delta$-preimages to $\iota$-opens, then the extension of $g$ is also 
$(\iota,\iota)$-continuous.
\end{proposition}

\proof
Note that $h$ lifts to a homomorphism of the canonical extensions by 
Theorem~\ref{thm:addop} parts (3) and (4).
Let $U$ be $\iota$-open in $B^\delta$. Then $(h^\delta\circ f^\sigma)^{-1}(U)$ is 
$\delta$-open in $(A^\delta)^n$ since $f^\sigma$ is $(\delta,\iota)$-continuous by 
assumption and $h^\delta$ is $(\iota,\iota)$-continuous by 
Theorem~\ref{thm:addop}(2). Now $h^\delta\circ f^\sigma=g^\sigma\circ (h^\delta)^{[n]}$ 
since $h$ lifts to a homomorphism of the canonical extensions. It follows that 
$(g^\sigma\circ (h^\delta)^{[n]})^{-1}(U)=((h^\delta)^{[n]})^{-1}\circ (g^\sigma)^{-1}(U)$ 
is $\delta$-open in $(A^\delta)^n$. We now use the lemma to conclude that the 
lifting $(h^\delta)^{[n]}=(h^{[n]})^\delta$ of the surjective homomorphism 
$h^{[n]}:A^n\twoheadrightarrow B^n$ which is obtained by doing $h$ in each 
coordinate, is a $(\delta,\delta)$-open mapping. We thus conclude that 
$(h^\delta)^{[n]}(((h^\delta)^{[n]})^{-1}\circ (g^\sigma)^{-1}(U))$ is $\delta$-open in 
$(B^\delta)^n$. Finally note that, as $(h^\delta)^{[n]}$ is surjective, 
$(h^\delta)^{[n]}(((h^\delta)^{[n]})^{-1}(S))=S$ for any subset of $(B^\delta)^n$.
We conclude that $(g^\sigma)^{-1}(U)$ is $\delta$-open in $(B^\delta)^n$ and thus
$g$ is smooth. 

For the statement on $(\iota,\iota)$-continuity, note that the openness of the map 
$(h^\delta)^{[n]}$ in the proof above is only needed on $(h^\delta)^{[n]}$-saturated 
opens and this is a consequence of the corresponding statement for $h^\delta$.
Thus, with the given assumptions, the same proof goes through for the 
$(\iota,\iota)$-continuity.
\qed\smallskip

A class of similar lattices with additional operations is called a class of lattice expansions.

\begin{corollary}
Let $\mathcal K$ be a class of lattice expansions for which the envelopes of the basic 
operations are universal. The operator $H$, taking homomorphic images of algebras, 
preserves smoothness.
\end{corollary}

\section{Canonical extensions in finitely generated varieties}\label{sec:fgv}

In this final section we illustrate the theory by giving a few consequences
for lattice expansions that lie within finitely generated varieties, 
varieties generated by a finite collection of finite algebras. 
These are simple 
consequences, mainly of the results in \cite{GeJo04} and \cite{GeHa01} but have 
not been published yet. The main result of \cite{GeJo04} (first published in \cite{GeJo99}) 
has as consequence that all finitely generated varieties of bounded distributive lattice 
expansions are canonical and in \cite{GeHa01} it was shown that this result goes through 
to finitely generated monotone bounded lattice varieties.

These results are based on two facts. First, the observation (also behind the 
famous J\'onsson Lemma of universal algebra) that any product of lattice expansions 
is isomorphic to a Boolean product of all the ultraproducts formed from the given 
product. And secondly, the following result which is central in \cite{GeJo99} and 
\cite{GeJo04} in its distributive lattice incarnation and is central in \cite{GeHa01} in 
its general form for arbitrary bounded lattices. We give the simple proof for arbitrary 
bounded lattices here for completeness.

\begin{theorem}\label{BP-lat}  
{\rm(Canonical extensions of Boolean products)}\ 
Let $(L_x)_{x\in X}$ be a family of bounded lattices. 
If $L\leq\prod_{X}L_x$ is a Boolean product, then 
$L^\delta = \prod_{X}L_x^\delta$.
\end{theorem}

\proof 
 We first show that the inclusion of $L$ into $\prod_XL_x^\delta $ given by
 the composition of the inclusion of $L$ into $\prod_XL_x$ followed by the 
 coordinate-wise embedding of $\prod_XL_x$ into $\prod_XL_x^\delta$ 
 yields  a canonical extension. 
As each $L_x^\delta $ is complete, the product $\prod_XL_x^\delta $ is a 
complete lattice. Suppose $x\in X$ and $p\in L_x^\delta $ is a filter element.
Define $u_{x,p}\in\prod_XL_x^\delta $ by setting $u_{x,p}(x)=p$ and 
$u_{x,p}(y)=0$ for $y\neq x$. We first show that $u_{x,p}$ is a meet in 
$\prod_XL_x^\delta $ of elements from $L$. It then follows that every 
element of $\prod_XL_x^\delta $ is a join of meets of elements of $L$, 
and by a dual argument, a meet of joins of elements of $L$. 

To show that $u_{x,p}$ is a meet of elements of $L$, note first that 
$p$ is a meet in $L_x^\delta $ of a family $S$ of elements of $L_x$. 
As $L\leq\prod_XL_x$ is subdirect, for each $s\in S$ there is some 
$u_s\in L$ with $u_s(x)=s$. Using the Patching Property, 
for each clopen neighbourhood $N$ of $x$, and each $s\in S$, 
we have $u_s|N\cup 0|N^c$ is an element of $L$. Then, the 
meet of 
$$
\{\, (u_s|N\cup 0|N^c) \mid s\in S, \ x\in N\text{ clopen}\,\}
$$ 
is equal to $u_{x,p}$. This shows that the inclusion of $L$ into 
$\prod_XL_x^\delta $ is dense.

Finally we show that the inclusion of $L$ into $\prod_XL_x^\delta $ is 
compact. Suppose that $S$ is a filter of $L$, $T$ is an ideal of $L$, and 
$\bigwedge S\leq\bigvee T$. For each $x\in X$ let 
$S_x =\{\, u(x) \mid u\in S\,\} $ and let $T_x=\{\, v(x)\mid v\in T\,\} $. Then 
$\bigwedge S_x\leq\bigvee T_x$ for each $x\in X$. As $L_x^\delta $ is a 
canonical extension of $L_x$, $S_x\cap T_x\neq\emptyset$, hence there 
are $u_x\in S$ and $v_x\in T$ with $u_x(x)=v_x(x)$. As equalisers in a 
Boolean product are clopen, $u_x$ and $v_x$ agree on some clopen 
neighbourhood $N_x$ of $x$. Then, as $X$ is compact, and 
$\{\, N_x\mid x\in X\,\}$ is an open 
cover of $X$, there is a finite family $x_1,\ldots,x_n$ with $N_{x_1},
\ldots,N_{x_n}$ a cover of $X$. We assume, without loss of generality,
 that $N_{x_1},\ldots,N_{x_n}$ are pairwise disjoint. Let $w$ be the 
function which agrees with $u_{x_i}$, hence also with $v_{x_i}$, 
on $N_i$ for $i=1,\ldots, n$. By the Patching Property, 
$w$ is an element of $L$. Also, $w$ is the join of the $n$ functions 
agreeing with $u_{x_i}$ on $N_{x_i}$ and defined to be $0$ elsewhere, 
hence $w$ is in the ideal $S$. Similarly $w$ is the meet of the $n$ 
functions agreeing with $v_{x_i}$ on $N_{x_i}$ and $1$ elsewhere, 
hence $w$ is in the filter $T$. Thus, $S\cap T\neq\emptyset$. 
This shows that the inclusion of $L$ into $\prod_XL_x^\delta $ 
is compact.
\qed\smallskip 

It is a fundamental universal algebraic fact that if a class $\mathcal K$ 
generates the variety $\mathcal V$, then ${\mathcal V}=HSP({\mathcal K})$
where $H,S,P$ are the operators closing a class under homomorphic images, 
subalgebras, and products, respectively. By the above mentioned observation, 
$P({\mathcal K})=P_BP_\mu({\mathcal K})$ where $P_B$ and $P_\mu$ are the 
operators closing a class under Boolean products and ultraproducts, respectively.
Since an ultraproduct of a single finite structure is always isomorphic to the structure
itself, it follows that for a finite lattice expansion $A$, ${\mathcal V}(A)=HSP_B(A)$. 
Many theorems, including the main canonicity theorems of \cite{GeJo04,GeHa01}
are proved by showing that $H$, $S$, and $P_B$ all three preserve canonicity. These 
three operators preserve many other nice properties and that is what we want to illustrate 
here. 

We start with a somewhat technical proposition drawing on work in domain theory.
The conclusion of the proposition identifies what is at stake here. An upper, respectively 
lower, tooth in a poset is the upset, respectively downset, of a finite subset. A perfect 
lattice is a complete lattice in which the completely join irreducibles are join-dense and
the completely meet irreducibles are meet-dense.

\begin{proposition}\label{prop:star}
Let $C$  be a perfect lattice with the following properties:
\begin{align*}
(\bigstar)\hskip.2cm\qquad &\forall p\in J^\infty(C)\quad ({\uparrow}p)^c={\downarrow} M_p\quad
\mbox{ where }M_p\subseteq M^\infty(C)\mbox{ is finite};
\\
(\bigstar)^\partial\qquad &\forall m\in M^\infty(C)\quad ({\downarrow}m)^c={\uparrow} J_m\quad
\mbox{ where }J_m\subseteq J^\infty(C)\mbox{ is finite}.
\end{align*}
Then $C$ is doubly algebraic and the Scott and the upper topologies on $C$ are equal and 
this topology is spectral. Dually, the dual Scott and the lower topologies on $C$ are equal and 
this topology is spectral as well. The bases of compact-opens of these two topologies come in 
complementary pairs of upper and lower teeth and the join of the two topologies makes $C$ 
into a Priestley space.
\end{proposition}

\proof
We first show that $C$ is algebraic. Denote the finite join closure of $J^\infty(C)$ by 
$J^\infty_\omega(C)$ and the finite meet closure of $M^\infty(C)$ by $M^\infty_\omega(C)$
and note that if $(\bigstar)$ and $(\bigstar)^\partial$ hold for elements in $J^\infty(C)$ and 
$M^\infty(C)$ then they also hold for elements of $J^\infty_\omega(C)$ and 
$M^\infty_\omega(C)$ since, e.g., $\bigvee_{i=1}^np_i\nleq u$ if and only if $p_i\nleq u$
for some $i$ with $1\leq i\leq n$. We will now show that each $k\in J^\infty_\omega(C)$
is compact in $C$. Let $U\subseteq({\uparrow}k)^c$ be directed. Then for each $u\in U$ there is 
$m\in M_k$ with $u\leq m$. We claim that in fact there is a single $m\in M_k$ with $U\leq m$.
To see this, suppose that for each $m\in M_k$ there is a $u_m\in U$ with $u_m\nleq m$.
Since $U$ is directed, there is $u\in U$ that is above each element of the finite subset 
$\{u_m\mid m\in M_k\}$ of $U$. But then $u\nleq m$ for each $m\in M_k$ which is a in 
contradiction with our assumptions. Note that this is a general argument showing that if a
directed set is contained in a lower tooth then it is below one of the generators of the tooth.
Now $U\leq m$ implies $\bigvee U\leq m$ so that $\bigvee U\neq k$ and ${\uparrow}k$ is 
compact. Further, as $C$ is perfect, for each $u\in C$
\begin{align*}
u &=\bigvee\{p\in J^\infty(C)\mid p\leq u\}\\
   &=\bigvee\{k\in J_\omega^\infty(C)\mid k\leq u\}
\end{align*}
where the second join is directed and thus $C$ is algebraic. As a consequence the Scott 
topology as well as the lower topology are spectral. We now show that the Scott topology 
is equal to the upper 
topology. It is always the case that the Scott topology contains the upper topology. Let
$U$ be Scott open and let $u\in U$. Then, as  $u$ is the directed join of 
$\{k\in J_\omega^\infty(C)\mid k\leq u\}$, there is a $k\in J_\omega^\infty(C)$ with
$u\geq k\in U$, or equivalently, $u\in {\uparrow}k\subseteq U$. As 
$({\uparrow}k)^c={\downarrow} M_k=\bigcup_{m\in M_k}{\downarrow}m$ we have
${\uparrow}k=\bigcap_{m\in M_k}({\downarrow}m)^c$ which is open in the interval topology
since $M_k$ is finite. Thus $U$ is the union of sets that are open in the interval topology
and we conclude that the two topologies agree. The rest follows by order duality.
\qed\smallskip

We will show that the canonical extension of any lattice lying in a finitely generated variety 
satisfies the hypothesis of the above proposition -- and thus also its conclusion. This shows
that working in lattice expansions based on lattices lying in finitely generated varieties of 
lattices essentially brings about the same advantages as working on distributive lattice 
expansions (for which the underlying lattice lies in the lattice variety generated by the 
two element lattice). As explained above, the strategy in proving this is to show that any 
finite lattice $A$
satisfies the hypothesis of the proposition and then move through the operators $P_B,S, H$.
First note that the canonical extension of any lattice is a perfect lattice so we just need to 
prove that the conditions $(\bigstar)$ and $(\bigstar)^\partial$ hold. Also, it is clear that any 
finite lattice satisfies the conditions. The only detail that may be worth comment is the 
observation that, in any lattice, an element $m$ which is maximal with respect to not 
being greater than or equal to some other element $k$ necessarily must be completely
meet irreducible since $m<a$ implies $k\leq a$.

\begin{lemma}
Let $A$ be a finite lattice, $B\in P_B(A)$. Then $B^\delta$ satisfies the conditions $(\bigstar)$ 
and $(\bigstar)^\partial$.
\end{lemma}

\proof
By Theorem~\ref{BP-lat} we have $B^\delta=A^X$ and it is straight forward to verify that 
$J^\infty(A^X)=\{\pi_x^\flat(p)\mid x\in X \mbox{ and } p\in J(A)\}$ and 
$M^\infty(A^X)=\{\pi_x^\sharp(m)\mid x\in X \mbox{ and } m\in M(A)\}$. The condition 
$(\bigstar)$ clearly holds since, for each $x\in X$ and $p\in J(A)$ the set
$({\uparrow}\pi_x^\flat(p))^c\cap M^\infty(A^X)=\{\pi_x^\sharp(m)\mid p\nleq m\in M(A)\}$  which 
is finite. By order duality $(\bigstar)^\partial$ holds as well.
\qed\smallskip

\begin{lemma}
Let $A$ be a finite lattice, $C\in S(P_B(A))$. Then $C^\delta$ satisfies the conditions $(\bigstar)$ 
and $(\bigstar)^\partial$.
\end{lemma}

\proof
If $C\in S(P_B(A))$ then $C\hookrightarrow B\hookrightarrow A^X$ where the second 
embedding is a Boolean product. Consequently $C^\delta\hookrightarrow B^\delta=A^X$
where the embedding is a complete lattice embedding. That is, we may assume that 
$D:=C^\delta$ is a perfect lattice which is a complete sublattice of $A^X$.  Now let $x\in X$.
Note that $\pi_x(D)$ is a sublattice of the finite lattice $A$. Consider the restricted projection 
$\pi_x{\upharpoonright}D:D\to\pi_x(D)$. It is a complete lattice homomorphism and has right 
and left adjoints. We claim that
\begin{align*}
J^\infty(D)&=\{(\pi_x{\upharpoonright}D)^\flat(p)\mid x\in X,p\in J(\pi_x(D))\}\\
M^\infty(D)&=\{(\pi_x{\upharpoonright}D)^\sharp(m)\mid  x\in X,m\in M(\pi_x(D))\}.
\end{align*}
We first show that $(\pi_x{\upharpoonright}D)^\flat(p)$ is completely join irreducible in $D$ for 
each $x\in X$ and $p\in J(\pi_x(D))$. To this end, let ${\mathcal E}\subseteq D$ with 
$e< (\pi_x{\upharpoonright}D)^\flat(p)$ for each $e\in {\mathcal E}$. Thus at least 
$e_x\leq((\pi_x{\upharpoonright}D)^\flat(p))_x=p$. However, if $e_x=p$ then $p\leq e_x$ and 
thus $(\pi_x{\upharpoonright}D)^\flat(p)\leq e$ which is not the case. So in fact, $e_x<p$ for 
each $e\in\mathcal E$. Since $p\in J(\pi_x(D))=J^\infty(\pi_x(D))$ it follows that 
$(\bigvee{\mathcal E})_x=\bigvee\{e_x\mid e\in{\mathcal E}\}<p$ and thus 
$\bigvee{\mathcal E}\neq (\pi_x{\upharpoonright}D)^\flat(p)$ so that the latter has been 
proved to be completely join irreducible. Since $((\pi_x{\upharpoonright}D)^\flat(p))_x=p$
it is clear that for each $d\in D$ we have
\[
d=\bigvee\{(\pi_x{\upharpoonright}D)^\flat(p)\mid x\in X, p\in J(\pi_x(D)), p\leq d_x\}
\]
so that the $(\pi_x{\upharpoonright}D)^\flat(p)$ must account for all the completely join 
irreducibles in $D$. The statement about completely meet irreducibles follows by order
duality.
Finally, let $x\in X$, $p\in J(\pi_x(D))$, and $d\in D$, then 
\begin{align*}
(\pi_x{\upharpoonright}D)^\flat(p)\nleq d\quad &\iff \quad p\nleq d_x\\
                      &\iff \quad d_x\in {\downarrow}M_p\\
                      &\iff \quad d\in {\downarrow}\{(\pi_x{\upharpoonright}D)^\sharp(m)\mid m\in M_p\}
\end{align*}
where $M_p$ is the set of maximal elements of $({\uparrow}p)^c$ in $\pi_x(D))$. Thus 
$(\bigstar)$ holds and by order duality $(\bigstar)^\partial$ also holds and we have proved 
the lemma.
\qed\smallskip

\begin{lemma}
Let $D$ be a complete lattice satisfying the conditions $(\bigstar)$. Further, let $E$ be a 
complete homomorphic image of $D$. Then $E$ also satisfies $(\bigstar)$. The same 
holds for $(\bigstar)^\partial$.
\end{lemma}

\proof
Let $D$ and $E$ be complete lattices, $h:D\twoheadrightarrow E$ a surjective complete
lattice homomorphism. Further, let $q\in J^\infty(E)$ and $e\in E$ with $q\nleq e$. Since 
$h$ is completely meet preserving it has a lower adjoint $h^\flat:E\to D$ given by
\[
\forall e\in E\ \forall d\in D\quad (h^\flat(e)\leq d\ \iff\ e\leq h(d))
\]
As $h$ is surjective it is not hard to see that $h^\flat$ carries completely join 
irreducible elements to completely join irreducible elements. Thus $h^\flat(q)\in J^\infty(D)$
and it follows by $(\bigstar)$ that $({\uparrow}h^\flat(q))^c={\downarrow}M$ for some finite 
subset $M$ of $M^\infty(D)$. Surjectivity of $h$ also implies that there is $d\in D$ with 
$h(d)=e$ and $q\nleq e=h(d)$ implies $h^\flat(q)\nleq d$ by the adjunction property.
Thus there is an $m\in M$ with $d\leq m$. Since $h$ is order preserving then 
$e=h(d)\leq h(m)$ so that $({\uparrow}q)^c={\downarrow}h(M)$. The set $h(M)$ is finite
and thus each element of $h(M)$ is below a maximal one and we have 
$({\uparrow}q)^c={\downarrow}\max(h(M))$. Since the elements of $\max(h(M))$ are 
also maximal in $({\uparrow}q)^c$ they are necessarily completely meet irreducible. 
The hypotheses are self dual so clearly, the dual condition $(\bigstar)^\partial$ is also
preserved.
\qed\smallskip

\begin{remark}
Let $A$ be a finite lattice and let $n$ be such that 
\[
\forall B\in S(A)\ \forall p\in J(B) \quad |\max(({\uparrow}p)^c)|\leq n
\]
(such an $n$ exists since $A$ is finite and only has finitely many subalgebras)
then
\[
\forall E\in HSP_B(A)={\mathcal V}(A)\ \forall p\in J^\infty(E^\delta) \quad |\max(({\uparrow}p)^c)|\leq n.
\]
This follows easily by looking at the proofs of the three lemmas.

Note also that if we start from any class $\mathcal K$ of finite lattices (not necessarily 
of bounded size) our lemmas still go through, so the algebras in $HSP_B({\mathcal K})$
satisfy $(\bigstar)$ and $(\bigstar)^\partial$ and thus also the conclusion of 
Proposition~\ref{prop:star}. This class is of course not necessarily a variety.
\end{remark}

We reiterate what we have achieved:

\begin{theorem}\label{thm:fgv}
Let $A$ be a finite lattice and let $E\in{\mathcal V}(A)$ then $E^\delta$ is doubly algebraic 
and the Scott and the upper topologies on $E^\delta$ are equal and this topology is spectral. 
Dually, the dual Scott and the lower topologies on $E^\delta$ are equal and this topology is 
spectral as well. The bases of compact-opens of these two topologies come in complementary 
pairs of upper and lower teeth and the join of the two topologies makes $C$ into a Priestley space.
\end{theorem}

Using the above result, we can prove the following result which is closely related to the result 
in \cite{Ha06}. Note though that we do not need the restriction to monotone additional
operations. The connection between canonical extension, profinite completion and topology is 
studied in further detail in \cite{VoPhD}. 
\begin{theorem}
Canonical extension is functorial on any finitely generated variety of lattice expansions and the 
canonical extension of all operations are continuous in the interval(=double Scott) topology. 
This implies that all basic operations on all the algebras in such a variety are smooth and that 
all the canonical extensions are Stone topological algebras in their interval(=double Scott) topologies.
\end{theorem}

\proof
Note first that by the above result combined with Proposition~\ref{prop:univenv}(2) the envelopes
of any maps between lattices lying in finitely generated varieties are universal so that the results
of the previous section may be applied. Our strategy is then to show, at each level of generation 
(through $P_B,S$ and $H$), that the additional operation lifts to an $(\iota,\iota)$-continuous map. 
It then follows by Proposition~\ref{prop1} that homomorphisms lift to the canonical extensions and 
thus that canonical extension is functorial on finitely generated varieties.

Let $A$ be a finite lattice, and let $B\leq A^X$ be a Boolean product. Without loss of generality, we
consider just one basic operation $f:A^n\to A$ on $A$. We know that $B^\delta = A^X$. Also, since 
the interval topology on bounded lattices is productive \cite{AlFr66} and $A$ is finite, the interval 
topology on $A^X$ is simply the product topology for $A$ with the discrete topology. Clearly then 
the map $f^{[X]}$ which is just $f$ coordinate-wise is interval continuous and extends $f^B$ since 
this map is coordinate-wise $f$ as well. By Theorem~\ref{thm:addop}(1), it follows that $f^B$ is 
smooth and that $(f^B)^\delta$ is equal to $f^{[X]}$.

Now let $C$ be in $S(P_B((A,f))$. Then $(C,f^C)\hookrightarrow (B,f^B)\leq (A^X,f^{[X]})$ where the 
latter is a Boolean product and thus $D:=C^\delta$ is a complete sublattice of $B^\delta=A^X$. By  
Theorem~\ref{thm:fgv}, the upper topology, $\iota^{\uparrow}$, on $A^X$ is generated by the 
subbasis consisting of the sets ${\uparrow}\pi_x^\flat(p)$ for $x\in X$ and $p\in J(A)$ whereas the 
upper topology on $D=C^\delta$ is generated by the subbasis consisting of the sets 
${\uparrow}_D(\pi_x{\upharpoonright}D)^\flat(q)$ for $x\in X$ and $q\in J(\pi_x(D))$. Note that for 
$x\in X$ and $p\in J(A)$ we have 
\begin{align*}
{\uparrow}\pi_x^\flat(p)\cap D &={\uparrow}_D(\pi_x{\upharpoonright}D)^\flat(a)\\
                                               &=\bigcap\{{\uparrow}_D(\pi_x{\upharpoonright}D)^\flat(q)\mid a\geq q\in J(\pi_x(D))\}
\end{align*}
where $a=\bigwedge\{a'\in \pi_x(D)\mid p\leq a'\}$. That is, the interval topology on $D$ is the subspace 
topology inherited from $A^X$. Secondly, we show that $(f^C)^\sigma$ must be the coordinate-wise map 
$f^{[X]}{\upharpoonright}D$. Let $(u_1,\ldots,u_n)\in D^n$ and $x\in X$. Then 
$U=\{(v_1,\ldots,v_n)\in D^n\mid (v_i)_x=(u_i)_x\mbox{ for each }i\}$  is open in the interval topology and 
thus in the $\delta$ topology on $D^n$. For any $U'$ open in the $\delta$ topology on $D^n$
with $U'\subseteq U$ we have $\pi_x(f^C(U'\cap C^n))=\{f^{[X]}((u_1)_x,\ldots,(u_n)_x)\}$ since $f^C$ is $f$
coordinate-wise. It follows that lower (and upper) envelope(s) of $f^C$ is the coordinate-wise map $f$.
 Finally putting these two things together we see that $(f^C)^\sigma$ is equal to the restriction of the continuous
 map $f^{[X]}$ to the subspace $D$ of $A^X$ and thus $f^C$ is smooth and $(f^C)^\delta$ is continuous in 
 the interval topology as required. 

To complete the proof, let $(E,f^E)$ be in $H(S(P_B((A,f))))$. Then there is $(C,f^C)\in S(P_B((A,f)))$ and a
surjective homomorphism $h:(C,f^C)\twoheadrightarrow (E,f^E)$. By our proof in the previous paragraph, 
$f^C$ is smooth and in fact $(f^C)^\delta$ is $(\iota,\iota)$-continuous. Thus Proposition~\ref{prop2} allows
us to conclude the same of $f^E$ provided $h^\delta$ can be shown to send $h^\delta$-saturated 
$\iota$-open sets to $\iota$-open sets. To this end,
let $U$ be an $\iota$-open $h^\delta$-saturated subset of $C^\delta$ and let $W=h^\delta(U)$. By 
Theorem~\ref{thm:fgv}, the interval topologies on these lattices are the Scott topologies and thus we just 
need to show that $W$ is inaccessible by directed joins. Let $D$ be directed subset of $E^\delta$ and 
suppose $\bigvee D\in W$. Since $h^\delta:C^\delta\to E^\delta$ is a complete homomorphism, it has a 
lower adjoint $(h^\delta)^\flat:E^\delta\to C^\delta$ which is necessarily join preserving. Thus 
$(h^\delta)^\flat(D)$ is directed in $C^\delta$ and $\bigvee(h^\delta)^\flat(D)=(h^\delta)^\flat(\bigvee D)$.
Furthermore, since $U$ is $h^\delta$-saturated and 
$h^\delta((h^\delta)^\flat(\bigvee D)=\bigvee D\in W=h^\delta(U)$, it follows that $(h^\delta)^\flat(\bigvee D)\in U$.
Now, since $U$ is $\iota$-open and thus Scott open, it follows that there is a $d\in D$ with $(h^\delta)^\flat(d)\in U$.
But then $d\in W$ and we have proved that $W$ is Scott open.
 \qed

\bibliographystyle{amsplain}     

\end{document}